\documentclass[12pt,reqno,psamsfonts]{amsart}
\usepackage[alphabetic]{amsrefs}
\usepackage{feynmp}
\usepackage{latexsym}
\usepackage{amsfonts}
\usepackage{amssymb}
\usepackage{amscd}
\usepackage{amsbsy}
\usepackage{amsmath}
\usepackage{bbm}
\usepackage{euscript}
\usepackage{mathrsfs}
\usepackage{yhmath}
\usepackage[dvips]{graphicx}
\usepackage{epsfig}    
\usepackage[all]{xy}
\xyoption{knot}
\xyoption{graph}
\xyoption{tile}
\xyoption{arc}

\newcommand{\raa}{\rightarrow}
\newcommand{\xbulx}{\bullet}

\newcommand{\ZZ}{{\mathbb Z}}

\newcommand{\CC}{{\mathbb C}}
\renewcommand{\AA}{{\mathbb A}}

\newcommand{\bpartial}{{\bar\partial}}

\newcommand{\cC}{{\mathcal C}}
\newcommand{\obul}{{\Omega^{0,\bullet}}}
\newcommand{\KK}{{\mathsf K}}
\newcommand{\Spec}{{\rm Spec}}

\newcommand{\Z}{\mathcal Y}

\newcommand{\sD}{{\mathsf D}}
\newcommand{\sT}{{\mathsf T}}
\newcommand{\Sym}{{\rm Sym}}

\newcommand{\Dperf}{{\mathfrak D}}

\newcommand{\Vect}{{\rm Vect}}

\newcommand{\frC}{{\mathfrak C}}
\newcommand{\wA}{{\mathcal A}}
\newcommand{\sDb}{{\sD^{\rm b}}}

\def\cf{{\it cf.\ }}
\def\rhs{{\it r.h.s.\ }}

\def\End{\mathop{{\rm End}}\nolimits}

\def\Hom{\mathop{{\rm Hom}}\nolimits}
\def\Ext{\mathop{{\rm Ext}}\nolimits}

\def\deg{ \mathop{{\rm deg}}\nolimits }
\def\p{^{\prime}}

\def\del{ \partial }

\def\End{\mathop{{\rm End}}\nolimits}

\def\lrbc#1{ \left( #1 \right) }
\def\lrbs#1{ \left[ #1 \right] }

\def\corr#1{ \langle #1 \rangle }

\def\xmapta#1{ \mathop{{\longrightarrow}}^{#1} }

\def\xumap#1#2#3{ \begin{CD} #1 @>{#2}>> #3 \end{CD} }

\parskip=\medskipamount                 
\thicklines                     
\def\inbar{\vrule height1.5ex width.4pt depth0pt}
\def\IC{\relax\,\hbox{$\inbar\kern-.3em{\rm C}$}}
\def\IN{\relax{\rm I\kern-.18em N}}
\def\IQ{\relax\,\hbox{$\inbar\kern-.3em{\rm Q}$}}
\def\IR{\relax{\rm I\kern-.18em R}}
\def\ZZ{\relax{\sf Z\kern-.4em Z}}

  \def\cC{{\cal C}} 
\def\cE{{\cal E}}
\def\cF{{\cal F}}   
   \def\cM{{\cal M}}
 \def\cO{{\cal O}} \def\cP{{\cal P}} 
   
\newtheorem{theorem}{Theorem}[section]

\newtheorem{corollary}[theorem]{Corollary}
\newtheorem{conjecture}[theorem]{Conjecture}

\catcode`\@=11

\newif\if@fewtab\@fewtabtrue

\catcode`\@=11

\newif\if@fewtab\@fewtabtrue

{\count255=\time\divide\count255 by 60
\xdef\hourmin{\number\count255} \multiply\count255
by-60\advance\count255 by\time
\xdef\hourmin{\hourmin:\ifnum\count255<10 0\fi\the\count255}}
\def\ps@draft{\let\@mkboth\@gobbletwo
    \def\@oddhead{}
    \def\@oddfoot
      {\hbox to 7 cm{\footnotesize {\em Draft of \jobname:} \draftdate
       \hfil}\hskip -7cm\hfil\rm\thepage \hfil}
    \def\@evenhead{}\let\@evenfoot\@oddfoot}

\def\ceqno{\global\@fewtabfalse
    \ifcase\@eqcnt \def\@tempa{& & &}\or \def\@tempa{& &}
      \or \def\@tempa{&}
      \or\def\@tempa{}\fi\@tempa
{\rm(\theequation)}}

\def\aeqno#1{\global\@fewtabfalse
    \ifcase\@eqcnt \def\@tempa{& & &}\or \def\@tempa{& &}
      \or \def\@tempa{&}
      \or\def\@tempa{}\fi\@tempa
{\rm(\theequation,#1)}}

\def\label#1{\ifnum\draftcontrol=1
 \global\def\draftnote{$\scriptstyle #1$}\fi
 \@bsphack\if@filesw {\let\thepage\relax
   \def\protect{\noexpand\noexpand\noexpand}%
\xdef\@gtempa{\write\@auxout{\string
      \newlabel{#1}{{\@currentlabel}{\thepage}}}}}\@gtempa
   \if@nobreak \ifvmode\nobreak\fi\fi\fi
  \@esphack}

\def\alabel#1#2{\label{#1}\global\@fewtabfalse
    \ifcase\@eqcnt \def\@tempa{& & &}\or \def\@tempa{& &}
      \or \def\@tempa{&}
      \or\def\@tempa{}\fi\@tempa
{\hbox to 3cm{\phantom{\rm(\theequation,#2)} \draftnote
\hfil}\hskip -3cm {\rm(\theequation,#2)}}}

\def\clabel#1{\label{#1}\global\@fewtabfalse
    \ifcase\@eqcnt \def\@tempa{& & &}\or \def\@tempa{& &}
      \or \def\@tempa{&}
      \or\def\@tempa{}\fi\@tempa
{\hbox to 3cm{\phantom{\rm(\theequation)} \draftnote \hfil}\hskip
-3cm{\rm(\theequation)}}}

\def\eqnarray{\def\draftnote{{}}\global\@fewtabtrue
\stepcounter{equation}\let\@currentlabel=\theequation
\global\@eqnswtrue
\global\@eqcnt\z@\tabskip\@centering\let\\=\@eqncr
$$\halign to \displaywidth\bgroup\@eqnsel\hskip\@centering\@eqcnt\z@
  $\displaystyle\tabskip\z@{##}$&\global\@eqcnt\@ne
  \hskip 1\arraycolsep \hfil$\displaystyle{##}$\hfil
  &\global\@eqcnt\tw@ \hskip 1\arraycolsep
$\displaystyle\tabskip\z@{##}$ \hfil
\tabskip\@centering&\global\@eqcnt\thr@@\llap{##}\tabskip\z@ \cr}

\def\endeqnarray{\@@eqncr\egroup
      \global\advance\c@equation\m@ne$$\global\@ignoretrue}

\def\@eqnnum{\hbox to 3cm{\phantom{\rm(\theequation)} \draftnote
                         \hfil}\hskip -3cm {\rm(\theequation)}}

\def\@@eqncr{\let\@tempa\relax
    \ifcase\@eqcnt \def\@tempa{& & &}\or \def\@tempa{& &}
      \or \def\@tempa{&}
      \or\def\@tempa{}
\fi\@tempa \if@eqnsw \if@fewtab\@eqnnum\fi
\stepcounter{equation}\fi\global
\@eqnswtrue\global\@eqcnt\z@\global\@fewtabtrue\cr}

\def\draftcite#1{\ifnum\draftcontrol=1#1\else{}\fi}

\def\@lbibitem[#1]#2{\item{}\hskip -3cm \hbox to 2cm
{\hfil$\scriptstyle\draftcite{#2}$}\hskip
1cm[\@biblabel{#1}]\if@filesw
     {\def\protect##1{\string ##1\space}\immediate
      \write\@auxout{\string\bibcite{#2}{#1}}}\fi\ignorespaces}

\def\@bibitem#1{\item\hskip -3cm \hbox to 2cm
{\hfil $\scriptstyle\draftcite{#1}$}\hskip 1cm \if@filesw
\immediate\write\@auxout
       {\string\bibcite{#1}{\the\value{\@listctr}}}\fi\ignorespaces}

\def\draftdate{\number\month/\number\day/\number\year\ \ \ \hourmin }
 \global\def\draftcontrol{0}
\catcode`\@=12

\def\theequation{{\thesection.\arabic{equation}}}

\def\qq{\begin{eqnarray}}
\def\qqq{\end{eqnarray}}
\def\ee{\begin{eqnarray}}
\def\eee{\end{eqnarray}}
\def\rx#1{~(\ref{#1})}

\def\ex#1{eq.\hspace*{-3pt}\rx{#1}}
\def\eex#1{eqs.\hspace*{-3pt}\rx{#1}}
\def\cx#1{~\cite{#1}}
\def\rw#1{~\ref{#1}}

\def\xlee#1{ \begin{eqnarray} \label{#1} }
\def\xeee{ \end{eqnarray} }
\def\ylee#1{ \begin{eqnarray}\nonumber }
\def\yeee{ \end{eqnarray} }
\def\zlee#1{ \begin{displaymath} }
\def\zleee{ \end{displaymath} }
\def\wlee#1{ $ }
\def\weee{ $ }

\newlength{\shiftwidth}
\def\shift#1{&&\hbox to \shiftwidth{\hfill $\displaystyle#1$}}
\newlength{\sshiftwidth}
\def\sshift#1{\lefteqn{\hbox to
\sshiftwidth{\hfill$\displaystyle#1$}}}

\def\qbezier{\bezier{120}}
\def\DottedCircle{
\bezier{4}(0.966,-0.259)(1.04,0)(0.966,0.259)
\bezier{4}(0.966,0.259)(0.897,0.518)(0.707,0.707)
\bezier{4}(0.707,0.707)(0.518,0.897)(0.259,0.966)
\bezier{4}(0.259,0.966)(0,1.04)(-0.259,0.966)
\bezier{4}(-0.259,0.966)(-0.518,0.897)(-0.707,0.707)
\bezier{4}(-0.707,0.707)(-0.897,0.518)(-0.966,0.259)
\bezier{4}(-0.966,0.259)(-1.04,0)(-0.966,-0.259)
\bezier{4}(-0.966,-0.259)(-0.897,-0.518)(-0.707,-0.707)
\bezier{4}(-0.707,-0.707)(-0.518,-0.897)(-0.259,-0.966)
\bezier{4}(-0.259,-0.966)(0,-1.04)(0.259,-0.966)
\bezier{4}(0.259,-0.966)(0.518,-0.897)(0.707,-0.707)
\bezier{4}(0.707,-0.707)(0.897,-0.518)(0.966,-0.259) }
\def\Endpoint[#1]{
\ifcase#1 \put(1,0){\circle*{0.15}}
\or\put(0.866,0.5){\circle*{0.15}}
\or\put(0.5,0.866){\circle*{0.15}} \or\put(0,1){\circle*{0.15}}
\or\put(-0.5,0.866){\circle*{0.15}}
\or\put(-0.866,0.5){\circle*{0.15}} \or\put(-1,0){\circle*{0.15}}
\or\put(-0.866,-0.5){\circle*{0.15}}
\or\put(-0.5,-0.866){\circle*{0.15}} \or\put(0,-1){\circle*{0.15}}
\or\put(0.5,-0.866){\circle*{0.15}}
\or\put(0.866,-0.5){\circle*{0.15}} \fi}
\def\Arc[#1]{
\thicklines         
\ifcase#1 \bezier{25}(0.966,-0.259)(1.04,0)(0.966,0.259) \or
\bezier{25}(0.966,0.259)(0.897,0.518)(0.707,0.707) \or
\bezier{25}(0.707,0.707)(0.518,0.897)(0.259,0.966) \or
\bezier{25}(0.259,0.966)(0,1.04)(-0.259,0.966) \or
\bezier{25}(-0.259,0.966)(-0.518,0.897)(-0.707,0.707) \or
\bezier{25}(-0.707,0.707)(-0.897,0.518)(-0.966,0.259) \or
\bezier{25}(-0.966,0.259)(-1.04,0)(-0.966,-0.259) \or
\bezier{25}(-0.966,-0.259)(-0.897,-0.518)(-0.707,-0.707) \or
\bezier{25}(-0.707,-0.707)(-0.518,-0.897)(-0.259,-0.966) \or
\bezier{25}(-0.259,-0.966)(0,-1.04)(0.259,-0.966) \or
\bezier{25}(0.259,-0.966)(0.518,-0.897)(0.707,-0.707) \or
\bezier{25}(0.707,-0.707)(0.897,-0.518)(0.966,-0.259) \fi}
\def\DottedArc[#1]{
\ifcase#1 \bezier{4}(0.966,-0.259)(1.04,0)(0.966,0.259) \or
\bezier{4}(0.966,0.259)(0.897,0.518)(0.707,0.707) \or
\bezier{4}(0.707,0.707)(0.518,0.897)(0.259,0.966) \or
\bezier{4}(0.259,0.966)(0,1.04)(-0.259,0.966) \or
\bezier{4}(-0.259,0.966)(-0.518,0.897)(-0.707,0.707) \or
\bezier{4}(-0.707,0.707)(-0.897,0.518)(-0.966,0.259) \or
\bezier{4}(-0.966,0.259)(-1.04,0)(-0.966,-0.259) \or
\bezier{4}(-0.966,-0.259)(-0.897,-0.518)(-0.707,-0.707) \or
\bezier{4}(-0.707,-0.707)(-0.518,-0.897)(-0.259,-0.966) \or
\bezier{4}(-0.259,-0.966)(0,-1.04)(0.259,-0.966) \or
\bezier{4}(0.259,-0.966)(0.518,-0.897)(0.707,-0.707) \or
\bezier{4}(0.707,-0.707)(0.897,-0.518)(0.966,-0.259) \fi}
\def\Chord[#1,#2]{
\thinlines \ifnum#1>#2\Chord[#2,#1] \else\ifnum#1<#2 \ifcase#1
\ifcase#2 \or\qbezier(1,0)(0.516,0.138)(0.866,0.5)
\or\qbezier(1,0)(0.45,0.26)(0.5,0.866)
\or\qbezier(1,0)(0.327,0.327)(0,1)
\or\qbezier(1,0)(0.179,0.311)(-0.5,0.866)
\or\qbezier(1,0)(0.0536,0.2)(-0.866,0.5) \or\put(1, 0){\line(-2,
0){2}} \or\qbezier(1,0)(0.0536,-0.2)(-0.866,-0.5)
\or\qbezier(1,0)(0.179,-0.311)(-0.5,-0.866)
\or\qbezier(1,0)(0.327,-0.327)(0,-1)
\or\qbezier(1,0)(0.45,-0.26)(0.5,-0.866)
\or\qbezier(1,0)(0.516,-0.138)(0.866,-0.5) \fi \or\ifcase#2\or
\or\qbezier(0.866,0.5)(0.378,0.378)(0.5,0.866)
\or\qbezier(0.866,0.5)(0.26,0.45)(0,1)
\or\qbezier(0.866,0.5)(0.12,0.446)(-0.5,0.866)
\or\qbezier(0.866,0.5)(0,0.359)(-0.866,0.5)
\or\qbezier(0.866,0.5)(-0.0536,0.2)(-1,0) \or\put(0.866,
0.5){\line(-5, -3){1.73}}
\or\qbezier(0.866,0.5)(0.146,-0.146)(-0.5,-0.866)
\or\qbezier(0.866,0.5)(0.311,-0.179)(0,-1)
\or\qbezier(0.866,0.5)(0.446,-0.12)(0.5,-0.866)
\or\qbezier(0.866,0.5)(0.52,0)(0.866,-0.5) \fi \or\ifcase#2\or\or
\or\qbezier(0.5,0.866)(0.138,0.516)(0,1)
\or\qbezier(0.5,0.866)(0,0.52)(-0.5,0.866)
\or\qbezier(0.5,0.866)(-0.12,0.446)(-0.866,0.5)
\or\qbezier(0.5,0.866)(-0.179,0.311)(-1,0)
\or\qbezier(0.5,0.866)(-0.146,0.146)(-0.866,-0.5) \or\put(0.5,
0.866){\line(-3, -5){1}} \or\qbezier(0.5,0.866)(0.2,-0.0536)(0,-1)
\or\qbezier(0.5,0.866)(0.359,0)(0.5,-0.866)
\or\qbezier(0.5,0.866)(0.446,0.12)(0.866,-0.5) \fi
\or\ifcase#2\or\or\or \or\qbezier(0,1.)(-0.138,0.516)(-0.5,0.866)
\or\qbezier(0,1.)(-0.26,0.45)(-0.866,0.5)
\or\qbezier(0,1.)(-0.327,0.327)(-1,0)
\or\qbezier(0,1.)(-0.311,0.179)(-0.866,-0.5)
\or\qbezier(0,1.)(-0.2,0.0536)(-0.5,-0.866) \or\put(0, 1){\line(0,
-2){2}} \or\qbezier(0,1.)(0.2,0.0536)(0.5,-0.866)
\or\qbezier(0,1.)(0.311,0.179)(0.866,-0.5) \fi
\or\ifcase#2\or\or\or\or
\or\qbezier(-0.5,0.866)(-0.378,0.378)(-0.866,0.5)
\or\qbezier(-0.5,0.866)(-0.45,0.26)(-1,0)
\or\qbezier(-0.5,0.866)(-0.446,0.12)(-0.866,-0.5)
\or\qbezier(-0.5,0.866)(-0.359,0)(-0.5,-0.866)
\or\qbezier(-0.5,0.866)(-0.2,-0.0536)(0,-1) \or\put(-0.5,
0.866){\line(3, -5){1}}
\or\qbezier(-0.5,0.866)(0.146,0.146)(0.866,-0.5) \fi
\or\ifcase#2\or\or\or\or\or
\or\qbezier(-0.866,0.5)(-0.516,0.138)(-1,0)
\or\qbezier(-0.866,0.5)(-0.52,0)(-0.866,-0.5)
\or\qbezier(-0.866,0.5)(-0.446,-0.12)(-0.5,-0.866)
\or\qbezier(-0.866,0.5)(-0.311,-0.179)(0,-1)
\or\qbezier(-0.866,0.5)(-0.146,-0.146)(0.5,-0.866) \or\put(-0.866,
0.5){\line(5, -3){1.73}} \fi \or\ifcase#2\or\or\or\or\or\or
\or\qbezier(-1,0)(-0.516,-0.138)(-0.866,-0.5)
\or\qbezier(-1,0)(-0.45,-0.26)(-0.5,-0.866)
\or\qbezier(-1,0)(-0.327,-0.327)(0,-1)
\or\qbezier(-1,0)(-0.179,-0.311)(0.5,-0.866)
\or\qbezier(-1,0)(-0.0536,-0.2)(0.866,-0.5) \fi
\or\ifcase#2\or\or\or\or\or\or\or
\or\qbezier(-0.866,-0.5)(-0.378,-0.378)(-0.5,-0.866)
\or\qbezier(-0.866,-0.5)(-0.26,-0.45)(0,-1)
\or\qbezier(-0.866,-0.5)(-0.12,-0.446)(0.5,-0.866)
\or\qbezier(-0.866,-0.5)(0,-0.359)(0.866,-0.5) \fi
\or\ifcase#2\or\or\or\or\or\or\or\or
\or\qbezier(-0.5,-0.866)(-0.138,-0.516)(0,-1)
\or\qbezier(-0.5,-0.866)(0,-0.52)(0.5,-0.866)
\or\qbezier(-0.5,-0.866)(0.12,-0.446)(0.866,-0.5) \fi
\or\ifcase#2\or\or\or\or\or\or\or\or\or
\or\qbezier(0,-1.)(0.138,-0.516)(0.5,-0.866)
\or\qbezier(0,-1.)(0.26,-0.45)(0.866,-0.5) \fi
\or\ifcase#2\or\or\or\or\or\or\or\or\or\or
\or\qbezier(0.5,-0.866)(0.378,-0.378)(0.866,-0.5) \fi\fi\fi\fi}
\def\FullChord[#1,#2]{
\Endpoint[#1] \Endpoint[#2] \Arc[#1] \Arc[#2] \Chord[#1,#2] }
\def\EndChord[#1,#2]{
\Endpoint[#1] \Endpoint[#2] \Chord[#1,#2] }
\def\Picture#1{
\begin{picture}(2,1)(-1,-0.167)
#1
\end{picture}
}
\def\DottedChordDiagram[#1,#2]{
\Picture{\DottedCircle \FullChord[#1,#2]} }
\def\ZZ{ \mathbb{Z} }
\def\IQ{ \mathbb{Q} }
\def\IC{ \mathbb{C} }
\def\IR{ \mathbb{R} }

\def\bfx{ \mathbf{x} }

\def\cC{ \mathcal{C} }

\def\cE{ \mathcal{E} }
\def\cF{ \mathcal{F} }

\def\hlf{ {1\over 2} }

\def\xrl{ \bar{\xl} }

\def\xE{ E }

\def\yY{ Y }
\def\yZ{ Z }

\def\xee{\begin{equation}}
\def\yee{\end{equation}}

\def\xom{ \omega }

\def\Tng{ \mathrm{T} }

\def\aee{ $$ }
\def\aeee{ $$ }

\def\xlabel#1{ \nonumber }

\def\xId{ \mathbbm{1} }

\def\rnk{\mathop{\mathrm{rank}} }

\def\crvng{curving}
\def\crvngp{\crvng\ polynomial}

\def\ddf{divided difference}
\def\ddfs{\ddf s}

\def\ctL{ \mathsf{L} }
\def\ctLL{ \ddot{\ctL} }
\def\ctLLL{ \dddot{\ctL} }
\def\ctLLv#1{ \ctLL #1 }
\def\ctLLbv#1{ \ctLLv{(#1)} }
\def\ctLLXsom{ \ctLLv{\Xsom} }
\def\ctLLX{ \ctLLv{(\xX)} }
\def\ctLLTsU{ \ctLLbv{\TsxcA} }
\def\ctLLTsCn{ \ctLLbv{\Ts\IC^n} }

\def\ctLLTUSk{ \ctLLv{ \mtdfTUSk } }

\def\qhlm{quasi-holomorphic}
\def\qhlmvb{\qhlm\ vector bundle}
\def\cqhlmvb{curved \qhlmvb}

\def\ctLLcrv#1#2{ \ctLL\brB{ (\xXv{#1},-\somv{#1})\times(\xXv{#2},\somv{#2}) } }
\def\ctLLcrot{ \ctLLcrv{1}{2} }

\def\ctLLXsomv#1{ \ctLL\Xsomv{#1} }
\def\ctLLXsomo{ \ctLLXsomv{1} }
\def\ctLLXsomt{ \ctLLXsomv{2} }

\def\xwd{well-defined}

\def\rmA{ \mathrm{A} }
\def\Ainf{ \rmA_\infty }
\def\tAinf{$\Ainf$}
\def\Ainfalg{\tAinf-algebra}
\def\fDAinfalg{filtered Dolbeault \Ainfalg}
\def\FDAi{FD\tAinf}

\def\Khl{ K\"{a}hler }
\def\dR{de Rham}

\def\xpBd{`pseudo-Beltrami' differential}

\def\xwk{weak}
\def\xacrv{curved}

\def\xrca{relative curved differential graded polynomial algebra}
\def\tcrca{2-category of \xrca s}
\def\yCDG{CDG}

\def\smcl{semi-classical}
\def\smclgrd{\smcl\ grading}
\def\smcldgr{\smcl\ degree}

\def\xrl{relative}
\def\xrlb{relatively balanced}
\def\xrlnn{relatively non-negative}

\def\Dlb{Dolbeault}
\def\Dlbf{\Dlb\ filtration}

\def\Knr{Knorrer}
\def\Knrp{\Knr\ periodicity}

\def\NSb{Schouten–-Nijenhuis bracket}

\def\CM{Maurer-Cartan}
\def\CMe{\CM\ equation}
\def\CMlt{\CM\ element}

\def\Ldr{Legendre}
\def\Ldrt{\Ldr\ transform}

\def\ooc{$(1,1)$-curvature}

\def\xPoi{Poisson}
\def\xGrs{Gerstenhaber}
\def\xLie{Lie}
\def\DLa{differential \xLie\ algebra}
\def\DPa{differential \xPoi\ algebra}
\def\DGa{differential \xGrs\ algebra}
\def\PGa{\xPoi-\xGrs\ algebra}
\def\PGa{\xLie-\xGrs\ algebra}
\def\DPGa{differential \PGa}

\def\rlcn{relative connection}
\def\darc{\DLa\ of \rlcn s}

\def\pgbvv#1#2{ [ #1,#2 ]_{\rm LG } }
\def\pgbdd{ \pgbvv{\cdot}{\cdot} }

\def\pgD{ D }
\def\pgDv#1{ \pgD_{#1} }
\def\pgDa{ \pgDv{\pga} }

\def\pga{ \alpha }
\def\pgae{ \dfev{\pga} }
\def\pgav#1{ \pga_{|#1} }
\def\pgai{ \pgav{i} }
\def\pgao{ \pgav{1} }
\def\pgat{ \pgav{2} }
\def\pgb{ \phi }
\def\pgdg{ d_{[\;]} }

\def\zdpa{ \mathscr{P} }
\def\zdpaXsom{ \zdpa\Xsom }
\def\zdpaTsU{ \zdpa(\TsU) }
\def\zdpaT{ \zdpa_{\Ts} }
\def\zdpaTv#1{ \zdpaT(#1) }
\def\zdpaTU{ \zdpaTv{\xcA} }

\def\zdga{ \mathscr{G} }
\def\zdgav#1{ \zdga(#1) }
\def\zdgaU{ \zdgav{\xcA} }
\def\zdgalv#1{ \zdga_{#1} }
\def\zdgalW{ \zdgalv{\xcW} }
\def\zdgalWU{ \zdgalW(\xcA) }

\def\zdrc{ \mathscr{C} }
\def\zdrcEnb{ \zdrc\adgmE }

\def\cdlv#1{ [#1] }
\def\cdlv#1{ \check{#1} }
\def\cpgao{ \cdlv{\pga}_{|1} }

\def\tKmf{Koszul matrix factorization}

\def\xtran{2-translation}
\def\ytran{2-translation}

\def\Mlc{Micro-local}
\def\mlc{micro-local}
\def\mlcz{micro-localization}
\def\prshf{presheaf}
\def\mlcs{\mlc\ \prshf}

\def\zctng{cotangent}

\def\zctngb{\zctng\ bundle}
\def\zdctngb{deformed \zctngb}

\def\zaffn{affine}

\def\zsaff{symplectic \zaffn}

\def\sd{self-dual}

\def\rb{ \mathrm{b} }
\def\rDb{ \xrDs^\rb }
\def\prftpv#1{ #1_{\Zt} }
\def\prftapv#1{ #1^{\mathrm{a}}_{\Zt} }
\def\rDprf{ \prftpv{\xrDs} } 
\def\rDprfa{ \prftapv{\xrDs} }

\def\rDprfbv#1{ \rDprf #1  }
\def\rDprfbA{ \rDprfbv{\ycdgaA} }

\def\rDprfcA{ \rDprf(\xdlA) }
\def\rDprfcAdv#1{ \rDprfcA\ctrnv{#1} }
\def\rDprfcAdn{ \rDprfcAdv{n} }

\def\rDprfcAAbx{ \rDprf(\xdlA\otimes\xdlAbx) }

\def\rDprfU{ \rDprf(\xcA) }

\def\xrDs{ \mathsf{D} }
\def\xrDsd{ \ddot{\xrDs} }
\def\xrDsdd{ \dddot{\xrDs} }
\def\rDsprf{ \prftpv{\xrDs} }
\def\rDDprf{ \prftpv{\xrDsd} }
\def\rDDDprf{ \prftpv{\xrDsdd} }
\def\rDDprfv#1{ \rDDprf(#1) }
\def\rDDprfU{ \rDDprf(\xcA) }
\def\rDDprfUxVy{ \ctLLv(\UxVy) }
\def\rDDprfaUxVy{ \ctLLv(\UxVy) }
\def\rDDprfaUxVypp{ \ctLLv(\UxVypp) }
\def\rDDprfaUpxVy{ \ctLLv(\UpxVy) }
\def\rDDprfaVmyUx{ \ctLLv(\VmyUx) }
\def\rDDprfaVyUmx{ \ctLLv(\VyUmx) }

\def\rDDprfUSk{ \rDDprf(\xcA,\hdf) }
\def\rDDprfYoSk{ \rDDprf(\yYo,\hdf) }
\def\rDDprfUSek{ \rDDprf(\xcA,\dfe\hdf) }

\def\rDDprfUhae{ \rDDprf(\xcA,\ehdf) }
\def\rDDprfUp{ \rDDprf(\xcAp) }

\def\rDprfY{ \rDprf(\yY) }

\def\rDsprfa{ \prftapv{\xrDs} }
\def\rDDprfa{ \prftapv{\xrDsd} }
\def\rDDDprfa{ \prftapv{\xrDsdd} }
\def\rDDprfav#1{ \rDDprfa(#1) }
\def\rDDprfaU{ \rDDprfa(\xcA) }

\def\rDDprfaUp{ \rDDprfa(\xcAp) }

\def\rDprfvv#1#2{ \rDprf(#1,#2) }

\def\rDprfUW{ \rDprfvv{\xcA}{\xcW} }
\def\rDprfUWv#1{ \rDprfvv{\xcAv{#1}}{\xcWv{#1}} }
\def\rDprfUWo{ \rDprfUWv{1} }
\def\rDprfUWt{ \rDprfUWv{2} }

\def\rDprfcUWv#1{ \rDprfvv{\xcUv{#1}}{\xcWv{#1}} }
\def\rDprfcUWo{ \rDprfcUWv{1} }
\def\rDprfcUWt{ \rDprfcUWv{2} }

\def\rDsprfv#1{ \rDsprf(#1) }
\def\rDsprfU{ \rDsprfv{\xcA} }

\def\rDsprfYot{ \rDsprfv{\yYot} }
\def\rDsprfcYoct{ \rDsprfv{\ycYoct} }

\def\rDsprfvv#1#2{ \rDsprf(#1,#2) }
\def\rDsprfsvv#1#2{ \rDsprf(#1;#2) }
\def\rDsprfavv#1#2{ \rDsprfa(#1,#2) }
\def\rDsprfBvv#1#2{ \rDsprf\Big(#1,#2\Big) }

\def\rDsprfUW{ \rDsprfvv{\xcA}{\xcW} }
\def\rDsprfUWtmo{ \rDsprfvv{\xcA}{\xcWt-\xcWo} }

\def\rDsprfUd{ \rDsprfvv{\xcA}{\xdfm} }
\def\rDsprfUdot{ \rDsprfsvv{\xcA}{\xdfmot} }

\def\rDsprfYadoct{ \rDsprfsvv{\yYoct}{\adfmot} }
\def\rDsprfYoxdoct{ \rDsprfsvv{\yYo}{\xdfmot} }

\def\rDsprfUeWot{ \rDsprfvv{\xcAekot}{\xcWot} }
\def\rDsprfUeWoh{ \rDsprfvv{\xcAekoh}{\xcWoh} }
\def\rDsprfUeWth{ \rDsprfvv{\xcAekth}{\xcWth} }

\def\rDsprfUmWv#1{ \rDsprfvv{\dfmUmv{#1}}{\xcWve{#1}} }
\def\rDsprfUmWot{ \rDsprfUmWv{12} }

\def\rDsprfUWav#1{ \rDsprfvv{\xcA}{\xcWv{#1}} }
\def\rDsprfUWao{ \rDsprfUWav{1} }
\def\rDsprfUWat{ \rDsprfUWav{2} }
\def\rDsprfaUW{ \rDsprfavv{\xcA}{\xcW} }
\def\rDsprfaUWot{ \rDsprfavv{\xcA}{\xcWot} }

\def\rDsprfUWot{ \rDsprfUWav{12} }

\def\rDsprfmvv#1#2#3{ \rDsprfvv{#1}{#2}\ytrnv{#3} }

\def\rDsprfmUWL{ \rDsprfmvv{\xcA}{\xcW}{\vcL} }

\def\Fk{Fukaya}
\def\FkF{\Fk-Floer}
\def\FkFc{\FkF\ category}
\def\rFuk{ \mathrm{F} }

\def\TQFT{TQFT}

\def\pb{pull-back}
\def\pf{push-forward}
\def\ccm{curved complex manifold}

\def\lgr{lagrangian}
\def\hlgr{\lgr}
\def\hlgrsm{\hlgr\ submanifold}
\def\hlsm{holomorphic symplectic}

\def\hlsmm{\hlsm\ manifold}
\def\hlsms{\hlsm\ structure}
\def\fCY{Calabi-Yau}

\def\tpd{2-periodic}

\def\tpdpr{\tpd\ \xper}
\def\tpdprc{\tpdpr\ category}

\def\hlfb{holomorphic fibration}
\def\fbl{fiber bundle}

\def\xah{$(0,\bullet)$}
\def\xahf{\xah-form}
\def\xahfs{\xahf s}

\def\xau{augment}
\def\xaug{augmented}
\def\xAug{Augmented}

\def\pdcat{perfect derived category}
\def\Ztpdcat{\tpd\ \pdcat}

\def\cdga{curved differential graded algebra}
\def\CDGA{CDGA}

\def\Ztdgm{\Ztgrdd\ differential  module}
\def\DGM{DGM}
\def\DGM{GDM}
\def\ZtDGM{$\Zt$-\DGM}

\def\corr{correspondence}
\def\corrs{\corr s}

\def\xxtf{2-functor}
\def\cfun{correspondence \xxtf}
\def\Ltf{Legendre \xxtf}

\def\xper{perfect}
\def\xprm{\xper\ module}
\def\xprmi{\xprm}

\def\Ac{Atiyah class}
\def\cAc{curved \Ac}

\def\xad{admissible}

\def\FM{Fourier-Mukai}
\def\FMtr{\FM\ transform}

\def\opfib{1-point fibration}
\def\copfib{curved \opfib}

\def\xep{embedding-projection}
\def\xeps{\xep s}

\def\cfb{curved complex manifold}
\def\cfbs{\cfb s}
\def\xcfb{curved fibration}
\def\xcfbs{\xcfb s}

\def\psmon{pseudo-monoidal}

\def\xgd{clean}
\def\gdint{\xgd\ intersection}
\def\ygd{\xgd}

\def\zord{order}

\def\xanl{analytic }

\def\zN{ d }
\def\zn{ n }

\def\xmf{matrix factorization}
\def\xmfs{\xmf s}
\def\cxmf{category of \xmf}

\def\tZtgdcs{\tpd\ category of a \cfb}

\def\tgsp{target space}
\def\sgmd{sigma-model}

\def\tsgmd{TSM}

\def\crvng{curving}

\def\xX{ X }
\def\xXv#1{ \xX_{#1}}
\def\xXo{ \xXv{1} }
\def\xXt{ \xXv{2} }
\def\xXh{ \xXv{3} }

\def\xXp{ \xX\p }

\def\Tng{ \mathrm{T} }
\def\Ts{ \Tvs }
\def\bTng{ \bar{\Tng} }
\def\bTs{ \dulv{\bTng} }
\def\Tsv#1{ \Ts #1 }

\def\TsU{ \Tsv{\xcA} }

\def\xcU{ \mathcal{U} }
\def\xcUp{ \xcU\p }
\def\xcUv#1{ \xcU_{#1} }
\def\xcUo{ \xcUv{1} }
\def\xcUt{ \xcUv{2} }
\def\xcUh{ \xcUv{3} }
\def\xcUot{ \xcUv{12} }
\def\xcUth{ \xcUv{23} }

\def\TsxcA{ \Ts\xcA }

\def\Nrm{ \mathrm{N} }

\def\Vfb{ V }
\def\Vfbv#1{ \Vfb_{#1} }

\def\xcA{ U }
\def\xcAp{ \xcA\p }
\def\xcAv#1{ \xcA_{#1} }
\def\xcAz{ \xcAv{0} }
\def\xcAo{ \xcAv{1} }
\def\xcAt{ \xcAv{2} }
\def\xcAh{ \xcAv{3} }

\def\xcV{ V }
\def\xcVp{ \xcV\p }
\def\xcVv#1{ \xcV_{#1} }

\def\xcrtv#1#2{ #1_{#2} }
\def\xcAWz{ \xcrtv{\xcA}{\xcWz} }
\def\xcAW{ \xcrtv{\xcA}{\xcW} }
\def\xcAWot{ \xcrtv{\xcA}{\xcWot} }

\def\xcW{ W }
\def\xtt{\mathrm{tot}}
\def\xcWtt{\xcW^{\xtt} }

\def\xcWttv#1{ \xcWtt_{#1} }
\def\xcWttt{ \xcWttv{2} }
\def\xcWttot{ \xcWttv{12} }
\def\xcWp{ \xcW\p }
\def\xcWv#1{ \xcW_{#1} }
\def\xcWo{ \xcWv{1} }
\def\xcWi{ \xcWv{i} }
\def\xcWot{ \xcWv{12} }
\def\xcWoh{ \xcWv{13} }
\def\xcWth{ \xcWv{23} }

\def\xcWt{ \xcWv{2} }
\def\xcWh{ \xcWv{3} }
\def\xcWot{ \xcWv{12} }
\def\xcWth{ \xcWv{23} }
\def\xcWz{ \xcW_{0} }
\def\xcWpl{ \xcW_{+} }
\def\xE{ E }

\def\xEv#1{ \xE^{#1} }
\def\xEdv#1{ \xE_{#1} }
\def\xEo{ \xEdv{1} }
\def\xEt{ \xEdv{2} }
\def\xEh{ \xEdv{3} }
\def\xEi{ \xEdv{i} }
\def\xEot{ \xEdv{12} }

\def\xEth{ \xEdv{23} }

\def\xEoh{ \xEdv{13} }

\def\xEp{ \xE\p }

\def\scO{ \mathscr{O} }
\def\strs#1{ \scO_{#1} }
\def\strsU{ \strs{\xcA} }
\def\strsY{ \strs{\yY} }

\def\xcp{ p }
\def\xcpv#1{ \xcp_{#1} }

\def\xcpot{ \xcpv{12} }

\def\xcE{ \cE }
\def\xcEv#1{ \xcE_{#1} }
\def\xcEo{ \xcEv{1} }
\def\xcEt{ \xcEv{2} }
\def\xcEot{ \xcEv{12} }
\def\xcEotp{ \xcEot\p }
\def\xcEth{ \xcEv{23} }

\def\xbOm{ {\Omega} }
\def\xbOmxvv#1#2{ \xbOm^{#1,#2} }

\def\xbOmxvb#1{ \xbOmxvv{#1}{\bullet} }
\def\xbOmxve#1{ \xbOmxvv{#1}{\yev} }

\def\xbOmxvb#1{ \xbOmxvv{#1}{\bullet} }
\def\xbOmoe{ \xbOmxve{1} }

\def\xbOmbz{ \xbOmxvb{0} }
\def\xbOmv#1{ \xbOm^{0,#1} }
\def\xbOmz{ \xbOmv{0} }
\def\xbOmb{ \xbOmv{\bullet} }
\def\xbOme{ \xbOmv{\yev} }
\def\xbOmev#1{ \xbOme(#1) }
\def\xbOmeU{ \xbOmev{\xcA} }

\def\xbOmecU{ \xbOmev{\xcU} }
\def\xbOmecUp{ \xbOmev{\xcUp} }

\def\xbOmov#1{ \xbOmo(#1) }
\def\xbOmbv#1{ \xbOmb(#1) }
\def\xbOmBv#1{ \xbOmb\big(#1\big) }
\def\xbOmbE{ \xbOmbv{\xE} }

\def\xbOmbzeE{ \xbOmBv{\End\xE } }

\def\xbOmbX{ \xbOmbv{\xX} }
\def\xbOmbU{ \xbOmbv{\xcA} }
\def\xbOmbUp{ \xbOmbv{\xcAp} }

\def\xbOmbU{ \xbOmbv{\xcA} }

\def\xbOmvv#1#2{ \xbOmv{#1}\lrbc{#2} }

\def\xbOmov#1{ \xbOmvv{1}{#1} }
\def\xbOmoU{ \xbOmov{\xcA} }

\def\xbOmbUUST{ \xbOmbv{\xcA,\SbTU} }

\def\xbOmbX{ \xbOmbv{\xX} }
\def\xbOmbUST{ \xbOmbv{\TsU} }

\def\xbOmbUUWT{ \xbOmbv{\xcA,\WbTU} }

\def\xbOmxob{ \xbOmxvb{1} }

\def\xbOmoeE{ \xbOmoe\brb{\End \xE} }

\def\xrDv#1{ \xrDs^{\mathrm{b}}(#1) }
\def\xrDX{ \xrDv{\xX} }

\def\cC{ \mathcal{C} }

\def\mdfv#1{ \mathcal{#1} }
\def\cE{ \mdfv{E} }

\def\cEv#1{ \cE_{#1} }
\def\cEo{ \cEv{1} }
\def\cEt{ \cEv{2} }

\def\cEot{ \cEv{12} }
\def\cEth{ \cEv{23} }

\def\cM{ \mdfv{M} }
\def\cP{ \mdfv{P} }
\def\cPp{ \cP\p }
\def\cPv#1{ \cP_{#1} }
\def\cPo{ \cPv{1} }
\def\cPt{ \cPv{2} }

\def\cM{ \mdfv{M} }
\def\cMp{ \cM\p }
\def\cMv#1{ \cM_{#1} }
\def\cMo{ \cMv{1} }
\def\cMt{ \cMv{2} }

\def\cfM{ \mdfv{M} }

\def\cfMv#1{ \cfM_{#1} }

\def\cfMot{ \cfMv{12} }
\def\cfMth{ \cfMv{23} }

\def\cfMyot{ \cfMv{\ayo,\ayt} }

\def\cP{ \mdfv{P} }
\def\cPp{ \cP\p }

\def\cWv#1{ \cW^{#1} }
\def\cWv#1{ \xvprvv{#1}{\xcW} }
\def\cWvv#1#2{ \cWv{\,#1}_{\,#2} }
\def\cWvv#1#2{ \xvprvv{#1}{\xcWv{#2} } }
\def\cWx{ \cWv{\bax} }
\def\cWy{ \cWv{\bay} }

\def\cWu{ \cWv{\bau} }
\def\cWe{ \cWv{\varnothing} }

\def\cWxo{ \cWvv{\bax}{1} }
\def\cWyt{ \cWvv{\bay}{2} }
\def\cWzh{ \cWvv{\baz}{3} }

\def\cWyo{ \cWvv{\bay}{1} }
\def\cWzt{ \cWvv{\baz}{2} }
\def\cWuh{ \cWvv{\bau}{3} }
\def\cWzot{ \cWvv{\baz}{12} }

\def\nWvv#1#2{ \nvprvv{#1}{\xcW} }
\def\nWzot{ \nWvv{\baz}{12} }

\def\cWuo{ \cWvv{\bau}{1} }
\def\cWwt{ \cWvv{\baw}{2} }

\def\Dlv#1{ \Delta_{#1} }
\def\DlX{ \Dlv{\xX} }

\def\xiosv#1{ \iota_{#1} }
\def\xiosot{ \xiosv{12} }
\def\xiosth{ \xiosv{23} }
\def\xiosoh{ \xiosv{13} }

\def\Mpsvv#1#2{ #1\rightarrow #2 }
\def\MpsMX{ \Mpsvv{M}{\xX} }

\def\xspt{space-time}
\def\xwsh{world-sheet}
\def\xtsp{target-space}

\def\rmT{ \mathrm{T} }

\def\rmH{ \mathrm{H} }

\def\xzth#1{ \hat{#1} }
\def\yev{ \xzth{0} }
\def\yod{ \xzth{1} }

\def\dlb{ \bar{\del} }

\def\dlbv#1{ \dlb_{#1} }

\def\ycdgav#1{ \big( #1, \nbb, \xcW \big) }
\def\ycdgaBv#1{ \Big( #1, \dlb, \xcW \Big) }
\def\ycdgaBU{ \ycdgaBv{\xbOmbU} }

\def\ycdgaA{ \ycdgav{\xdlA} }
\def\ycdgaAv#1{ \big(\xdlA,\nbb,\xcW_{#1} \big) }
\def\ycdgaAo{ \ycdgaAv{1} }
\def\ycdgaAt{ \ycdgaAv{2} }
\def\ycdgaAot{ \big(\xdlA,\nbb,\xcWo+\xcWt \big) }
\def\ycdgaAd{ \big(\xdlA,\nbb,-\xcW\big) }

\def\ycdgaBUv#1{ \Big(\xbOmbU,\dlb,#1 \Big) }
\def\ycdgaBUl{ \ycdgaBUv{\xdfm} }

\def\otdlA{ \otimes_{\xdlA} }
\def\otCbx{ \otimes_{\ICbx} }
\def\otCby{ \otimes_{\ICby} }
\def\otCbxz{ \otimes_{\ICv{\bax,\baz}} }

\def\nbbA{ \nbbv{\xdlA} }

\def\maug{\mathrm{aug}}
\def\maugv#1{ #1^{\maug} }

\def\xchk#1{ \check{#1} }

\def\chmP{ \xchk{\xmP} }

\def\xcMF{ \mathsf{MF} }
\def\xcMFd{ \ddot{\xcMF} }
\def\xcMFdd{ \dddot{\xcMF} }
\def\xcMFv#1{ \xcMF_{#1} }

\def\xcMFnz{ \xcMFvv{\varnothing}{0} }
\def\xcMFyot{ \xcMFvv{y_1,y_2}{y_1^2+y_2^2} }

\def\xcMFa{ \maugv{\mathsf{MF}} }
\def\xcMFad{ \maugv{\ddot{\xcMF}} }

\def\xcMFvv#1#2{ \xcMF(#1;#2) }
\def\xcMFbxW{ \xcMFvv{\bax}{\xcW} }
\def\xcMFbxpW{ \xcMFvv{\baxp}{\xcW} }
\def\xcMFbxWo{ \xcMFvv{\bax}{\xcWo} }
\def\xcMFbxWt{ \xcMFvv{\bax}{\xcWt} }

\def\xcMFavv#1#2{ \xcMFa_{#1;#2} }
\def\xcMFabxW{ \xcMFavv{\bax}{\xcW} }

\def\xcMFdv#1{ \xcMFd(#1) }
\def\xcMFdbx{ \xcMFdv{\bax} }

\def\xcMFdby{ \xcMFdv{\bay} }

\def\xcMFdbxy{ \xcMFdv{\bax,\bay} }
\def\xcMFdbxz{ \xcMFdv{\bax,\baz} }
\def\xcMFdbyz{ \xcMFdv{\bay,\baz} }

\def\xcMFadv#1{ \xcMFad(#1) }
\def\xcMFadbx{ \xcMFadv{\bax} }

\def\xcMFbxyWt{ \xcMFvv{\bax,\ayo,\ayt}{\xcW(\bax) + \ayo^2 +
\ayt^2} }
\def\xcMFbxpyWt{ \xcMFvv{\baxp,\ayo,\ayt}{\xcW(\baxp) + \ayo^2 +
\ayt^2} }

\def\ayots{ \ayo^2 + \ayt^2 }

\def\xcMFA{ \mathsf{MFA} }
\def\xcMFAd{ \ddot{\xcMFA} }

\def\xcMFAvv#1#2{ \xcMFA(#1;#2) }
\def\xcMFAbxW{ \xcMFAvv{\bax}{\xcW} }

\def\xcMFAdv#1{ \xcMFAd(#1) }
\def\xcMFAdbx{ \xcMFAdv{\bax} }

\def\shC{ C }

\def\xO{ O }
\def\xOv#1{ \xO_{#1} }

\def\ycY{ \mathcal{Y} }

\def\ycYov#1{ \ycY_{#1} }

\def\yZ{ Z }
\def\yY{ Y }

\def\ycYp{ \ycY\p }
\def\yYp{ \yY\p }

\def\ycYav#1{ \ycY_{#1} }
\def\ycYo{ \ycYav{1} }
\def\ycYt{ \ycYav{2} }

\def\ycYi{ \ycYav{i} }
\def\ycYj{ \ycYav{j} }

\def\yYv#1{ \yY_{#1} }
\def\yYo{ \yYv{1} }
\def\yYt{ \yYv{2} }
\def\yYh{ \yYv{3} }
\def\yYi{ \yYv{i} }
\def\yYj{ \yYv{j} }
\def\yYot{ \yYv{12} }
\def\yYij{ \yYv{ij} }

\def\yYoct{ \yYo\cap\yYt }
\def\yYocct{ \yYo\cap\yYt }

\def\ycYav#1{ \ycY_{#1} }

\def\ycYaot{ \ycYav{12} }
\def\ycYath{ \ycYav{23} }

\def\xfbr#1#2#3{ \xymatrix@C=1.5pc@R=1.5pc{#1\ar[r] & #2 \ar[d] \\ & #3} }
\def\xfbrv#1#2#3#4{ \xymatrix@C=1.5pc@R=1.5pc{#1\ar[r] & #2 \ar[d]^{#4} \\ & #3} }

\def\som{ \omega }
\def\somv#1{ \som_{#1} }
\def\somo{ \somv{1} }
\def\somt{ \somv{2} }

\def\Xsom{ (\xX,\som) }
\def\Xsomm{ (\xX,-\som) }
\def\Xsomv#1{ (\xXv{#1},\somv{#1}) }
\def\Xsomo{ \Xsomv{1} }
\def\Xsomt{ \Xsomv{2} }
\def\Xsomh{ \Xsomv{3} }

\def\ycapv#1{ \times_{#1} }
\def\ycapX{ \ycapv{\xX} }
\def\ycapXott{ \ycapv{(\xXott)} }

\def\ycapU{ \ycapv{\xcA} }

\def\ycYoct{ \ycYo \ycapX \ycYt }

\def\zcap{ \boxtimes }
\def\zcap{ \tilde{\times} }
\def\zcapv#1{ \zcap_{#1} }

\def\zcapU{ \zcapv{\xcA} }

\def\xMor{\mathop{ \mathrm{Hom} }}

\def\yepi{ \pi }
\def\yepiu{ \yepi^\ast }
\def\yepiuv#1{ \yepiu_{#1} }
\def\yepiuo{ \yepiuv{1} }
\def\yepiut{ \yepiuv{2} }
\def\yepiuot{ \yepiuv{12} }
\def\yepiuth{ \yepiuv{23} }

\def\yepidv#1{ \yepi_{#1,\ast} }

\def\yepidt{ \yepidv{2} }
\def\yepidoh{ \yepidv{13} }

\def\yepiv#1{ \yepi_{#1} }
\def\yepio{ \yepiv{1} }
\def\yepit{ \yepiv{2} }
\def\yepiij{ \yepiv{ij} }
\def\yepiot{ \yepiv{12} }
\def\yepith{ \yepiv{23} }
\def\yepioh{ \yepiv{13} }

\def\spprt{support}

\def\xsupp{ \mathop{\mathrm{supp}} }

\def\Zt{ \ZZ_2 }
\def\Ztgrdd{$\Zt$-graded}
\def\Ztgrdng{$\Zt$-grading}
\def\Zgrdd{$\ZZ$-graded}
\def\Zgrdng{$\ZZ$-grading}
\def\Ztdgr{$\Zt$-degree}
\def\Zdgr{$\ZZ$-degree}

\def\dgZt{ \deg_{\Zt} }

\def\xPh{ \Phi }
\def\xPhv#1{ \xPh\lrbs{#1} }

\def\xPhYot{ \xPhv{\ycYaot} }
\def\dPhYot{ \dPhv{\ycYaot} }

\def\xPhcF{ \xPhv{\ycF} }

\def\xPhuv#1{ \xPh_{#1} }
\def\xPhuot{ \xPhuv{12} }
\def\xPhuotv#1{ \xPhuot\lrbs{#1} }
\def\xPhuotWzot{ \xPhuotv{\nWzot} }

\def\dPh{ \ddot{\Phi} }
\def\dPhv#1{ \dPh\lrbs{#1} }
\def\dPhWzot{ \dPhv{\nWzot} }

\def\xPhUWot{ \dPhv{\cmfnUWot} }

\def\xPhMot{ \xPhv{\cfMot} }

\def\xPhEot{ \xPhv{\xcEot} }
\def\xPhEpot{ \xPhv{\xcEot\p} }
\def\xPhEth{ \xPhv{\xcEth} }
\def\xPhEpth{ \xPhv{\xcEth\p} }

\def\xvprvv#1#2{ \lrbc{#1\,;#2} }
\def\xvprbvv#1#2{ \brb{#1\,;#2} }

\def\nvprvv#1#2{ #1\,;#2 }

\def\xrarv#1{ \xrightarrow{\;\;#1\;\;} }

\def\any{ \forall }

\def\xdz{ d_0 }
\def\xdo{ d_1 }

\def\tpcmpv#1#2{ \xymatrix{ #1_0\ar@<0.5ex>[r]^{\xdz} & \#1_1\ar@<0.5ex>[l]^{\xdo}
} }

\def\xobj{object}
\def\xobjv#1{$#1$-\xobj}

\def\zobj{\xobjv{0}}

\def\xIdv#1{ \xId_{#1} }

\def\xIdXsom{ \xIdv{\Xsom} }

\def\xIdbx{ \xIdv{\xlobx} }

\def\xopt{ 1-\mathrm{pt} }
\def\xoptv#1{ #1_{\xopt} }
\def\xXopt{ \xoptv{\xX} }

\def\xcAopt{ \xoptv{\xcA} }

\def\cO{ \mathscr{O} }

\def\Trf{ \mathcal{T} }
\def\Trfv#1{ \Trf_{#1} }
\def\TrfY{ \Trfv{\yY} }

\def\TrfDlX{ \Trfv{\DlX} }

\def\TrfU{ \Trfv{\xcA} }

\def\xopt{ \mathrm{1-pt} }
\def\xopt{ \text{1-pt} }
\def\xXopt{ \xX_{\xopt} }

\def\xtn{ \otimes }
\def\xtnv#1{ \xtn }

\def\xXo{ \xX_1 }
\def\xXt{ \xX_2 }
\def\xXott{ \xXo\times\xXt }

\def\HmB#1{ \Hom\Big( #1 \Big) }
\def\Hmb#1{ \Hom\big( #1 \big) }
\def\Hmo#1{ \Hom( #1 ) }

\def\mF{ F }
\def\mFua{ \mF^\ast }
\def\mFda{ \mF_{\ast} }

\def\dF{ \ddot{\mF} }
\def\dFua{ \dF^\ast }
\def\dFda{ \dF_{\ast} }

\def\lvar#1#2{ #1_{1},\ldots,#1_{#2} }

\def\ax{ x }

\def\ay{ y }
\def\ayv#1{ \ay_{#1} }
\def\ayo{ \ayv{1} }
\def\ayt{ \ayv{2} }

\def\az{ z }

\def\au{ u }
\def\aw{ w }
\def\av{ v }

\def\ap{ p }

\def\aq{ q }

\def\aa{ a }

\def\af{ f }
\def\ag{ g }

\def\afcPv#1{ \af_{\cP,#1} }

\def\xR{ \ICbx }
\def\xR{ R }

\def\xRtt{ \xR^{\xtt} }

\def\ICbxv#1{ \ICbx_{#1} }
\def\ICbxo{ \ICbxv{\yev} }
\def\ICbxz{ \ICbxv{\yod} }

\def\xlov#1{ #1 }
\def\xlobx{ \xlov{\bax} }
\def\xloby{ \xlov{\bay} }
\def\xlobxy{ \xlov{\bax,\bay} }
\def\xlobe{ \xlov{\varnothing} }

\def\bax{ \mathbf{\ax} }
\def\baxp{ \bax\p }
\def\baxpp{ \bax\pp }
\def\bay{ \mathbf{\ay} }
\def\baz{ \mathbf{\az} }
\def\bau{ \mathbf{\au} }
\def\baup{ \bau\p }
\def\baw{ \mathbf{\aw} }

\def\baf{ \mathbf{\af} }
\def\baa{ \mathbf{\aa} }

\def\bap{ \mathbf{\ap} }
\def\baq{ \mathbf{\aq} }

\def\bafv#1{ \baf_{#1} }

\def\bafcP{ \bafv{\cP} }
\def\bafcPp{ \bafv{\cPp} }

\def\Koz{ \mathrm{K} }
\def\xKov#1{ \Koz\lrbc{#1} }

\def\xKobp{ \xKov{\bap} }
\def\xKmfv#1#2{ \Koz(#1;#2) }
\def\xKmfbpq{ \xKmfv{\bap}{\baq} }

\def\xKpvv#1#2{ \mathrm{K}_{#1}(#2) }
\def\xKpWttbp{ \xKpvv{\xcWtt}{\bap} }

\def\vrp#1{ #1^{\prime} }
\def\vrpp#1{ #1^{\prime\prime} }

\def\baxp{ \vrp{\bax} }
\def\baxpp{ \vrpp{\bax} }

\def\baup{ \vrp{\bau} }

\def\bazp{ \vrp{\baz} }
\def\bazpp{ \vrpp{\baz} }

\def\bawpp{ \vrpp{\baw} }

\def\ICv#1{ \IC[#1] }
\def\ICbx{ \ICv{\bax} }
\def\ICby{ \ICv{\bay} }
\def\ICbxy{ \ICv{\bax,\bay} }

\def\ICav#1{ \IC_{#1} }

\def\ICa{ \ICav{\aa} }

\def\Tngx{ \Tng_x }

\def\smu{ \alpha }

\def\pIvv#1#2{
\xy
(-5,0)*{\bullet};(5,0)*{\bullet}**@{-},
(-5,4)*{\obo},(5,4)*{\obop}
\endxy
}

\def\xdlL{ \mathfrak{L} }
\def\xdlLv#1{ \xdlL_{#1} }
\def\xdlLa{ \xdlLv{\pga} }

\def\xdlA{ \mathcal{A} }
\def\xdlAuv#1{ \xdlA^{#1} }

\def\xdlAe{ \xdlAuv{\yev} }
\def\xdlAo{ \xdlAuv{\yod} }
\def\xdlAZz{ \xdlAuv{0} }

\def\shx{ y }
\def\shbx{ \bay }

\def\xdlAv#1{ \xdlA_{#1} }
\def\xdlAx{ \xdlAv{\shx^2} }
\def\xdlAbx{ \xdlAv{\shbx^2} }

\def\nb{ \nabla }
\def\nbv#1{ \nb_{\!\!#1} }
\def\nbE{ \nbv{\xE} }
\def\nbEp{ \nbE\p }
\def\nbEot{ \nbv{\xEot} }
\def\nbEth{ \nbv{\xEth} }

\def\nbb{ \bar{\nabla} }
\def\nbbv#1{ \nbb_{\!\!#1} }
\def\nbboth{ \nbbv{123} }
\def\nbbE{ \nbbv{\xE} }

\def\nbbEi{ \nbbv{\xEi} }

\def\xdgmvv#1#2{ \big(\xbOmbv{#1},\nbbv{#1} \big) }

\def\xdgmv#1{ \xdgmvv{#1}{#1} } 

\def\xdgsv#1{ \big(#1,\nbbv{#1} \big) }
\def\xdgmE{ \xdgmv{\xE} }

\def\xdgmP{ \xdgsv{\xmP} }

\def\xdgmM{ \xdgsv{\xmM} }

\def\adgmv#1{ (#1,\nbbv{#1}) }
\def\adgmE{ \adgmv{\xE} }
\def\adgmEv#1{ \adgmE_{#1} }
\def\admgEzdfm{ \adgmEv{\zdfm} }

\def\adgmpE{ (\xE,\nbbE\p) }

\def\adgmev#1{ (#1,\nbbv{#1}) }
\def\adgmeE{ \adgmev{\xE} }
\def\adgmeEp{ \adgmev{\xEp} }
\def\adgmeEo{ \adgmev{\xEo} }
\def\adgmeEt{ \adgmev{\xEt} }
\def\adgmeEh{ \adgmev{\xEh} }
\def\adgmeEot{ \adgmev{\xEot} }

\def\adgmeEth{ \adgmev{\xEth} }

\def\nbbM{ \nbbv{\xmM} }
\def\nbbMxo{ \nbbv{\xmMo} }
\def\nbbMxt{ \nbbv{\xmMt} }

\def\nbbMv#1{ \nbbv{\xmM^{#1}} }
\def\nbbMe{ \nbbMv{\yev} }
\def\nbbMo{ \nbbMv{\yod} }

\def\xdD{ D }
\def\xdDv#1{ \xdD_{#1} }
\def\xmdv#1{ \big( #1, \xdDv{#1} \big) }

\def\xmdM{ \xmdv{\xmM} }

\def\xdDM{ \xdDv{\xmM} }

\def\EndICv#1{ \End_{\ICv{#1}} }
\def\EndICbx{ \EndICv{\bax} }

\def\xIdv#1{ \xId_{#1} }

\def\xmP{ P }

\def\xmM{ M }
\def\xmMv#1{ \xmM_{#1} }
\def\xmMo{ \xmMv{1} }
\def\xmMt{ \xmMv{2} }

\def\xmMv#1{ \xmM_{#1} }

\def\xmMuv#1{ \xmM^{#1} }

\def\xmPv#1{ \xmP_{#1} }
\def\xmPo{ \xmPv{1} }
\def\xmPt{ \xmPv{2} }

\def\cmfv#1#2{ (#1,#2) }
\def\cmfBv#1#2{ \Big(#1,#2\Big) }

\def\cmfUW{ \cmfv{\xcA}{\xcW} }
\def\cmfUWv#1{ \cmfv{\xcA_{#1}}{\xcW_{#1}} }
\def\cmfUWo{ \cmfUWv{1} }
\def\cmfUWt{ \cmfUWv{2} }
\def\cmfUWh{ \cmfUWv{3} }

\def\cmfcUW{ \cmfv{\xcU}{\xcW} }
\def\cmfUW{ \cmfv{\xcA}{\xcW} }
\def\cmfUWv#1{ \cmfv{\xcAv{#1}}{\xcWv{#1}} }
\def\cmfUWo{ \cmfUWv{1} }
\def\cmfUWt{ \cmfUWv{2} }
\def\cmfcUWv#1{ \cmfv{\xcUv{#1}}{\xcWv{#1}} }
\def\cmfcUWo{ \cmfcUWv{1} }
\def\cmfcUWt{ \cmfcUWv{2} }
\def\cmfcUWh{ \cmfcUWv{3} }
\def\cmfcUWot{ \cmfcUWv{12} }
\def\cmfcUWth{ \cmfcUWv{23} }

\def\cmfcUWp{ \cmfv{\xcUp}{\xcWp} }

\def\cmfcUmWv#1{ \cmfv{\xcUv{#1}}{-\xcWv{#1}} }
\def\cmfcUmWt{ \cmfcUmWv{2} }

\def\nmfv#1#2{ #1,#2 }

\def\cmfnUWv#1{ \nmfv{\xcUv{#1}}{\xcWv{#1}} }

\def\cmfnUWot{ \cmfnUWv{12} }

\def\rmO{ \mathrm{O} }
\def\OnC{ \rmO(n,\IC) }

\def\ycF{ \cF }
\def\ycE{ \cE }
\def\ycEp{ \ycE\p }
\def\ycL{ \mathcal{L} }

\def\ycEv#1{ \ycE_{#1} }
\def\ycEo{ \ycEv{1} }

\def\ycFv#1{ \ycF_{#1} }
\def\ycFot{ \ycFv{12} }
\def\ycFth{ \ycFv{23} }

\def\ay{ y }

\def\Exttv#1{ \Ext^{#1} }
\def\Exttev{ \Exttv{\yev} }
\def\Exttod{ \Exttv{\yod} }

\def\Extub{ \Ext^\bullet }

\def\xcWq{ \xcW_{\mathrm{q}} }

\def\vcL{ L }

\def\vcLv#1{ \vcL_{#1} }

\def\vcLz{ \vcLv{0} }

\def\vcLot{ \vcLv{12} }

\def\vccL{ Q }
\def\vccLv#1{ \vccL_{#1} }
\def\vccLi{ \vccLv{i} }
\def\vccLo{ \vccLv{1} }
\def\vccLt{ \vccLv{2} }

\def\vccLd{ \dulv{\vccL} }
\def\vccLdv#1{ \vccLd_{#1} }
\def\vccLdi{ \vccLdv{i} }

\def\xfv#1{ f_{\mathrm{b},#1} }
\def\fvL{ \xfv{\vcL} }
\def\xafv#1{ f_{#1} }
\def\fvcE{ \xafv{\ycE} }
\def\fvcEp{ \xafv{\ycEp} }

\def\yobv#1{ \brb{#1,\xafv{#1}} }
\def\yobcE{ \yobv{\ycE} }
\def\yobcEp{ \yobv{\ycEp} }

\def\brB#1{ \Big( #1 \Big) }
\def\brb#1{ \big( #1 \big) }
\def\bro#1{ (#1) }

\def\shlf{ \tfrac{1}{2} }

\def\ypv#1{ p_{#1} }
\def\ypY{ \ypv{\ycY} }
\def\ypYi{ \ypv{\ycYi} }

\def\xcWlc{ \xcW_{\mathrm{lc}} }

\def\vcLY{ \vcLv{\ycY} }
\def\vcLYv#1{ \vcLv{\ycYov{#1}} }
\def\vcLYo{ \vcLYv{1} }
\def\vcLYt{ \vcLYv{2} }

\def\vcLYi{ \vcLYv{i} }

\def\vLY{ \vcLv{\yY} }
\def\vLYv#1{ \vcLv{\yYv{#1}} }
\def\vLYo{ \vLYv{1} }
\def\vLYt{ \vLYv{2} }

\def\govv#1#2{ \lrbc{#1,#2} }
\def\goYL{ \govv{\ycY}{\vcLY} }
\def\goYLv#1{ \govv{\ycYov{#1}}{\vcLYv{#1}} }
\def\goYLo{ \goYLv{1} }
\def\goYLt{ \goYLv{2} }

\def\gYL{ \govv{\yY}{\vLY} }
\def\gYLv#1{ \govv{\yYv{#1}}{\vLYv{#1}} }
\def\gYLo{ \gYLv{1} }
\def\gYLt{ \gYLv{2} }

\def\Hdlb{ \rmH_{\dlb} }
\def\HnbbE{ \rmH_{\nbbE} }
\def\Hdlbv#1{ \Hdlb^{#1} }

\def\Hdlbo{ \Hdlbv{1} }

\def\rmS{ \mathrm{S} }

\def\TU{ \Tng\xcA }

\def\mtdfv#1#2{ (#1)_{#2} }

\def\mtdfTUSk{ \mtdfv{\TsU}{\hdf} }
\def\mtdfTYSk{ \mtdfv{\TsY}{\hdf} }
\def\mtdfTYoSk{ \mtdfv{\TsY_1}{\hdf} }
\def\mtdfTUh{ \mtdfv{\TsU}{\hdf} }

\def\bdmu{ \mu }
\def\bdmuv#1{ \bdmu_{#1} }
\def\bdmuot{ \bdmuv{12} }
\def\bdmuth{ \bdmuv{23} }
\def\bdmuoh{ \bdmuv{13} }

\def\tdmupr{ \tilde{\bdmu}_\approx }

\def\bdmuoha#1{ \bdmuoh \spsmb #1 }

\def\bdprv#1{ #1_{\approx} }
\def\bdmupr{ \bdprv{\bdmu} }
\def\bdnupr{ \bdprv{\bdnu} }
\def\bdxipr{ \bdprv{\bdxi} }

\def\bdnu{ \nu }
\def\bdnuv#1{ \bdnu_{#1} }
\def\bdnuot{ \bdnuv{12} }

\def\bdxi{ \xi }
\def\bdxiv#1{ \bdxi_{#1} }
\def\bdxiot{ \bdxiv{12} }

\def\Lbrv#1{ [ #1 ] }
\def\Pbrv#1{ \{ #1 \} }

\def\hdf{\varkappa}
\def\hdfv#1{ \hdf_{|#1} }
\def\hdfi{ \hdfv{i} }

\def\hdfo{ \hdfv{1} }
\def\hdftw{ \hdfv{2} }
\def\hcdf{ \cdlv{\hdf} }
\def\hcdfv#1{ \hcdf_{|#1} }
\def\hcdfo{ \hcdfv{1} }

\def\hdfe{ \dfev{\hdf} }
\def\ehdf{ \dfe\hdf }

\def\hdfbt{ \beta }
\def\hdfgm{ \gamma }

\def\ccrR{ \check{\crR} }
\def\cdfbt{ \check{\beta} }
\def\cdfbtY{ \cdfbt_{\yY} }
\def\tdfbt{ \tilde{\beta} }
\def\tdfbtv#1{ \tdfbt_{#1} }
\def\tdfbtY{ \tdfbtv{\yY} }
\def\tdfbtYo{ \tdfbtv{\yYo} }
\def\tdfbtYt{ \tdfbtv{\yYt} }

\def\tdfbtU{ \tdfbtv{\xcA} }

\def\Sta{ \mathrm{S} }
\def\Stav#1{ \Sta^{#1} }
\def\Sb{ \Stav{\bullet} }
\def\Stat{ \Stav{2} }

\def\Stai{ \Stav{i} }

\def\SbTU{ \Sb\TU }

\def\Wb{ \wedge^{\bullet} }
\def\WbTU{ \Wb\TU }

\def\ydah{$(0,\yod)$}
\def\ydh{$(1,0)$}

\def\acF{ F }
\def\hacF{ \hat{F} }
\def\hacF{ \check{F} }
\def\acFv#1{ \acF_{#1} }
\def\hacFv#1{ \hacF_{#1} }
\def\hacFEo{ \hacFv{\xEo} }
\def\hacFEt{ \hacFv{\xEt} }

\def\acFpv#1{ \acF\p_{#1} }
\def\acFE{ \acFv{\xE} }
\def\hacFE{ \hacFv{\xE} }
\def\acFpE{ \acFpv{\xE} }
\def\acFEp\acFv{\xEp}

\def\acFEo{ \acFv{\xEo} }
\def\acFEt{ \acFv{\xEt} }
\def\acFEh{ \acFv{\xEh} }
\def\acFEi{ \acFv{\xEi} }
\def\acFEot{ \acFv{\xEot} }

\def\acFEth{ \acFv{\xEth} }

\def\dOv#1{ O_{#1} }
\def\dOk{ \dOv{k} }

\def\dfe{ \epsilon }
\def\dfev#1{ #1_{\dfe} }

\def\dfzzv#1{ #1_0 }
\def\dfov#1{ #1_1 }

\def\acWv#1{ \xcW_{|#1} }
\def\acWz{ \acWv{0} }
\def\acWo{ \acWv{1} }

\def\xcWe{ \dfev{\xcW} }

\def\xcWz{ \dfzv{\xcW} }
\def\xcWzz{ \dfzzv{\xcW} }
\def\xcWo{ \dfov{\xcW} }

\def\xcWvv#1#2{ \xcWv{#1|#2} }
\def\xcWve#1{ \xcWvv{#1}{\dfe} }
\def\xcWvz#1{ \xcWvv{#1}{0} }

\def\xcWoz{ \xcWvz{1} }
\def\xcWtz{ \xcWvz{2} }

\def\xcWotz{ \xcWvz{12} }
\def\xcWthz{ \xcWvz{23} }

\def\xcWoto{ \xcWvo{12} }

\def\dfmvv#1#2{ #1_{#2} }

\def\dfmUmv#1{ \dfmvv{\xcA}{\bdmuv{#1}} }

\def\nbbvv#1#2{ \nbb_{#1\,|\, #2} }
\def\nbbvv#1#2{ \nbb_{#1| #2} }

\def\nbbzv#1{ \nbbvv{#1}{0} }
\def\nbbov#1{ \nbbvv{#1}{1} }

\def\nbbzE{ \nbbzv{\xE} }
\def\nbboE{ \nbbov{\xE} }

\def\nbbzEot{ \nbbzv{\xEot} }

\def\nbbzEoh{ \nbbzv{\xEoh} }

\def\nbbzcxE{ \nbbzv{\cxrE} }
\def\nbbocxE{ \nbbov{\cxrE} }

\def\sea{ a }
\def\seb{ b }

\def\msgm{ \sigma }
\def\msgmv#1{ \msgm_{#1} }

\def\msgmbv#1{ \msgmv{|#1} }
\def\msgmbz{ \msgmbv{0} }
\def\msgmbo{ \msgmbv{1} }
\def\msgmvv#1#2{ \msgmv{#1\,|\,#2} }
\def\msgmve#1{ \msgmvv{#1}{\dfe} }
\def\msgmvz#1{ \msgmvv{#1}{0} }

\def\msgmot{ \msgmv{12} }
\def\msgmote{ \msgmve{12} }
\def\msgmotz{ \msgmvz{12} }

\def\msgmth{ \msgmv{23} }
\def\msgmthe{ \msgmve{23} }
\def\msgmthz{ \msgmvz{23} }

\def\xIdE{ \xIdv{\xE} }

\def\cxrv#1{ #1 }
\def\cxrA{ \cxrv{A} }
\def\cxrmM{ \cxrv{\xmM} }
\def\cxrE{ \cxrv{\xE} }

\def\shlf{ \tfrac{1}{2} }

\def\TrfUZ{ (\TrfU,0) }
\def\TrfUW{ (\TrfU,\xcW) }

\def\yYW{ \yY_{\xcW} }
\def\yYWv#1{ \yY_{\xcWv{#1} }}
\def\yYWo{ \yYWv{1} }
\def\yYWt{ \yYWv{2} }
\def\yYWi{ \yYWv{i} }

\def\yYUW{ \yY_{\cmfcUW} }

\def\obo{B}

\def\obop{ \obo\p }

\def\xcd#1#2{ #1\cdot #2 }
\def\pdq{ \xcd{\bap}{\baq} }

\def\rDbU{ \rDb(\xcA) }

\def\ttrnv#1{ [#1]_2 }
\def\ttrno{ \ttrnv{1} }
\def\ttrnt{ \ttrnv{2} }
\def\ttrnn{ \ttrnv{n} }

\def\ytrnv#1{ [ #1 ]_{\mathrm{tw}} }

\def\ytrnL{ \ytrnv{\vcL} }
\def\ytrnLot{ \ytrnv{\vcLot} }
\def\ytrnB{ \ytrnv{\bnB} }

\def\ztrnv#1{ \ytrnv{#1}\p }

\def\ztrnL{ \ztrnv{\vcL} }

\def\spsmb{ \,\llcorner\, }

\def\nboth{ \zdfmv{12,23} }

\def\aA{ \mathscr{A} }

\def\xmi{ i }
\def\xmj{ j }
\def\xmn{ n }

\def\xmvmlt#1{$#1$-multiplication}
\def\xmnmlt{\xmvmlt{\xmn}}

\def\xmvact#1{$#1$-action}
\def\xmnact{\xmvact{\xmn}}

\def\xmm{ a }
\def\xmmbf{ \mathbf{\xmm} }
\def\xmmv#1{ \xmm_{#1} }
\def\xmmz{ \xmmv{0} }
\def\xmmo{ \xmmv{1} }
\def\xmmt{ \xmmv{2} }
\def\xmmn{ \xmmv{\xmn} }

\def\zmmv#1{ \xmm^{#1} }
\def\zmmM{ \zmmv{\xmM} }
\def\zmmMv#1{ \zmmM_{#1} }

\def\zmmMo{ \zmmMv{1} }
\def\zmmMt{ \zmmMv{2} }
\def\zmmMn{ \zmmMv{\xmn} }

\def\zmmbfv#1{ \xmmbf^{#1} }
\def\zmmbfM{ \zmmbfv{\xmM} }

\def\zmmE{ \zmmv{\xE} }
\def\zmmEv#1{ \zmmE_{#1} }

\def\zmmEn{ \zmmEv{\xmn} }

\def\zmmbfE{ \zmmbfv{\xE} }

\def\zmmuMv#1{ a^{\xmMv{#1}} }
\def\zmmuMo{ \zmmuMv{1} }
\def\zmmuMt{ \zmmuMv{2} }

\def\zmf{f}
\def\zmfv#1{ \zmf_{#1} }
\def\zmfo{ \zmfv{1} }
\def\zmft{ \zmfv{2} }
\def\zmfn{ \zmfv{\xmn} }
\def\zmbf{ \mathbf{\zmf} }

\def\xmmuz{ \xmm^{(0)} }
\def\xmmbfuz{ \xmmbf^{(0)} }
\def\xmmuzv#1{ \xmmuz_{#1} }
\def\xmmuzo{ \xmmuzv{1} }
\def\xmmuzt{ \xmmuzv{2} }
\def\xmmuzn{ \xmmuzv{\xmn}}

\def\xmmuzw#1{ \xmm^{#1,(0)} }
\def\xmmuzE{ \xmmuzw{\xE} }
\def\xmmuzEv#1{ \xmmuzE_{#1} }
\def\xmmuzEo{ \xmmuzEv{1} }
\def\xmmuzEt{ \xmmuzEv{2} }
\def\xmmuzEn{ \xmmuzEv{\xmn}}

\def\xmmbfuzw#1{ \xmmbf^{#1,(0)} }
\def\xmmbfuzE{ \xmmbfuzw{\xE} }

\def\xdfm{ \lambda }

\def\xdfmv#1{ \xdfm_{#1} }
\def\xdfmi{ \xdfmv{i} }
\def\xdfmz{ \xdfmv{0} }
\def\xdfmo{ \xdfmv{1} }

\def\xdfmot{ \xdfmv{12} }

\def\adfmv#1{ \xdfm_{\cap,#1} }
\def\adfmot{ \adfmv{12} }
\def\adfmotv#1{ \adfmv{12,#1} }
\def\adfmoti{ \adfmotv{i} }
\def\adfmotz{ \adfmotv{0} }
\def\adfmoto{ \adfmotv{1} }

\def\xdfmotv#1{ \xdfmv{12,#1} }
\def\xdfmoti{ \xdfmotv{\xmi} }
\def\xdfmotn{ \xdfmotv{\xmn} }
\def\xdfmotnmo{ \xdfmotv{\xmn-1} }
\def\xdfmotz{ \xdfmotv{0} }
\def\xdfmoto{ \xdfmotv{1} }
\def\xdfmott{ \xdfmotv{2} }
\def\xdfmoth{ \xdfmotv{3} }

\def\dlbW{ \dlbv{_{\xcW}} }

\def\ydfmotv#1{ \bdmuot \spsmb #1 }

\def\zdfm{ \xcW }
\def\zdfmv#1{ \zdfm_{#1} }
\def\zdfmo{ \zdfmv{1} }

\def\dgDlb{ \deg_{\mathrm{Dlb}} }
\def\dgDlbv#1{ \dgDlb #1 }

\def\yYviv#1#2{ \yYv{#1}\cap\yYv{#2} }
\def\yYoit{ \yYviv{1}{2} }
\def\yYtih{ \yYviv{2}{3} }
\def\yYoih{ \yYviv{1}{3} }
\def\yYothi{ \yYo\cap\yYt\cap\yYh }

\def\bfp{ \mathbf{p} }

\def\Cr{ \mathrm{Cr} }
\def\Crv#1{ \Cr_{#1} }
\def\CrW{ \Crv{\xcW} }
\def\rCrW{ |_{\CrW} }
\def\CrWtmo{ \Crv{\xcWt-\xcWo} }

\def\yrl{\mathrm{rel}}
\def\dgrl{ \deg_{\yrl} }
\def\ydgm{ m }
\def\ydgn{ n }

\def\tHom{\widetilde{\Hom}}

\def\dfib{ \del_{\mathrm{vrt}} }

\def\xaV{ V }

\def\xbV{ \xaV^\prime }
\def\xaa{ \alpha }
\def\xav{ v }
\def\xavv#1{ \xav_{#1} }
\def\xavo{ \xavv{1} }
\def\xavt{ \xavv{2} }
\def\xavh{ \xavv{3} }

\def\fsd{ \del_{\mathrm{s}} }
\def\ssd{ \fsd^2 }

\def\xcmv#1{ |_{#1} }

\def\xcmi{ \xcmv{i} }

\def\tcW{ \widetilde{\xcW} }

\def\hcW{ \widehat{\xcW} }
\def\hcWv#1{ \hcW_{#1} }
\def\hcWo{ \hcWv{1} }
\def\hcWt{ \hcWv{2} }
\def\tmu{ \tilde{\mu} }

\def\xcWov#1{ \xcW_{1|#1} }
\def\xcWoz{ \xcWov{0} }

\def\xcWtv#1{ \xcW_{2|#1} }
\def\xcWtz{ \xcWtv{0} }

\def\xcWotv#1{ \xcW_{12|#1} }
\def\xcWotz{ \xcWotv{0} }
\def\xcWoto{ \xcWotv{1} }

\def\hcWot{ \hcW_{12} }

\def\vrx{ x }
\def\vrxv#1{ \vrx^{#1} }
\def\vrxI{ \vrxv{\inI} }
\def\vrxJ{ \vrxv{\inJ} }
\def\vrxK{ \vrxv{\inK} }

\def\inI{ I }
\def\inJ{ J }
\def\inK{ K }
\def\inL{ L }
\def\inM{ M }
\def\vrxI{ \vrxv{\inI} }

\def\dlv#1{ \del_{#1} }
\def\dlI{ \dlv{\inI} }
\def\dlJ{ \dlv{\inJ} }
\def\dlK{ \dlv{\inK} }
\def\dlL{ \dlv{\inL} }

\def\nbcv#1{ \nabla_{#1} }

\def\nbcK{ \nbcv{\inK} }
\def\nbcL{ \nbcv{\inL} }

\def\crR{ R }

\def\OWv#1{ O(\xcW^{#1}) }

\def\OWo{ O(\xcW) }
\def\OWt{ \OWv{2} }
\def\OWh{ \OWv{3} }
\def\OWn{ \OWv{n} }

\def\xcAekot{ \xcA_{\dfe\bdmuot} }
\def\xcAekth{ \xcA_{\dfe\bdmuth} }
\def\xcAekoh{ \xcA_{\dfe\bdmuoh} }

\def\xcEbz{ \xcE_{|0} }
\def\xcEvbz#1{ \xcE_{#1|0} }
\def\xcEobz{ \xcEvbz{1} }
\def\xcEtbz{ \xcEvbz{2} }

\def\shdfv#1{ \hdf_{#1} }
\def\shdfij{ \fsd\shdfv{ij} }
\def\shdfijk{ \ssd\shdfv{ijk} }
\def\shdfot{ \fsd\shdfv{12} }
\def\shdfoth{ \ssd\shdfv{123} }

\def\dulv#1{ #1\check }
\def\bnB{ B }

\def\tbw{ {\textstyle \bigwedge} }
\def\xtbw{ \tbw }
\def\wdtpv#1{ {\textstyle \bigwedge}^{\mathrm{ top}} #1 }
\def\wdtpB{ \wdtpv{\bnB} }
\def\Cr{ \mathop{\mathrm{Crit}} }
\def\Crv#1{ \Cr(#1) }

\def\CrW{ \Crv{\xcW} }

\def\edfn{ := }

\def\cnKv#1{ K_{#1} }
\def\cnKY{ \cnKv{\yY} }
\def\cnKU{ \cnKv{\xcA} }
\def\cnKYot{ \cnKv{\yYot} }
\def\cnKYi{ \cnKv{\yYi} }
\def\cnKYo{ \cnKv{\yYo} }
\def\cnKYt{ \cnKv{\yYt} }
\def\cnKYij{ \cnKv{\yYij} }

\def\yYoti{ \yYot }
\def\yYoct{ \yYoti }

\def\dulv#1{ #1^{\vee} }
\def\xdulv#1{ \dulv{#1} }
\def\ott{ ^{\otimes 2} }

\def\Vd{ \dulv{\xaV} }

\def\TY{ \Tng\yY }
\def\NY{ \Nrm\yY }
\def\TsY{ \Ts\yY }

\def\TsUhdf{ (\TsU)_\hdf }

\def\xfu{ u }
\def\xfVv#1{ \Vfbv{#1} }
\def\xfVu{ \xfVv{\xfu} }
\def\crcUv#1{ \xcU_{#1} }
\def\crcUW{ \crcUv{\xcW} }

\def\mpf{ f }
\def\zzf{ g }
\def\zzfv#1{ \zzf_{#1} }
\def\zzfo{ \zzfv{1} }
\def\zzft{ \zzfv{2} }
\def\zzfot{ \zzfv{12} }
\def\zzfi{ \zzfv{i} }

\def\yobYWv#1{ (\yYWv{#1},\zzfv{#1}^\ast\vcLz) }
\def\yobYWo{ \yobYWv{1} }
\def\yobYWt{ \yobYWv{2} }

\def\xcX{ \mathcal{X} }
\def\xcXv#1{ \xcX_{#1} }
\def\xcXo{ \xcXv{1} }
\def\xcXt{ \xcXv{2} }

\def\xcXopt{ \xcXv{\xopt} }

\def\xs{ s }
\def\xsv#1{ \xs_{#1} }
\def\xso{ \xsv{1} }
\def\xst{ \xsv{2} }
\def\xsopt{ \xsv{\xopt} }

\def\zNc{$\zN$-category}
\def\zNcmo{$(\zN-1)$-category}
\def\zNcmos{$(\zN-1)$-categories}
\def\zdlv#1{ #1^{\hve} }

\def\HomC{ \Hom_{\cC} }

\def\cCov#1{ \cC_{#1} }
\def\cCX{ \cCov{\xcX} }
\def\cCXo{ \cCov{\xcXo} }
\def\cCXt{ \cCov{\xcXt} }

\def\Tvs{ \rmT^\vee }

\def\yZv#1{ \yZ_{#1} }
\def\yZo{ \yZv{1} }
\def\yZt{ \yZv{2} }

\def\xzobj{object}
\def\xzobjs{\xzobj s}

\def\zdgz#1{ |#1|_{\ZZ} }
\def\zdgt#1{ |#1| }

\def\xblnk{ - }

\def\IdbxW{ \xIdv{\bax;\xcW} }

\def\dLt{ \ddot{\Lambda} }
\def\dLtp{ \dLt_+ }
\def\dLtm{ \dLt_- }
\def\dLtpm{ \dLt_\pm }

\def\dId{ \ddot{\xId} }
\def\dIdv#1{ \dId_{#1} }
\def\Idba{ \dIdv{\baa} }

\def\ocW{ \xcW }
\def\ocWv#1{ \ocW_{#1} }
\def\ocWo{ \ocWv{1} }
\def\ocWt{ \ocWv{2} }

\def\HmMFbx{ \Hom_{\xcMFdbx} }
\def\HmDDU{ \Hom_{\rDDprfU} }
\def\HmDDaU{ \Hom_{\rDDprfaU} }

\def\mgv#1{ g }

\def\dgdW{ \deg_{\del\xcW} }
\def\spdg{ \deg_{\mathrm{sc}} }

\def\xblnk{ - }

\def\tnsct{ \tau }
\def\bnbl{ \boldsymbol{\nabla} }
\def\Tfl{ \mathcal{T} }
\def\Tflvv#1#2{ \Tfl^{#1}_{#2} }
\def\ixk{ k }
\def\ixl{ l }
\def\Tflkl{ \Tflvv{\ixk}{\ixl} }
\def\TflklU{ \Tflkl(\xcA) }
\def\Tflbb{ \Tflvv{\bullet}{\bullet} }
\def\TflbbU{ \Tflbb(\xcA) }
\def\xTfl{ \Tfl }
\def\xTflU{ \xTfl(\xcA) }
\def\Tflnb{ \Tfl_{\bnbl} }
\def\Tflnbot{ \Tfl_{\bnbl,12} }
\def\Tflnv#1{ \Tflnb[#1] }

\def\dgbl{ \deg_{\mathrm{bal}} }
\def\dgdf{ \deg_{\mathrm{tot}} }

\def\eE{ E }
\def\eEv#1{ \eE_{#1} }
\def\eEo{ \eEv{1} }
\def\eEt{ \eEv{2} }
\def\eEh{ \eEv{3} }

\def\dfal{ \alpha }
\def\dfal{ \zeta }
\def\dfalv#1{ \dfal_{#1} }
\def\dfalot{ \dfalv{12} }
\def\dfalotv#1{ \dfalv{12,#1} }
\def\dfaloto{ \dfalotv{1} }
\def\dfalott{ \dfalotv{2} }

\def\dfaloti{ \dfalotv{i} }

\def\baxp{ \bax\p }
\def\bayp{ \bay\p }

\def\xcApv#1{ \xcA\p_{#1} }

\def\xcAbx{ \xcAv{\bax} }
\def\xcApbx{ \xcAp_{\bax} }
\def\xcAmbx{ \xcAv{-\bax} }
\def\xcApbxp{ \xcApv{\baxp} }
\def\xcApmbxp{ \xcApv{-\baxp} }

\def\ICnbv#1{ \IC^{\zn}_{#1} }
\def\ICnbax{ \ICnbv{\bax} }
\def\ICnbay{ \ICnbv{\bay} }
\def\xcVv#1{ \xcV_{#1} }

\def\xcVby{ \xcVv{\bay} }
\def\xcVmby{ \xcVv{-\bay} }
\def\xcVpv#1{ \xcV\p_{#1} }

\def\xcVpbyp{ \xcVpv{\bayp} }

\def\xcAxv#1{ \xcAv{#1,\bax_{#1}} }
\def\xcVyv#1{ \xcVv{#1,\bay_{#1}} }
\def\xcAxo{ \xcAxv{1} }
\def\xcAxt{ \xcAxv{2} }
\def\xcAxh{ \xcAxv{3} }
\def\xcVyo{ \xcVyv{1} }
\def\xcVyt{ \xcVyv{2} }
\def\xcVyh{ \xcVyv{3} }

\def\UxVyo{ \xcAxo\times\xcVyo }
\def\UxVyt{ \xcAxt\times\xcVyt }
\def\UxVyh{ \xcAxh\times\xcVyh }

\def\TsUbx{ \Ts\xcAbx }

\def\UxVy{ \xcAbx\times\xcVby }
\def\UpxVy{ \xcApbx\times\xcVby }
\def\UmxVmy{ \xcAmbx\times\xcVmby }

\def\VmyUx{ \xcVmby\times\xcAbx }
\def\VyUmx{ \xcVby\times\xcAmbx }
\def\UxVypp{ \xcApbxp\times\xcVpbyp }
\def\sUV{ \xcA\times\xcV }
\def\sUVpp{ \xcAp\times\xcVp }
\def\zsq{ q }
\def\zsp{ p }
\def\qfib{$\zsq$-fibration}
\def\pfib{$\zsp$-fibration}
\def\qemb{$\zsq$-embedding}

\def\smrc{symplectic rectangle}
\def\rcch{rectangular chart}

\def\sfq{ \varepsilon }
\def\rfq{ \sfq_{\zsq} }
\def\sfqv#1{ \sfq_{#1} }
\def\sfqth{ \sfqv{23} }
\def\sfqot{ \sfqv{12} }

\def\hve{ \diamondsuit }

\def\xtsm{ \tau }
\def\ddxtsm{ \ddot{\xtsm} }
\def\xtsmo{ \xtsm_1 }
\def\ddxtsmo{ \ddxtsm_1 }

\def\TsUv#1{ \Ts\xcAv{#1} }
\def\TsUo{ \TsUv{1} }
\def\TsUt{ \TsUv{2} }

\def\etfe{ \ddot{\mathsf{E}} }
\def\etfe{ \dPh_{\cong} }
\def\etfev#1{ \dPh_{\cong,#1} }
\def\etfeot{ \etfev{12} }
\def\etfeo{ \etfev{1} }
\def\etfet{ \etfev{2} }

\def\prms{ \sigma }

\def\rsf{ \dPh_{ \mathrm{r} } }
\def\rsfv#1{ \dPh_{ \mathrm{r},#1 } }
\def\rsfo{ \rsfv{1} }
\def\rsft{ \rsfv{2} }
\def\rsfe{ \rsfv{\sfq} }

\def\xemq{ \hookrightarrow_q }

\def\mch{ f }

\def\prsf{ \ddot{\mathfrak{P}} }
\def\prsfv#1{ \prsf #1 }
\def\prsfXsom{ \prsfv{\Xsom} }

\def\zzO{ O }
\def\zzOv#1{ \zzO_{#1} }
\def\zzOf{ \zzOv{\mch} }
\def\zzOfp{ \zzOv{\mch\p} }

\def\zzOo{ \zzOv{1} }
\def\zzOt{ \zzOv{2} }

\def\zzOfv#1{ \zzOv{#1,\mch_{#1} } }
\def\zzOfo{ \zzOfv{1} }
\def\zzOft{ \zzOfv{2} }

\def\ftnt#1{ \footnote{#1} }

\def\hsnv#1{ \del^2 #1 }
\def\hsnW{ \hsnv{\xcW} }
\def\hsnWb{ (\hsnW) }

\def\xEndsp{ \End_{\rDDprfUSk}(\zsU) }

\def\nbbEav#1{ \nbbv{\xE,#1} }
\def\nbbEai{ \nbbEav{i} }
\def\nbbEaz{ \nbbEav{0} }
\def\nbbEao{ \nbbEav{1} }
\def\nbbEat{ \nbbEav{2} }

\def\zdfm{ \zeta }
\def\zdfmv#1{ \zdfm_{#1} }
\def\zdfmo{ \zdfmv{1} }
\def\zdfmt{ \zdfmv{2} }
\def\zdfmn{ \zdfmv{\xmn} }
\def\zdfmi{ \zdfmv{i} }
\def\zdfmot{ \zdfmv{12} }
\def\zdfmoth{ \zdfmv{123} }
\def\zdfmothv#1{ \zdfmv{123,#1} }
\def\zdfmothn{ \zdfmothv{\xmn} }
\def\zdfmothp{ \zdfmoth\p }
\def\zdfmothnp{ \zdfmothn\p }

\def\tzdfm{ \tilde{\zdfm} }
\def\tzdfmv#1{ \tzdfm_{#1} }
\def\tzdfmoth{ \tzdfmv{123,\xmn} }

\def\zdfmotv#1{ \zdfmv{12,#1} }
\def\zdfmoti{ \zdfmotv{i} }
\def\zdfmotn{ \zdfmotv{\xmn} }
\def\zdfmoto{ \zdfmotv{1} }
\def\zdfmott{ \zdfmotv{2} }

\def\zdfmdv#1{ \zdfm[#1] }
\def\zdfmdvv#1#2{ \zdfm_{#1}[#2] }
\def\zdfmdnv#1{ \zdfmdvv{\xmn}{#1} }

\def\zas{ \alpha }
\def\zasv#1{ \zas_{#1} }
\def\zasoth{ \zasv{123} }
\def\zasothv#1{ \zasv{123,#1} }
\def\zasothi{ \zasothv{i} }
\def\zasothn{ \zasothv{\xmn} }
\def\zasotht{ \zasothv{2} }

\def\ycrc{ \circ }

\def\xIdEm{ \xIdv{\xEo\otimes\xEt\otimes\xEh} }

\def\adlt{ \delta }
\def\adltv#1{ \adlt_{#1} }
\def\adltoth{ \adltv{123} }
\def\adltothzn{ \adltoth\zdfmn }

\def\xO{ \mathcal{O} }
\def\xOv#1{ \xO_{#1} }
\def\xOot{ \xOv{12} }

\def\xfob{fibration object}

\def\xkk{ k }

\def\zsU{ 0_{\xcA} }

\def\btrno{ \ttrno }
\def\btrnt{ \ttrnt }
\def\btrnn{ \ttrnn }
\def\btrnv#1{ \ttrnv{#1} }

\def\ctrnv#1{ \btrnv{#1}\p }

\numberwithin{equation}{section}

\title[3D TFT and symplectic algebraic geometry II]{Three-dimensional
topological field theory and symplectic algebraic geometry II}
\author[A.~Kapustin]{Anton Kapustin}
\address{A.~Kapustin\\
California Institute of Technology, 452-48\\
Pasadena, CA 91125 }
\email{kapustin@theory.caltech.edu}
\thanks{The work of A.K. was supported in part by the DOE grant DE-FG03-92-ER40701}
\author[L.~Rozansky]{Lev Rozansky}
\address{
L.~Rozansky\\
Department of Mathematics\\
University of North Carolina at Chapel Hill\\
CB \# 3250, Phillips Hall\\
Chapel Hill, NC 27599
}
\email{rozansky@math.unc.edu}
\thanks{The work of L.R. was supported in part by the NSF grant DMS-0808974}

\dedicatory{To our parents}

\begin{document}
\maketitle
\begin{abstract}
Motivated by the path integral analysis\cx{KRS1} of boundary
conditions in a 3-dimensional topological sigma-model, we suggest
a definition of the 2-category $\ctLLX$ associated with a
holomorphic symplectic manifold $\xX$ and study its properties. The simplest objects of $\ctLLX$ are holomorphic lagrangian submanifolds $\yY\subset \xX$. We pay special attention to the case when $\xX$ is the total space of the cotangent bundle of a complex manifold $U$ or a deformation thereof. In the latter case the endomorphism category of the zero section is a monoidal category which is an $A_\infty$ deformation of the 2-periodic derived category of $U$. 
%
\end{abstract}
\tableofcontents
\section{Introduction}
\subsection{Sigma-models and categories}
\label{ss.intr}

Let $M$ be a real $\zN$-dimensional manifold and let
$\xcX=(\xX,\xs)$ be a pair in which $\xX$ is a real
manifold and $\xs$ is a geometric structure on $\xX$ such as a complex
structure or a symplectic structure.
A $\zN$-dimensional topological \sgmd\ (\tsgmd) with a
\emph{\xspt} (also known as the \xwsh\ or
the world-volume) $M$ and a
\xtsp\ $\xcX$
is a quantum field theory based on a
path integral over the infinite-dimensional space of maps $\MpsMX$. The measure on the
space of maps is determined by the structure $\xs$.

We are
interested in 2-dimensional and 3-dimensional \tsgmd s.

Path-integral based arguments
suggest that if manifolds $\xcX$ with a certain type of structure serve
as \xtsp s for $\zN$-dimensional \tsgmd\ then they form a
\zNc\ $\cC$ with special features. Let us briefly recall these features and
illustrate them by two well-known examples.

The category $\cC$ has a symmetric monoidal structure related to
the
cartesian product of manifolds:
\xlee{bag1.1}
\cC\times\cC\longrightarrow\cC,\qquad
(\xXo,\xso)\times(\xXt,\xst) = (\xXo\times\xXt,\xso\times\xst),
\xeee
where $\xso\times\xst$ is the natural structure on $\xXo\times\xXt$.
The monoidal structure has a unit element $\xcXopt=(\xXopt,\xsopt)$, where
$\xXopt$ is the manifold consisting of a single point and $\xsopt$
is the corresponding trivial structure. This element has the
property
\wlee{bag1.2}
\xcXopt\times\xcX=\xcX.
\weee
For a structured manifold $\xcX$ we define a \zNcmo\ of morphisms
\ylee{bag1.2a}
\cCX \edfn \HomC(\xcXopt,\xcX).
\yeee
This category is known in quantum field theory as the category of
boundary conditions of the \tsgmd\ related to $\xcX$.

The \zNc\ $\cC$ has a contravariant
duality functor
\xlee{bag1.3}
\cC\xmapta{\hve}\cC,\qquad
\zdlv{(\xX,\xs)} = (\xX,\zdlv{\xs} ),
\xeee
such that there is a canonical equivalence between \zNcmos\ of
morphisms:
\xlee{bag1.4}
\HomC(\xcXo,\xcXt) = \HomC(\xcXopt,\zdlv{\xcXo}\times\xcXt)
=\cCov{\zdlv{\xcXo}\times\xcXt}.
\xeee
This equivalence implies that an object
$\xOot\in\HomC(\xcXo,\xcXt)$ determines a functor between \zNcmos
\xlee{bag1.5}
\xPhv{\xOot}\colon\cCXo \longrightarrow\cCXt,
\xeee
which represents a composition of morphisms within $\cC$.
Moreover, a composition of morphisms of $\cC$ corresponds to the composition
of functors\rx{bag1.5}, so the structure of the \zNc\ $\cC$ is
determined by the boundary condition categories $\cCX$ and the
functors\rx{bag1.5}.

Recall two examples of this general construction for
$\zN=2$, that is, when $\cC$ is a 2-category. The first example is related to the A-model.
The structure $\xs$ is a symplectic structure (that is, $\xs$ is a symplectic form on $\xX$), the category of
boundary conditions $\cCX$ is the \FkFc\ $\rFuk(\xcX)$,
its simplest objects being lagrangian submanifolds of $\xX$,
the action of the duality functor is $\zdlv{\xs}=-\xs$, and the
functor $\xPhv{\xOot}$ is the lagrangian correspondence functor
determined by a lagrangian submanifold
$\xcEot\subset\zdlv{\xcXo}\times\xcXt$.

The second example of the 2-category $\cC$ comes from the B-model:
$\xX$ is a \fCY\ manifold,
$\xs$ is its complex structure, the category of boundary conditions
$\cCX$ is the bounded derived category of coherent sheaves
$\xrDv{\xcX}$, its simplest objects being complexes of holomorphic
vector bundles on $\xX$, the duality functor acts trivially:
$\zdlv{\xs}=\xs$, and the functor $\xPhv{\xOot}$ is the
Fourier-Mukai transform corresponding to the object $\xOot$.

\subsection{The 3-category of holomorphic symplectic manifolds}
\label{ss.hcatintr}

For $d=3$ a natural class of \tsgmd\ comes from the Rozansky-Witten model \cx{RW}.
In\cx{KRS1} we studied this \tsgmd\ and its 2-category of boundary conditions  from the path integral viewpoint. In this paper we
attempt to present a mathematical description of the 3-category $\ctLLL$  formed by
these theories and formulate conjectures about it.


Objects of $\ctLLL$ are \hlsmm s $\xcX=\Xsom$,
where $\xX$ is a complex manifold and
$\som\in\Omega^{2,0}$ is a \hlsm\ form: it is non-degenerate at every point of $\xX$
and $d \som=0$. If the symplectic form is canonical,
we abbreviate the notation $\Xsom$ down to $\xX$.
The monoidal structure\rx{bag1.1} comes from the product of
manifolds and the sum of their symplectic structures:
$\Xsomo\otimes\Xsomt = \brb{\xXo\times\xXt,\yepiuo(\somo) +
\yepiut(\somt) }$, where $\yepio$ and $\yepit$ are the projections
of $\xXo\times\xXt$ onto $\xXo$ and $\xXt$.
The duality
functor $\hve$ acts on objects by switching the sign of the
symplectic form: $\zdlv{\Xsom} = \Xsomm$.

The main purpose of the paper is to investigate the 2-category
$\cCX$ associated to a \hlsmm\ $\xcX=\Xsom$. We denote it
as
$\ctLLXsom$. The definition of $\ctLLXsom$ for a
general \hlsmm\ is rather complicated and requires a construction
of a \mlcs\ of 2-categories on $\xX$. Therefore we devote much of
this paper to special $\Xsom$ and present an attempt at a general definition
only in Section\rw{s.mcl}.
%

%

\subsection{Algebraic approach}

The first approach towards the description of
the 2-category $\ctLLXsom$
is based on the fact
that when $\xX$ is a
cotangent bundle of a complex manifold $\xcA$, the category
$\ctLLTsU$ can be described in terms of the properties of $\xcA$.
This description is algebraic in nature and there is no reference
to $\TsU$, so by looking at the definitions one would not see
directly that symplectic structure is involved.

In Section\rw{s.sct3} we study a `toy'  2-category
$\xcMFdbx$,
$\bax=\lvar{\ax}{n}$, which after a minor modification should be equivalent to the
2-category $\ctLLv{(\Ts\IC^n)}$
associated with a \zsaff\ space, that is, with the cotangent bundle
$\Ts\IC^n$.

Recall that for a polynomial $\xcW\in\ICbx$, an object of the category of matrix
factorizations $\xcMFbxW$ is a free finite rank \Ztgrdd\ $\ICbx$-module
$\xmM$ with a degree-1 endomorphism (called a curved differential) $\xdD$  satisfying the
condition $\xdD^2 = \xcW\,\xIdv{\xmM}$. The polynomial $\xcW$ is called \emph{a \crvng}.
A \crvng\ of a tensor product of two matrix factorizations over $\ICbx$ is
the sum of their \crvng s.

The simplest \xzobj s of $\xcMFdbx$ are polynomials $\xcW\in\ICbx$.
The category of morphisms between two polynomials $\xcWo$ and
$\xcWt$ is the
category of matrix factorizations of their
difference
\xlee{aeq1.7}
\Hom(\xcWo,\xcWt) =
\xcMFvv{\bax}{\xcWt-\xcWo},
\yeee
and the composition of morphisms comes from the tensor product of
matrix factorizations over $\ICbx$.


In Section\rw{s.sct4} we extend the algebraic construction of
$\ctLLv{(\Ts\IC^n)}$ to the cotangent bundle $\TsU$ of a complex
manifold $\xcA$ by defining algebraically a 2-category $\rDDprfU$
which is supposed to be equivalent to $\ctLLTsU$.
 Similar to the
\zaffn\ case, the simplest objects of $\rDDprfU$ are labeled
by holomorphic functions $\xcW$ on $\xcA$, and a category of
morphisms between two such functions is
the \xacrv\ version of the derived category of coherent sheaves
$\rDbU$:
\xlee{aeq1.8}
\Hom_{\rDDprfU}(\xcWo,\xcWt) = \rDprfvv{\xcA}{\xcWt-\xcWo}.
\xeee
This category is an analog of the category of matrix
factorizations\rx{aeq1.7} when the algebra $\ICbx$ is replaced by
the differential graded Dolbeault algebra $(\xbOmbU,\dlb)$.
A \xper\ object of $\rDprfUW$ is a \Ztgrdd\ vector bundle
$\xE\rightarrow\xcA$ with a \xacrv\ $(0,1)$ differential $\nbb$ such that
$\nbb^2 = \xcW\,\xIdv{\xE}$, and the composition of morphisms of the
type\rx{aeq1.8} comes from the tensor product of vector bundles.

\subsection{Geometric approach: the case of a cotangent bundle}

Path-integral analysis in\cx{KRS1} indicates that `in the classical approximation'
the category $\ctLLXsom$ should contain special `geometric'
objects. These objects are
holomorphic fibrations $\ycY\rightarrow\yY$, where
$\yY\subset\xX$ is a lagrangian submanifold.
Generally, fibration objects $\ycY\rightarrow\yY$ have to be
deformed because of quantum corrections, but this is unnecessary in two special
cases.
The first case is when $\ycY$ is a \opfib,
that is, the fiber is a point and the object is just the lagrangian submanifold $\yY$
itself. The second case is when $\xX$ is isomorphic to a cotangent bundle:
$\xX\cong\TsU$.

In Section\rw{s.sct6} we study
morphisms between \xfob s of $\ctLLXsom$ for
$\xX=\TsU$.
Here is an overview of our conjectures
for the simplest objects of $\ctLLTsU$ which are lagrangian
submanifolds.

%

We say that two holomorphic submanifolds $\yZo,\yZt\subset\xX$
have \emph{a \gdint} if  any point $x\in\yZo\cap\yZt$ has has an open
neighborhood $U_x\subset\xX$ which is isomorphic to a neighborhood
of $0$ in $\Tng_x\xX$,
so that $\yZo$ and $\yZt$ are identified with $\Tng_x\yZo$ and
$\Tng_x\yZt$. This condition implies that the intersection
$\yZo\cap\yZt$ is smooth.

Suppose that two lagrangian submanifolds $\yYo,\yYt\subset \xX=\TsU$ are
\fCY,  their intersection is \xgd,\ and
%
the difference of dimensions $\dim \yYo - \dim(\yYo\cap\yYt)$ is even. Then
the category of morphisms between them
becomes fairly simple:
\xlee{aeq1.2a}
\Hom_{\ctLLXsom}(\yYo,\yYt) = \rDprf(\yYo\cap\yYt).
\xeee
%
Moreover, if all intersections between three \fCY\ lagrangian
submanifolds $\yYo$, $\yYt$ and $\yYh$ are \xgd, then
%
%
%
the composition of morphisms $\cEot\in\Hom(\yYo,\yYt)$ and
$\cEth\in\Hom(\yYt,\yYh)$ is a combination of
pull-backs, tensor product and push-forward
\xlee{aeq1.3}
\cE_{23}\circ\cE_{12} = (\xiosoh)_\ast\Big(
\xiosot^\ast(\cE_{12})\otimes \xiosth^\ast(\cE_{23})
\Big),
\xeee
where $\xiosv{ij}$ are injections
\ylee{aeq1.4}
\xymatrix{ & \yYothi
\ar@{_{(}->}[dl]_{\xiosot}\ar@{^{(}->}[dr]^{\xiosth} \ar@{_{(}->}[d]^{\xiosoh}
\\
\yYoit &\yYoih& \yYtih.
}
\yeee
%
In particular, the endomorphism category of a lagrangian
submanifold $\yY\subset\xX$ is its 2-periodic category:
\xlee{bah7.1}
\End_{\ctLLXsom}(\yY) = \rDprf(\yY),
\xeee
and the composition of endomorphisms is just the
tensor product:
\xlee{bah7.2}
\cEo\circ\cEt=\cEo\otimes\cEt
\xeee
for any
$\cEo,\cEt\in\rDprf(\yY)$.

In Section\rw{ss.dlv} we study the geometric description of the
category $\ctLLXsom$ in case of a general \hlsmm\ $\Xsom$. We show
that general features remain the same as for the cotangent bundle
$\xX=\TsU$, but almost everything is
deformed. As we have mentioned, lagrangian submanifolds $\yY\subset\xX$ remain
objects of $\ctLLXsom$, but the categories of
morphisms\rx{aeq1.2a} and the composition rules\rx{aeq1.3} are
deformed. These \tAinf\ deformations are determined by the symplectic geometry
of tubular neighborhoods of the lagrangian
submanifolds involved, as summarized
in subsection\rw{ss.geomgen}.

\subsection{Relation between algebraic and geometric approaches}
\label{ss.relout}
In subsection\rw{ss.reltan} we relate the
algebraic and geometric descriptions of the
2-category of a cotangent bundle.
We describe the restriction of the equivalence functor
\xlee{bah5.1}
\etfe\colon \rDDprfaU
\xmapta{\cong}
\ctLLTsU
\xeee
(where the 2-category $\rDDprfaU$ is a slight modification of
$\rDDprfU$ defined in subsection\rw{ss.aug} )
to the objects of $\rDDprfaU$ whose images in $\ctLLTsU$ admit a
geometric description as holomorphic fibrations.

Let us sketch this relation for
$\xX=\Ts\IC^n$.
To a polynomial $\xcW\in\IC[\bax]$, which is an object of the
2-category $\xcMFdbx$, we associate the graph of its holomorphic
differential
\xlee{bag1.11}
\yYW = \{(\bfx,\bfp)\in\Ts\IC^n\,|\, \bfp =
\boldsymbol{\del}\xcW\}.
\xeee
$\yYW$ is a lagrangian submanifold of $\xX$ and it represents
the object of $\ctLLTsCn$ corresponding to $\xcW$.

Let $\CrW$ denote the critical locus of the polynomial $\xcW$:
$\CrW = \{ \bfx\in\IC^n\,|\,\del\xcW(\bfx)=0\}$.
We say that the critical locus is \xgd, if $\CrW$ is a smooth
manifold and the Hessian of $\xcW$ is non-degenerate in the normal
directions.
A category $\xcMFbxW$ `localizes' to $\CrW$, and if
$\CrW$ is \xgd then
we conjecture that $\xcMFbxW$ is equivalent to $\rDprf(\CrW)$ up
to a certain categorical `shift' explained in subsection\rw{ss.lcl}. The intersection
$\yYWo\cap\yYWt\subset\Ts\IC^n$ projects onto
$\CrWtmo\subset\IC^n$ and the projection establishes an
isomorphism between them. The intersection of $\yYWo$ and $\yYWt$ is \xgd\
exactly when the difference $\xcWt-\xcWo$
has a \xgd\ critical locus.
Hence in this case the categories
$\Hom_{\xcMFdbx}(\xcWo,\xcWt)$ and $\Hom_{\ctLLTsCn}(\yYWo,\yYWt)$
are equivalent (up to a shift).

\subsection{The 2-category of a deformed cotangent bundle}
\label{ss.tcdcb}
%

Path integral analysis in paper\cx{KRS1} suggests that similarly to
the derived category of coherent sheaves and in contrast to the Fukaya category, the 2-category
$\ctLLXsom$ is local: the category of morphisms between two
lagrangian submanifolds $\Hom_{\ctLLXsom}(\yYo,\yYt)$ is
determined by a small tubular neighborhood of the intersection
$\yYo\cap\yYt$ and the composition of morphisms connecting $\yYo$,
$\yYt$ and $\yYh$ is determined by a small tubular neighborhood of
the triple intersection $\yYo\cap\yYt\cap\yYh$. Therefore, if we
knew the properties of the 2-category $\ctLL$ of a tubular
neighborhood of a lagrangian submanifold $\yY\subset\xX$, we would
know the morphisms involving $\yY$ and their compositions.

In real symplectic geometry, a small tubular neighborhood of a
lagrangian submanifold is symplectomorphic to a small tubular
neighborhood of the zero section of its cotangent bundle. This is
no longer the case in holomorphic symplectic geometry: the holomorphic symplectic
structure of the cotangent
bundle $\TsY$ may have non-trivial deformations and a tubular
neighborhood of $\yY\subset\xX$ may be isomorphic to a tubular
neighborhood of the zero section within this deformed bundle. Thus
in order to apply the locality principle to the study of
$\ctLLXsom$, in Section\rw{ss.dlv} we explore the 2-category $\ctLL$ of a deformed
cotangent bundle of a complex manifold $\xcA$.

The best way for us to describe the 2-category
$\ctLLTsU$ is through its equivalence to the 2-category $\rDDprfU$.
Hence we describe the category $\ctLL$ of the deformed
cotangent bundle of $\xcA$ by constructing a deformation of
$\rDDprfU$. We assume that the deformation parameter $\hdf$ of
$\rDDprfU$ is the same parameter that describes the deformation of
the holomorphic complex structure of $\TsU$ and that the simplest
objects of the deformed category $\rDDprfUSk$ are functions $\xcW$
on $\xcA$, such that the graphs of their holomorphic differentials
$\del\xcW$ are lagrangian submanifolds of the deformed cotangent
bundle $\mtdfTUSk$.

We find that the deformation parameter $\hdf$ is an element of
$\xbOmov{\xcA,\Sb\TU}$, $\deg_{\Sb}\hdf\geq 2$, satisfying the \CMe\ $\dlb\hdf +
\shlf\Pbrv{\hdf,\hdf}=0$, where $\Pbrv{\xblnk,\xblnk}$ is the
Poisson-Schouten bracket on $\xbOmov{\xcA,\Sb\TU}$. The functions
$\xcW$ parameterizing the simplest objects of $\rDDprfUSk$ must
satisfy the equation $\dlb\xcW=\hdf(\del\xcW)$, where $\hdf$ is
regarded as a $(0,1)$ form on $\TsU$ taking values in
polynomial functions on the fibers of $\TsU$. The category of morphisms
between two objects $\xcWo$ and $\xcWt$ of $\rDDprfUSk$ turns out
to be an \tAinf-deformation of\rx{aeq1.8}:
\xlee{bah5.1a}
\Hom_{\rDDprfUSk}(\xcWo,\xcWt) =
 \rDsprfUdot,
\xeee
where the deformation parameter
$\xdfmot=\xcWt-\xcWo+\cdots\in\xbOmv{\xmi}(\xcA,\wedge^{\bullet}\TU)$
satisfies the \CMe\ $\dlb\xdfmot + \shlf\,[\xdfmot,\xdfmot] = 0$
and the bracket
$[\xblnk,\xblnk]$ is the Schouten bracket. The composition
of morphisms turns out to be a non-commutative deformation of the
tensor product that described the composition of morphisms within
$\rDDprfU$.

\subsection{Geometric approach: the general case}
\label{ss.geomgen}

The locality principle applied to the formula\rx{bah5.1a} says that
%
 for a general \hlsmm\ $\Xsom$
the category of morphisms between lagrangian submanifolds $\yYo$
and $\yYt$ having a \gdint is given by a deformation of
\ex{aeq1.2a}:
\xlee{bah5.3}
\Hom_{\ctLLXsom}(\yYo,\yYt) = \rDprf(\yYo\cap\yYt;\adfmot),
\xeee
where $\adfmot$ is a deformation parameter of a special form:
$\adfmot = \adfmotv{2}+\adfmotv{3}+\cdots$ and
$\adfmoti\in\xbOmv{\xmi}\brb{\xcA,\wedge^{\xmi}\Tng(\yYocct)}$.
Properties of $\rDDprfUSk$ suggest that if one of the exact
sequences $\Tng\yYi\rightarrow\Tng\xX|_{\yYi}\rightarrow\Nrm\yYi$
$i=1,2$
(for
example, the one with $i=1$)
 splits and the other lagrangian submanifold $\yYt$
can be presented as the graph of
$\del\xcW$ with $\xcA=\yYo$ as described in subsection\rw{ss.odc}
\ftnt{We expect that such a presentation exists if
the holomorphic bundle $\Tng\yYt|_{\yYocct}/\Tng(\yYocct)$ admits an $\OnC$
structure.},
then $\adfmot=0$, so that the formula\rx{aeq1.2a}
holds true.

The non-split nature of the exact sequence
$\TY\rightarrow\Tng\xX|_{\yY}\rightarrow\Nrm\yY$
is measured by a class
$\cdfbtY\in\Ext^1_{\yY}(\strsY,\Stat\TY)\subset\Ext^1_{\yY}(\NY,\TY)$, where
$\strsY$ is the structure sheaf of $\yY$ and we used the fact that
$\yY$ is lagrangian, so $\NY=\TsY$. If the exact sequence does not
split or equivalently $\cdfbtY\neq 0$,
%
then the category of endomorphisms of $\yY$ is
deformed:
\ylee{bah5.4}
\End_{\ctLLXsom}(\yY) = \rDprf(\yY;\xdfmv{\yY}),
\yeee
where
$\xdfmv{\yY} = \xdfmv{\yY,3} + \cdots$,
$\xdfmv{\yY,i}\in \xbOmv{\xmi}\brb{\xcA,\wedge^{\xmi}\TY}$
and the leading term $\xdfmv{\yY,3}$
is
quadratic in $\cdfbtY$ and linear in the Atiyah class $\ccrR$ of the tangent bundle
$\TY$ (\cf \ex{beq2.42b}).

To illustrate the deformation of the composition rule\rx{aeq1.3},
consider the case when $\yYo=\yYt=\yYh=\yY$
and $\cdfbtY\neq 0$. If the Atiyah class of the tangent bundle
$\ccrR$ is zero, then, according to subsection\rw{ss.gdgs}, the
endomorphism category of $\yY$ remains undeformed as in
\ex{bah7.1}, but the composition rule\rx{bah7.2} is deformed. For
example, if $\eEo,\eEt\in\rDprf(\yY)$ are two
holomorphic vector bundles
on $\yY$, then their composition is the deformed tensor product
\xlee{bah7.4}
\eEo\circ\eEt = (\eEo\otimes\eEt)_{\dfalot},
\xeee
where $\dfalot$ is a deformation parameter
$\dfalot=\dfaloto+\dfalott+\cdots$, $\dfaloti \in
\xbOmv{2i+1}\brb{\End(\eEo\otimes\eEt)}$, satisfying the \CMe\
$\dlb\dfalot + \shlf[\dfalot,\dfalot]=0$. The cohomology class
of the leading component of
$\dfalot$ is proportional to $\cdfbtY$ and to the Atiyah classes
of the bundles $\eEo$ and $\eEt$:
\ylee{bah7.5}
\check{\dfal}_{12,3}=
\cdfbtY\spsmb (\hacFEo\hacFEt).
\yeee
Note that $\check{\dfal}_{21,3}=-\check{\dfal}_{12,3}$, so the
composition in the category $\End_{\ctLLXsom}(\yY)$ is
non-commutative due to the deformation\rx{bah7.4}.

If $\cdfbtY=0$, then $\dfalot=0$ and the composition
rule\rx{bah7.2} remains undeformed, however the associator
isomorphism
\ylee{bah7.5a}
\zasoth\colon(\eEo\otimes\eEt)\otimes\eEh\longrightarrow
\eEo\otimes(\eEt\otimes\eEh)
\yeee
may be non-trivial.

\subsection{A \prshf\ definition of the 2-category $\ctLLXsom$}

In Section\rw{s.mcl} we sketch an approach to a rigorous
definition of the 2-category $\ctLLXsom$ as the category of global
sections of a certain \prshf\ $\prsfXsom$ of 2-categories defined
on $\xX$.

Let $\ICnbax$ be the affine space $\IC^\zn$ with standard
coordinates $\bax=\ax_1,\ldots,\ax_{\zn}$ and let
$\xcAbx\subset\ICnbax$ be an open subset inheriting the coordinates.
A \emph{\smrc} is a product of two such subsets $\UxVy$; it has a natural \hlsms\  $\som = \sum_{i=1}^{\zn}d\ay_i\wedge d\ax_i$.
A \emph{\rcch} is a symplectic embedding
\xlee{bai1.1}
\mch\colon\UxVy\rightarrow\xX.
\xeee
The images of these charts form an open cover of
$\xX$.

Since the coordinates $\bax$ on $\xcAbx$ provide a trivialization
of the cotangent bundle $\TsU$, there is a canonical symplectic embedding
$\UxVy\subset\TsU$.
To a \smrc\ $\UxVy$ we associate a 2-category $\ctLLbv{\UxVy}$ which is a
`\mlcz' of the 2-category $\rDDprfU$.
The
2-category $\ctLLbv{\UxVy}$ is a full subcategory of $\rDDprfU$
and
its simplest objects
are holomorphic functions $\xcW$ on $\xcA$
satisfying the condition that the associated lagrangian
submanifolds\rx{bag1.11} should lie within $\UxVy\subset\TsU$: for
any point $u\in\xcAbx$, the differential $v=\del_{\bax}\xcW\in\ICnbay$ should
belong to  $\xcVby\subset\ICnbay$.

Thus to a \rcch\rx{bai1.1} we associate the 2-category
$\rDDprfUxVy$. The structure of the \mlc\ \prshf\ $\prsfXsom$
comes from two types of functors defined in subsection\rw{ss.tctf}:
the restriction functor and the \Ldrt. A \smrc\ $\UxVy$ has
a lagrangian `\qfib' formed by subspaces $u\times\xcVby$, where $u\in\xcAbx$.
The restriction functor
$\rsfe\colon\rDDprfUxVy\longrightarrow\ctLLbv{\UxVypp}$ is
associated to a symplectic embedding
$\sfq\colon\UxVypp\hookrightarrow\UxVy$, which preserves the
\qfib. The \Ldrt\ is a special equivalence functor
$\dLtp\colon\rDDprfUxVy \longrightarrow \ctLLbv{\VyUmx}$, which
permutes the lagrangian fibrations $u\times\xcVby$, $u\in\xcAbx$ and
$\xcAbx\times v$, $v\in\xcVby$ of the \smrc\ $\UxVy$.

\subsection{Derived categorical sheaves}

In Section\rw{sec:sheaves} we discuss the relationship between
the RW model and the theory of derived categorical sheaves introduced by B.~Toen
and G.~Vezzosi \cite{ToVe}.  This relationship emerges when the target manifold $X$
is the cotangent bundle of a complex manifold $Y$. In this special case
one can promote the $\ZZ_2$ grading of the RW model to a $\ZZ$-grading
by declaring that natural fiber coordinates of the cotangent bundle
sit in cohomological degree $2$. Objects of the corresponding 2-category
of boundary conditions are naturally associated with sheaves of DG-categories
over $Y$, i.e. with derived categorical sheaves. More precisely,
objects of the kind mentioned above (complex fibrations over $Y$)
correspond to rather special sheaves of DG-categories. However,
we argue that more general sheaves of DG-categories, such as skyscraper
sheaves, can also be related to boundary conditions in the RW model if
we allow fibrations whose fibers are graded manifolds. Conjecturally,
the 2-category of boundary conditions in the $\ZZ$-graded RW model with
target $\TsY $ is a full sub-2-category of the 2-category of derived
categorical sheaves over $Y$. We perform some simple checks of this conjecture.

\subsection*{Acknowledgements}

L.R. is indebted to D.~Arinkin for many patient explanations of the
properties of coherent sheaves. He is also grateful to
V.~Ginzburg for numerous discussions and encouragement.
 A.K. would like to thank D.~Orlov
for the same. A.K. is also grateful to D.~Ben-Zvi, V.~Ostrik, and
L.~Positselski for advice. Both authors would like to thank Natalia Saulina for collaboration on Part I of the paper. The work of A.K. was supported in part by
the DOE grant DE-FG03-92-ER40701. The work of L.R. was supported by
the NSF grant DMS-0808974.

\section{The 3-category of affine spaces}
\label{s.sct3}


\subsection{\Ztpdcat\ of a \cdga}
\label{ztpdcat}

In this section we define the 2-category of boundary conditions corresponding to the RW model whose target is a complex symplectic vector spaces. We also describe the 3-category of all such RW models. Definitions of categories of morphisms between two \xzobjs\
of our 2-categories follow the same general pattern that we are
going to review in this subsection.
We follow closely the exposition of J.~Block\cx{JB}, replacing
\Zgrdng\ with \Ztgrdng\ when needed.

A commutative curved differential graded algebra (\CDGA) is a triple $\ycdgaA$, where $\xdlA$ is a
\Zgrdd\ associative commutative algebra
%
$\xdlA = \bigoplus_{i=0}^{\infty} \xdlAuv{i}$
%
with an associated \Ztgrdng\
\ee
\label{bae1.1a} 
\xdlA = \xdlAe\oplus\xdlAo,\quad
\xdlAe = \bigoplus_{i=0}^\infty \xdlAuv{2i},\quad
\xdlAo = \bigoplus_{i=1}^\infty \xdlAuv{2i+1},
\eee
$\nbb$ is its differential of (possibly inhomogeneous) odd degree not less than 1:
\ee
\xlabel{bae1.2}
\nbb^2=0,\qquad \nbb(\xdlAuv{i})\subset\bigoplus_{j=0}^\infty
\xdlAuv{i+2j+1},
\eee
and a \crvng\ $\xcW$ is a $\nbb$-closed element of $\xdlA$ of even \Ztdgr:
$\xcW\in\xdlAe$, $\;\;\nbb\xcW=0$.
%
%

We adopt the notations $\zdgz{\xblnk}$ and $\zdgt{\xblnk}$ for
\Zdgr\ and \Ztdgr\ respectively and we denote the elements of
$\Zt$ as $\yev$ and $\yod$.

A \Ztdgm\ (\ZtDGM) over a \CDGA\ $\ycdgaA$  is a pair $\cM=\xdgmM$,
where $\xmM$ is a \Ztgrdd\ module over $\xdlA$, while $\nbbM$ is
its curved differential: $\nbbM$ is a $\IC$-linear map
$\xmM\xrarv{\nbbM}\xmM$, $\zdgt{\nbbM}=\yod$,
satisfying the Leibnitz
identity
\ee
\xlabel{bae1.4}
\nbbM(am) = (\nbb a)\,m + (-1)^{\zdgt{a}}\, a\,(\nbbM m),\qquad
\any a\in\xdlA,\quad \any m\in\xmM,
\eee
and having the \crvng\ $\xcW$:
\ee
\xlabel{bae1.5}
\nbbM\circ\nbbM =\xcW\, \xId_{\xmM},
\eee
where $\xId_{\xmM}$ is the identity endomorphism of $\xmM$.
The module $\xmM$ can be rolled out into a \tpd\  twisted
complex
\begin{gather}
\xlabel{bae1.5a}
\xymatrix{
\cdots \ar[r]^-{\nbbMo} &
\xmMuv{\yev} \ar[r]^-{\nbbMe} &
\xmMuv{\yod} \ar[r]^-{\nbbMo} &
\xmMuv{\yev} \ar[r]^-{\nbbMe} &
\cdots
}
\\
\nonumber
\nbbMe\circ\nbbMo = \xcW \xId_{\xmMuv{\yod}},\qquad
\nbbMo\circ\nbbMe = \xcW \xId_{\xmMuv{\yev}},
\end{gather}
hence the name of the category.

Suppose that
two \Ztgrdd\ modules (or vector spaces) $M_1$ and $M_2$ have
endomorphisms $A_1$ and $A_2$ of a similar nature. Then for a
linear map $f\colon M_1\rightarrow M_2$ we use the commutator
notation for the following expression:
\ee
\label{bae1.5a1} 
[\cxrA,f] = A_2 f - (-1)^{\zdgt{A}\zdgt{f}} f A_1.
\eee

For two \ZtDGM s
$\cMo$ and $\cMt$
over a \CDGA\ $\ycdgaA$, the space
of homomorphisms $\Hom_{\xdlA}(\xmMo,\xmMt)$ has a differential $d$:
\ee
\xlabel{bae1.6}
df = [\nbbv{\cxrmM},f],\qquad f\in\Hom_{\xdlA}(\xmMo,\xmMt).
\eee
Thus \ZtDGM s are objects of a DG-category.

%
%


We define the tensor product of a \ZtDGM\ $\cMo$ over a \CDGA\ $\ycdgaAo$ and a \ZtDGM\
$\cMt$ over a \CDGA\ $\ycdgaAt$ as a \ZtDGM\ over a \CDGA\
$\ycdgaAot$ by the formula
\ee
\label{bae1.8} 
\cMo\otdlA\cMt = \big(\xmMo\otdlA\xmMt,\nbbMxo\otimes\xId_{\xmM} +
(-1)^{\zdgt{\xblnk}}\otimes\nbbMxt\big).
\eee
Also we define the dual \ZtDGM\ as
\ee
\label{bae1.9} 
\xdulv{\xdgmM} = \big(\xdulv{\xmM},\xdulv{\nbbM}\big).
\eee
Note that $\xdulv{\xdgmM}$ is a \ZtDGM\ over the \CDGA\ $\ycdgaAd$.

A \ZtDGM\ $\cP=\xdgmP$ is called \emph{\xper} if the $\xdlA$-module
$\xmP$ has the form
\ee
\label{bae1.10} 
\xmP = \chmP\otimes_{\xdlAZz} \xdlA,
\eee
where $\chmP$ is a projective \Ztgrdd\ module over $\xdlAZz$. The
\emph{\Ztpdcat} $\rDprfbA$ of a \CDGA\ $\ycdgaA$ is defined as a
graded category whose objects are \xper\ \ZtDGM s, and morphisms
are defined by
\ee
\label{bae1.7} 
\HmB{
\cPo,\cPt } = \rmH_d^\bullet\,\Big( \Hom_{\xdlA}(\xmPo,\xmPt)\Big).
\eee

One may enhance the category $\rDprfbA$ by
adding new `\xad' objects which are declared isomorphic to the
existing \xper\ objects according to the following rule. A \ZtDGM\
$\cM$
is called \emph{\xad}, if there
exists a \xper\ \ZtDGM\
$\cP$ such that for any \xper\ \ZtDGM\
$\cPp$
there is an isomorphism
%
\ee
\label{bae1.10a}
\Hom(\cPp,\cM) = \Hom(\cPp,\cP)
\eee
and for any other \xad\ \ZtDGM\ $\cMp$ we define
\ee
\xlabel{bae1.10a1}
\Ext(\cM,\cMp) \edfn \Hom(\cP,\cMp).
\eee

\subsection{Categories of matrix factorizations}

\subsubsection{Definition of the category}

A category of matrix factorizations is a particular case of a
\Ztpdcat\ defined in subsection\rw{ztpdcat}. For a finite set of
commuting variables
\ee
\label{bae1.11p} 
\bax=\lvar{\ax}{n}
\eee
consider the algebra of polynomial functions $\xdlA=\xdlAZz = \ICbx$ regarded as  $\ZZ$-graded CDGA placed in zero degree, with zero differential $\nbbA=0$, and the
\crvng\ being a polynomial $\xcW\in\ICbx$. Then the \cxmf\
$\xcMFbxW$ is the corresponding \Ztpdcat:
\ee
\label{bae1.11} 
\xcMFbxW \edfn \rDprf\big(\ICbx,0,\xcW\big).
\eee

According to the general definition, an object of $\xcMFbxW$ is a pair
$\cfM=\xmdM$, where $\xmM$ is a free \Ztgrdd\ $\ICbx$-module, while
$\xdDM$ is its curved differential:
\ee
\xlabel{bae1.12}
\xdDM\in\EndICbx(\xmM),\qquad\zdgt{\xdDM}=\yod,\qquad \xdDM^2 =
\xcW\,
\xIdv{\xmM}.
\eee
%
Morphisms between objects are defined by eq. \rx{bae1.7}. The
tensor product\rx{bae1.8} gives a functor
\ee
\label{bae1.12a} 
\xymatrix@C=1.5cm{\xcMFbxWo\times\xcMFbxWt \ar[r]^-{\otCbx} &
\xcMFvv{\bax}{\xcWo+\xcWt}}.
\eee

Let $\bay=\lvar{\ay}{k}$ be another list of variables, generally
of a different length. For $\xcWo\in\ICbx$ and $\xcWt\in\ICby$,
a \xmf\ $\cfMot\in\xcMFvv{\bax,\bay}{\xcWt-\xcWo}$
determines a functor
\ee
\label{bae1.12a1} 
\xymatrix@C=1.5cm{\xcMFbxWo \ar[r]^-{\xPhMot} &\xcMFvv{\bay}{\xcWt}}
\eee
which acts by taking a tensor product with $\cfMot$: for a \xmf\
$\cfM\in\xcMFbxWo$,
\ee
\label{bae1.12a2} 
\xPhMot(\xmM) = \cfM\otCbx\cfMot\in\xcMFvv{\bay}{\xcWt}.
\eee
Note that since $\xcWo$ cancels from the \crvng\ of the tensor
product\rx{bae1.12a2}, we can forget its $\ICbx$-module structure,
thus turning it into an object of $\xcMFvv{\bay}{\xcWt}$.

All categories $\xcMFbxW$ can be unified into a single 2-category
$\xcMFd$ of Landau-Ginzburg B-models with affine
\tgsp s along
the lines explained in subsection\rw{ss.intr}. This 2-category should be thought of as a 2-category of boundary conditions for the RW model whose target is a point.
An \xzobj\ of $\xcMFd$ is a pair $\cWx$,
$\xcW\in\ICbx$
or, equivalently, a category of matrix factorizations $\xcMFbxW$.
Morphisms between two \xzobj s also form a \cxmf
\def\Hmo#1{ \Hmb{#1} }
\ee
\label{bae1.12a3} 
\Hom_{\xcMFd}\brb{\cWxo,\cWyt}=\xcMFvv{\bax,\bay}{\xcWt-\xcWo},
\eee
and the composition of morphisms corresponds to the composition of
the functors\rx{bae1.12a1}: for two morphisms
%
$\cfMot\in\Hmo{ \cWxo,\cWyt }$ and
$\cfMth\in\Hmo{ \cWyt,\cWzh }$
%
we define
\ee
\xlabel{bae1.12a5}
\cfMth\circ\cfMot = \cfMth\otCby\cfMot.
\eee

\subsubsection{\tKmf s}

A Koszul complex corresponding to a list of polynomials
$\bap=\ap_1,\ldots,\ap_{\xkk}\in\ICbx$ is the tensor product of complexes
\ee
\label{bae1.22a} 
\xKobp=\bigotimes_{i=1}^{\xkk}{}
\Big(\xymatrix{
\ICbxo \ar[r]^{\ap_i}& \ICbxz }
\Big),
\eee
where $\ICbxv{i}$ denotes a rank-1 $\ICbx$-module of
\Ztdgr\ $i$. A Koszul \xmf\ is defined similarly:
for two sequences of polynomials $\bap,\baq\in\ICbx$ of equal length
$\xkk$, a Koszul \xmf\ $\xKmfbpq$ is a tensor product of
rank-$(1,1)$ \xmfs
\ee
\label{bae1.23} 
\xKmfbpq =
\begin{pmatrix}
\ap_1 & \aq_1
\\
\ap_2 & \aq_2
\\
\hdotsfor{2}
\\
\ap_k & \aq_k
\end{pmatrix}
=
\bigotimes_{i=1}^{\xkk}
\Big(\xymatrix{
\ICbxo \ar@<0.5ex>[r]^{\ap_i}& \ICbxz \ar@<0.5ex>[l]^{\aq_i}
}
\Big).
\eee
%
Obviously, $\xKmfbpq\in\xcMFvv{\bax}{\pdq}$, where
we use the notation
\ee
\xlabel{eab1.24a1}
\pdq = \sum_{i=1}^{\xkk} \ap_i\aq_i.
\eee

Suppose that an ideal $(\bap)\subset\ICbx$ is generated by a
regular sequence $\bap$ and $\xcW\in (\bap)$. Then the polynomial
$\xcW$ has a presentation $\xcW=\pdq$, where $\baq\in\ICbx$.
%
It is easy to check that for fixed $\xcW$, the isomorphism class of the
\xmf\rx{bae1.23} does not depend on the choice of the polynomials
$\baq$, so we can use an abbreviated notation
\ee
\label{e.koza} 
\xKpvv{\xcW}{\bap} = \xKmfbpq
\eee
%
for the \xmf\ $\xKmfbpq$ of \ex{bae1.23}.


The identity endofunctor of a \xmf\ category $\xcMFbxW$ can
be presented in the form\rx{bae1.12a1} with the help of a \tKmf.
%
%
For a list of
variables $\bax$, consider another list $\baxp$ of the same
length. The difference $\xcW(\baxp)-\xcW(\bax)$ belongs to the ideal
$(\baxp-\bax)\subset\ICv{\bax,\baxp}$, and the corresponding
\tKmf\
\ylee{e.koza1}
\IdbxW\edfn\xKpvv{\xcW(\baxp)-\xcW(\bax)}{\baxp-\bax}
\yeee
determines the
functor\rx{bae1.12a1}
\ee
\xlabel{e.kid}
\xymatrix
@C=2cm
{\xcMFbxW
\ar[r]^-{\xPhv{\IdbxW}}
&\xcMFbxpW},
\eee
which becomes the identity functor, if we identify the categories
$\xcMFbxW$ and $\xcMFbxpW$ by identifying the variables $\bax$ and
$\baxp$.

\subsubsection{\Knrp\ and the translation 2-functor}
\label{ss.ttran}

If the list of variables $\bax$ is empty and the  curving $\xcW$
is zero, then the corresponding category $\xcMFnz$
is equivalent to the category of \Ztgrdd\ vector spaces. The
latter has only two indecomposable objects: the 1-dimensional vector space
$\IC$ in degree 0 and its translation $\IC[\hat 1]$.

The category $\xcMFyot$ also has only two indecomposable objects:
the \tKmf\
\ylee{bag2.1a}
\cfMyot\edfn \xKpvv{\ayots}{\ayo-\sqrt{-1}\,\ayt}
\yeee
and its translation
$\cfMyot[\hat 1] = \xKpvv{\ayots}{\ayo+\sqrt{-1}\,\ayt}$.

The \Knrp\ theorem states that there is an equivalence of
categories
\ylee{bag2.1}
\xcMFnz \cong \xcMFyot,
\yeee
established by the functor
\ee
\xlabel{bag2.2}
\xymatrix@C=2.5cm{\xcMFnz
\ar[r]^-{\xPhv{\cfMyot}}
&\xcMFv{y_1,y_2;y_1^2 + y_2^2}}.
\eee
More generally, for any \crvngp\ $\xcW\in\ICbx$ there is an
equivalence of categories
\ee
\xlabel{bag2.3}
\xcMFbxyWt\cong\xcMFbxW,
\eee
established by the functor
\ee
\label{bag2.4}
\xymatrix@C=3.5cm{\xcMFbxW
\ar[r]^-{\xPhv{\IdbxW\otimes\cfMyot}}
&\xcMFbxpyWt}.
\eee

The translation 2-functor $\btrno\colon\xcMFd\rightarrow\xcMFd$ acts on a
matrix factorization category $\xcMFbxW$ by adding an extra variable
to the list $\bax$ and adding its square to the \crvngp\ $\xcW$:
\ee
\label{e.trn1}
\xcMFbxW\btrno \edfn \xcMFvv{\bax,\ay}{\xcW(\bax) + \ay^2}.
\eee
The action of $\btrno$ on categories of morphisms is provided by the
\Knrp\ functors\rx{bag2.4}. \Knrp\ also implies that the square of the translation
2-functor is equivalent to the identity 2-functor:
$$\btrnt\edfn \btrno\circ\btrno \simeq
\xIdv{\xcMFd}.$$

%
%
%
%


The definition\rx{e.trn1} of the translation 2-functor $\btrno$ can be
adapted to
the general setting of the \Ztpdcat\ of a \cdga\ (see
subsection\rw{ztpdcat}). Consider a special \CDGA\
$\xdlAx=(\IC[\shx],0,\shx^2)$, so that according to \ex{bae1.11}
$\rDprf(\xdlAx)=\xcMFvv{\shx}{\shx^2}$. For a \CDGA\ $\xdlA$ we
define the translation of its \Ztpdcat\ as
\ee
\label{bal1.1}
\rDprf(\xdlA)\btrno = \rDprf(\xdlA\otimes\xdlAx).
\eee
\Knrp\ establishes an equivalence of categories
\aee
\xlabel{e.knp3}
\rDprf(\xdlA)\btrnt \cong \rDprf(\xdlA).
\aeee

For a list of variables $\shbx=\lvar{\ay}{n}$
define a \CDGA\ $\xdlAbx=(\IC[\shbx],0,\shbx^2)$, so that by our definition
$\rDprfcA\btrnn=\rDprfcAAbx$.
%
The category $\rDprfcAAbx$ has an equivalent `intrinsic'
description in terms of objects and morphisms of $\rDprfcA$.
Namely, consider a category $\rDprfcAdn$, whose objects are pairs
$(\cP,\bafcP)$, where $\cP=\xdgmP$ is a \xper\ module of $\rDprfcA$,
while
$$\bafcP = \afcPv{1},\ldots\afcPv{n}\in\Exttod(\cP,\cP)$$
satisfies the property
$$\{\afcPv{i},\afcPv{j}\}= \delta_{ij}\,\xId_{\cP}.$$
Morphisms between two objects $(\cP,\bafcP)$ and
$(\cPp,\bafcPp)$ are $\rDprfcA$-morphisms
$\ag\in\Extub(\cP,\cPp)$ which
intertwine the lists $\bafcP$ and $\bafcPp$: $\ag\,\bafcP=\bafcPp\ag$.
An equivalence functor between the categories
$\rDprfcAdn$ and $\rDprfcAAbx$ maps an object $(\cP,\bafcP)$ of $\rDprfcAdn$
into an object $\big(\xmP\otimes\IC[\shbx],\nbbv{\xmP}+\xcd{\shbx}{\bafcP}
\big)$ of $\rDprfcAAbx$.

\subsection{The  \tcrca}

Fix a finite set of variables $\bax$ of length $n$. The 2-category of
relative curved differential graded (\yCDG) polynomial algebras
$\xcMFdbx$ is
a result of `fibering' the 2-category $\xcMFd$ over the algebra
$\ICbx$. One should regard this 2-category as a 2-category of boundary conditions for the RW model whose target is $\Ts \CC^n$. (These are not the most general boundary conditions: more general ones will be described in the next section.) Objects of $\xcMFdbx$
are pairs $\cWy$,
where
$\bay=\lvar{\ay}{k}$ is a list of `extra' variables of arbitrary length
and the \crvng\ $\xcW$ is an element of the algebra $\ICbxy$ over the ring $\ICbx$.
The category of morphisms between two \xzobj s is defined as
the \cxmf\ of the difference of \crvng s:
\ee
\label{bae1.13} 
\Hom_{\xcMFdbx}\brb{\cWyo,\cWzt} \edfn
\xcMFvv{\bax,\bay,\baz}{\xcWt(\bax,\baz)-\xcWo(\bax,\bay)}
\eee
(\cf \ex{bae1.12a3}).
The composition of morphisms between \xzobj s is given by the tensor
product functor\rx{bae1.12a}: for
$\cfMot\in\Hmo{\cWyo,\cWzt}$ and $\cfMth\in\Hmo{\cWzt,\cWuh}$
we define
\ee
\xlabel{bae1.15}
\cfMth\circ\cfMot =
\cfMth\otCbxz\cfMot\in\xcMFvv{\bax,\bay,\bau}{\xcWh-\xcWo}
\eee
(\cf \ex{bae1.12a2}).

The simplest objects of $\xcMFdbx$ are the ones without extra
variables, and we denote them as $\xcW \edfn \cWe$, where
$\xcW\in\ICbx$. The category of morphisms between such objects is
\ylee{bag3.0a}
\Hom_{\xcMFdbx}(\xcWo,\xcWt) = \xcMFvv{\bax}{\xcWt-\xcWo}.
\yeee

An object $\cWzot\in\xcMFdbxy$ determines a
\cfun
\ee
\label{bae1.16} 
\xymatrix@C=1.5cm{ \xcMFdbx \ar[r]^-{\dPhWzot} & \xcMFdby },
\eee
which turns an object $\cWu\in\xcMFdbx$ into an
object
\ee
\label{bae1.17}
\dPhWzot \cWu  = \xvprvv{\bax,\bau,\baz}{\xcW+\xcWot}.
\eee

The action of the \cfun\rx{bae1.16} on categories of morphisms
between \xzobj s is defined with the help of \tKmf s. For
\xzobj s $\cWuo,\cWwt\in\xcMFdbx$  we have to define the functor
%
%
\ee
\xlabel{bae1.18}
\xymatrix@C=2cm{
\Hmo{ \cWuo,\cWwt }
\ar[r]^-{\xPhuotWzot} &
\HmB{ \dPhWzot\cWuo,\dPhWzot\cWwt }
},
\eee
%
or more explicitly, according to \eex{bae1.13} and\rx{bae1.17},
\ee
\label{bae1.19} 
\xymatrix@C=2cm{
\xcMFvv{\bax,\bau,\baw}{\xcWt(\bax,\baw)-\xcWo(\bax,\bau)}
\ar[r]^-{
\xPhuotWzot
} &
\xcMFvv{\bay,\baxp,\baup,\bazp,\baxpp,\bawpp,\bazpp}
{\xcWttt}
},
\eee
where
\ee
\xlabel{bae1.20}
\xcWttt =
\xcWt(\baxpp,\bawpp)+ \xcWot(\baxpp,\bay,\bazpp)
-
\xcWo(\baxp,\baup) - \xcWot(\baxp,\bay,\bazp)
\eee
and for some lists of variables we used primed and double-primed
lists of the same length in order to change the names of
variables. The functor\rx{bae1.19} can be written in the form\rx{bae1.12a2}
\ylee{bae1.21}
\xPhuotWzot = \xPhv{\cfMot}
\yeee
%
%
for a certain matrix factorization
$\cfMot\in\xcMFvv{\xRtt}{\xcWttot},$
%
where
%
$$\xRtt =
\ICv{\bax,\bau,\baw,\bay,\baxp,\baup,\bazp,\baxpp,\bawpp,\bazpp}$$
%
and
\begin{multline}
\xlabel{bae1.22}
\xcWttot = \Big(\xcWt(\baxpp,\bawpp) - \xcWt(\bax,\baw)\Big)
-\Big(\xcWo(\baxp,\baup)-\xcWo(\bax,\bau)\Big)
\\
\nonumber
+ \Big(
\xcWot(\baxpp,\bay,\bazpp) - \xcWot(\baxp,\bay,\bazp)
\Big).
\end{multline}
The full formula for the \xmf\ $\cfMot$ is rather bulky, so we
describe it indirectly in terms of the \tKmf\rx{e.koza}.
Denote $$\bap =
(\baxpp-\bax,\bawpp-\baw,\baxp-\bax,\baup-\bau,\bazpp-\bazp).$$
Then $\xcWttot\in(\bap)$, and we set
\ee
\xlabel{bae1.27}
\cfMot=\xKpWttbp.
\eee

When the lists $\bax$ and $\bay$ have equal number of variables,
there exist two important equivalence \xxtf s of
type\rx{bae1.16}. The first one is the \xxtf
\xlee{bag3.1}
\Idba=\dPhv{\baa;\baa\cdot(\bay-\bax)}.
\xeee
If we identify $\xcMFdbx$ with
$\xcMFdby$ by identifying the variables $\bax$ and $\bay$, then
$\Idba$ becomes equivalent to the identity \xxtf. We leave the details of the
proof to the reader, while illustrating this statement by the
following example: after the identification of $\bax$ and $\bay$,
the action of $\Idba$ on an object $\ocW\in\xcMFdbx$ becomes
$\Idba(\ocW) = \xvprvv{\baxp,\baa}{\tcW}$, where $\tcW(\bax,\baxp,\baa)=\xcW(\baxp)+
\baa\cdot(\bax-\baxp)$, and the isomorphism between $\ocW$ and
$\Idba(\ocW)$ is established by the \tKmf\
$\xKpvv{\tcW-\xcW}{\baxp-\bax}$.

The second equivalence \xxtf\ is the `\Ldrt' and it has two
versions:
\xlee{bag3.1a}
\dLtp = \dPhv{\varnothing;\bax\cdot\bay},\qquad
\dLtm=\dPhv{\varnothing;-\bax\cdot\bay}.
\xeee
%
Again, we leave the verification of their
equivalence nature  to the reader. Note that the equivalence of morphism
categories $\Hom_{\xcMFdbx}\brb{\ocWo,\ocWt}$
and $\Hom_{\xcMFdby}\brb{\dLtpm(\ocWo),\dLtpm(\ocWt)}$ is a
corollary of the \Knrp.

The composition of two \Ldrt s with opposite signs is equivalent to the identity 2-functor:
\ylee{bag3.2}
\dLtp\circ\dLtm \simeq \dLtm\circ\dLtp \simeq \dId.
\yeee

The translation 2-functor is an equivalence 2-functor
$\ttrno\colon\xcMFdbx\rightarrow\xcMFdbx$, which acts on an object
$\cWy$ by adding a new variable $\aa$ to the extra variable list
and adding its square to the polynomial $\xcW$:
\xlee{bag2.1b}
\cWy\ttrno \edfn \xvprvv{\bay,\aa}{\xcW+\aa^2}.
\xeee
An equivalence
between the morphism categories
$\HmMFbx\brb{\xvprvv{\bay}{\xcWo},\xvprvv{\baz}{\xcWt}}$ and\\
$\HmMFbx\brb{\xvprvv{\bay}{\xcWo}\ttrno,\xvprvv{\baz}{\xcWt}\ttrno}$
is established by the \Knrp\
functor in view of an obvious
equivalence of categories
\ee
\label{bae1.27a6}
\Hom_{\xcMFdbx}\brB{
\xvprbvv{\bay}{\xcWo}\ttrno,\xvprbvv{\baz}{\xcWt}} =
\Hom_{\xcMFdbx}\brB{
\xvprbvv{\bay}{\xcWo},\xvprbvv{\baz}{\xcWt}}\btrno.
\eee
The \xtran\ by 2 is isomorphic to the identity endofunctor:
$\ttrnt\cong \xIdv{\xcMFdbx}$.

\subsection{The 3-category $\xcMFdd$ of polynomial algebras}

The 2-categories described in the previous subsection can be combined into a single 3-category. This 3-category should be thought of as the 3-category of RW models whose target spaces have the form $\Ts \CC^n$ for some nonnegative integer $n$.

An \xzobj\ of the 3-category $\xcMFdd$ is a list of variables
$\xlobx$.
The 2-category of morphisms between two \xzobj s $\xlobx,\xloby\in\xcMFdd$
is the correspondence 2-category $\xcMFdbxy$:
%
\ee
\xlabel{bae1.28}
\Hmb{\xlobx,\xloby} = \xcMFdbxy.
\eee
Each correspondence determines a 2-functor \rx{bae1.16}, and the
composition of correspondences as morphisms of $\xcMFdd$ is defined
to agree with the composition of the corresponding 2-functors. Namely, the
composition of two correspondences
%
$\xvprvv{\bau}{\xcWot}\in\xcMFdbxy$ and 
$\xvprvv{\baw}{\xcWth}\in\xcMFdbyz$
%
is the correspondence
\ee
\xlabel{bae1.30}
\xvprvv{\baw}{\xcWth}\circ\xvprvv{\bau}{\xcWot} =
\xvprvv{\bau,\baw,\bay}{\xcWot+\xcWth}\in\xcMFdbxz.
\eee

The identity endomorphism of an \xzobj\ $\xlobx$ can be represented by
the correspondence
\ee
\label{bae1.31} 
\xIdbx \simeq \xvprbvv{ \baa}{\xcd{\baa}{(\baxp-\bax)} 
}\in\xcMFdv{\bax,\baxp} = \End(\xlobx),
\eee
the lists $\baxp$ and $\baa$ having the same length as $\bax$
(\cf
\ex{bag3.1}).


The 3-category $\xcMFdd$ has a symmetric monoidal structure corresponding to
\ex{bag1.1}. The
product of \xzobj s corresponds to the concatenation of lists
\ee
\xlabel{bae1.31a}
\xlobx
\times
\xloby = \xlobxy
\eee
and the
unit element
is the empty list $\xlobe$. The duality endofunctor $\hve$ of
\ex{bag1.3} acts on $\xcMFdd$ as the identity.

\section{The 3-category of complex manifolds}
\label{s.sct4}
\label{ss.3catlv}

\subsection{The 2-periodic derived category of a curved complex manifold}
\label{tZtgdcs}

The \cxmf\rx{bae1.11} is based on a polynomial algebra $\ICbx$.
One can define a similar `analytic' category $\xcMFAbxW$ based on the algebra of holomorphic
functions on $\IC^n$.
The 2-periodic derived category of a curved complex manifold $\rDsprfUW$, which is the \Ztpdcat\ of its curved \Dlb\ algebra,
generalizes the category
$\xcMFAbxW$ from $\IC^n$ to a general complex manifold $\xcA$, so
that if $\xcA=\IC^n$, then there is an equivalence of categories
\ee
\xlabel{bae1.35e}
\rDsprfvv{\IC^n}{\xcW} \simeq \xcMFAbxW.
\eee
{}From the physical point of view $\rDsprfUW$ is the category of boundary conditions for the B-model with target $\xcA$ deformed by a curving $W$. In the special case when $W$ is a holomorphic function on $\xcA$, this theory is the topological Landau-Ginzburg model with target $\xcA$ and the superpotential $W$.

Following the general construction of subsection\rw{ztpdcat},
let us define the curved \Dlb\ algebra of
a complex manifold $\xcA$. The \Zgrdd\ \CDGA\ $\xdlA$ of \ex{bae1.1a}
is  the algebra of \xahfs\ $\xbOmbU$, the differential $\nbb$
is $\dlb$, and
%
the \crvng\ $\xcW$ is a
$\dlb$-closed even form:
%
\ee
\label{bae1.39} 
\nbb = \dlb,\qquad\xcW\in\xbOmeU,\qquad \dlb\xcW=0.
\eee

Let $\xE$ be a smooth (not necessarily holomorphic) \Ztgrdd\ vector
bundle over $\xcA$, and let $\xbOmbE$ be the space of \xahfs\
with values in $\xE$:
\ee
\xlabel{bae1.40}
\xbOmbE = \Gamma\big(\xE\otimes\wedge^\bullet\bTs\xX \big)
= \Gamma(\xE)\otimes_{\xbOmz(\xcA)}\xbOmbU,
\eee
where $\Gamma(\xE)$ is the space of sections of $\xE$.
The space $\xbOmbE$ is a \Ztgrdd\ module over $\xbOmbU$ of the
form\rx{bae1.10}.
%
A \emph{\cqhlmvb} is a pair $(\xE,\nbbE)$, where
$\nbbE$ is a curved \ydah-differential acting on $\xbOmbE$, that is,
$\nbbE$ is a $\IC$-linear operator 
\ee
\xlabel{bae1.41}
\nbbE\in\End_{\IC}\big(\xbOmbE\big),
\eee
which satisfies the following properties:
\begin{gather}
\label{bae1.42} 
\zdgt{\nbbE} = \yod,\\
\label{bae1.42a} 
\nbbE\,(\smu\wedge\sigma) = (\dlb\smu)\wedge \sigma +
(-1)^{\zdgt{\smu}}
\smu\wedge (\nbbE\,\sigma),
\\
\label{bae1.42a1} 
\nbbE^2\,\sigma = \xcW\wedge\sigma,
\end{gather}
where $\smu\in\xbOmbU$ and $\sigma\in\xbOmbE$.

We call the pair $(\xE,\nbbE)$ \qhlm, because even if $\xcW=0$, the bundle $\xE$ is not necessarily
holomorphic: if we split the differential $\nbbE$ according to the
Dolbeault degree:
\ylee{baj6.1}
\nbbE = \sum_{i=0}^{\dim\xcA} \nbbEai,\qquad
\nbbEai\colon \xbOmbE\longrightarrow\xbOmv{\bullet + i}(\xE),
\yeee
then $\nbbEao^2=-\{\nbbEaz,\nbbEat\}$ rather than $\nbbEao^2=0$,
so in general $\nbbEao$ does not determine a holomorphic structure on $\xE$.

If $(\xE,\nbbE)$ is a \cqhlmvb,
then the pair
\ee
\label{bae1.44} 
\cE=\xdgmE
\eee
is a \xper\ \ZtDGM\ over the \CDGA\ $\ycdgaBU$. In fact,
all \xper\ \ZtDGM s over this \CDGA\ originate in this way from vector bundles
with curved differentials.

The pair $\cmfUW$ will be called a \emph{\ccm}. We
define its 2-periodic derived category $\rDsprfUW$ as the
\Ztpdcat\ 
of its curved \Dlb\ algebra: 
\ee
\label{bae1.45} 
\rDsprfUW= \rDprf\ycdgaBU,
\eee
%
its \xper\ objects being the pairs\rx{bae1.44} and
morphisms defined according to the general formula\rx{bae1.7}.
The monoidal structure and the action of the duality functor
$\dulv{}$ also follow the general definitions\rx{bae1.8}
and\rx{bae1.9}.

For a \xper\ object\rx{bae1.44} we use an abbreviated
notation
\ee
\label{bae1.45b}
\cE=\adgmE.
\eee
If $\xcW=0$, then we will abbreviate the notation\rx{bae1.45} down
to
\ee
\label{bae1.45a} 
\rDprfU \edfn \rDprf(\xcA,0).
\eee
The latter category contains a full subcategory
which is equivalent to the bounded derived category of coherent
sheaves $\rDbU$:
%
\ee
\label{bae1.45cm1} 
\rDbU\hookrightarrow\rDprfU.
\eee
An object of $\xrDX$, represented by
a chain complex of holomorphic vector bundles
\ee
\label{bae1.45c}
\xymatrix{\xEv{0} \ar[r]^-{\sigma_0} &\xEv{1} \ar[r]^-{\sigma_1}
&\ldots \ar[r]^-{\sigma_{k-1}} & \xEv{k} },
\eee
corresponds to an object $(\xE,\dlb+\sigma)$ of $\rDprfU$, where $\xE$ is
the total
\Ztgrdd\ vector bundle
\ee
\xlabel{bae1.45c1}
\xE= \bigoplus_{i-\rm{even}} \xEv{i} \oplus\bigoplus_{i-\rm{odd}}
\xEv{i},
\eee
while $\dlb$ is the \ydah\ differential for holomorphic vector bundles and
$\sigma$ is the combined differential of the complex\rx{bae1.45c}:
%
$\sigma = \sum_{i=1}^k \sigma_i.$
%

The category $\rDsprfUW$ admits a certain `\Dlbf'. Consider
two \xper\ \ZtDGM s, which share the same vector bundle $\xE$:
$ \cE=\adgmE$ and $\cE\p=\adgmpE$.
We say
that the objects $\cE$ and $\cE\p$
are isomorphic up to \zord\ $k$, if the difference between their
connections is of higher degree as an element of $\xbOmbE$:
\ee
\xlabel{bae1.45a1} 
\nbbE\p-\nbbE=\dOk,\qquad\dOk\in\bigoplus_{i>k}\xbOmv{i}(\xE).
\eee

The tensor product of \xper\ \ZtDGM s over $\xbOmbU$ corresponds
to the tensor product of vector bundles:
\ee
\label{bae1.45a2} 
\cEo
\otimes
\cEt=\brb{
\xEo\otimes\xEt
,\nbbv{\xEo}+\nbbv{\xEt}}.
\eee
Since
\ee
\xlabel{bae1.45a3}
(\nbbv{\xEo} + \nbbv{\xEt})^2 = (\xcWo+\xcWt)\,\xIdv{\cEo\otimes\cEt},
\eee
the tensor product\rx{bae1.45a2} gives rise to a functor
\ee
\xlabel{bae1.45a4}
\rDsprfUWao \times \rDsprfUWat \xrarv{\otimes}
\rDsprfvv{\xcA}{\xcWo+\xcWt}.
\eee

For a holomorphic map $\mF\colon \xcAp\longrightarrow\xcA$ and for
a \xahf\ $\xcW$ of\rx{bae1.39} let $\mF^\ast(\xcW)$ denote its \pb\ to
$\xcAp$. We introduce a `derived' \pb\ functor $\mF^\ast$ and a
\pf\ functor $\mF_\ast$:
\ee
\xlabel{bae1.46}
\xymatrix@C=1cm{ \rDsprfBvv{\xcAp}{\mF^\ast(\xcW)}
\ar@<0.5ex>[r]^-{\mF_\ast}
& \;\rDsprfUW
\ar@<0.5ex>[l]^-{F^\ast}
}.
\eee
The \pb\ functor $\mF^\ast$ acts on \xper\ objects\rx{bae1.44} by pulling back
quasi-holomorphic vector bundles. The definition of the \pf\ functor $\mF_\ast$ is a bit
tricky. There is a \pb\ homomorphism of \CDGA s
\ee
\xlabel{bae1.47}
\xymatrix@C=1cm{
\Big(\xbOmbUp,\dlb,\mF^\ast(\xcW) \Big)& \ar[l]\ycdgaBU},
\eee
which turns a \xper\ \ZtDGM\ $\cE\p$ over the \CDGA\ $\Big(\xbOmbUp,\dlb,\mF^\ast(\xcW)
\Big)$ into a \ZtDGM\ over the \CDGA\ $\ycdgaBU$. We conjecture that the latter
\ZtDGM\ is \xad\ (\cf \ex{bae1.10a}) and use it as the definition of $\mF_\ast(\cE\p)$.

\subsection{The 2-category of \cfbs}
\label{ssctU}
The 2-category of \cfbs\ $\rDDprf$ is a generalization of the
2-category of \xanl\ matrix factorizations $\xcMFAd$. From the physical viewpoint, it is the 2-category of curved B-models, or equivalently the 2-category of boundary conditions for the RW model whose target is a point.
An object of $\rDDprf$ is a \ccm\ $\cmfUW$. A morphism between two
\ccm s $\cmfUWo$ and $\cmfUWt$ is a \xper\ (or an \xad) \ZtDGM
\ee
\xlabel{bae1.48}
\ycF\in \rDprf\cmfBv{ \xcAo\times\xcAt }{\;\yepiut(\xcWt)-\yepiuo(\xcWo)},
\eee
where $\yepio$ and $\yepit$ are projections
\ee
\label{bae1.49} 
\xymatrix@C=0.5pc{
& \xcAo\times\xcAt \ar[dl]_{\yepio} \ar[dr]^{\yepit}
\\
\xcAo &&\xcAt
}
\eee
Such an object $\ycF$ determines
a \FMtr\ functor
\ee
\xlabel{bae1.50}
\xymatrix@C=1.5cm{
\rDprfUWo \ar[r]^-{\xPhcF} & \rDprfUWt,
}
\eee
with the standard action on objects $\ycE\in\rDprfUWo$:
\ee
\label{bae1.51} 
\xPhcF(\ycE) = \yepidt\Big(\ycF\otimes\yepiuo(\ycE) \Big).
\eee
We define the composition of morphisms $\ycF$ as the
composition of the corresponding functors, so
for
%
$\ycFot\in\HmB{\cmfUWo,\cmfUWt}$ and
$\ycFth\in\HmB{\cmfUWt,\cmfUWh}$
%
the composition is
\ee
\label{bae1.53} 
\ycFth\circ\ycFot = \yepidoh \Big(
\yepiuot(\ycFot)\otimes\yepiuth(\ycFth)
\Big),
\eee
where the maps $\yepiij$ are the projections
\ee
\label{bae1.54} 
\xymatrix{ & \xcAo\times\xcAt\times\xcAh
\ar[dl]_{\yepiot}\ar[dr]^{\yepioh} \ar[d]^{\yepith}
\\
\xcAo\times\xcAt & \xcAt\times\xcAh &\xcAo\times\xcAh
}
\eee

\subsection{The 2-category of \xcfbs}
\label{xcfbs}
\def\xcA{ U }

Let us fix a complex manifold $\xcA$. A \hlfb\ over $\xcA$ is a
smooth \fbl
\ee
\label{bae1.55} 
\xfbrv{\Vfb}{\xcU}{\xcA}{\xcp}
\eee
where $\xcU$ is a complex manifold and the projection $\xcp$ is
holomorphic. The bundle\rx{bae1.55} does not have to be
holomorphic, so the complex structure of the fiber $\Vfb$ may be
different over different points of the base $\xcA$.

Let $\ycapU$ denote the fibered product of two fibrations with the
same base $\xcA$.
%
It has natural projections
\ee
\label{bae1.57} 
\xymatrix@C=0.5pc{
& \xcUo\ycapU\xcUt \ar[dl]_{\yepio} \ar[dr]^{\yepit}
\\
\xcUo &&\xcUt
}
\eee

The 2-category of \xcfbs\ $\rDDprfU$ is obtained by `fibering' the 2-category $\rDDprf$
over $\xcA$, the fibers $\Vfb$ being the analogs of the
manifolds $\xcA$
appearing in subsection\rw{ssctU}.
This category is also a generalization of the \xanl\ 2-category $\xcMFAdbx$:
\ee
\xlabel{bae1.57e}
\rDDprfv{\IC^n} = \xcMFAdbx.
\eee
{}From the physical viewpoint, $\rDDprfU$ is the 2-category of boundary conditions for the RW model with target $\sT U$.

An \xzobj\ of  
$\rDDprfU$ is a \xcfb\
$\cmfcUW$ where, according to the definition\rx{bae1.39},
$\xcW\in\xbOmecU$ and $\dlb\xcW=0$. The category of morphisms
between two \xzobj s is the \tpd\ derived category of their curved fibered
product:
\ee
\label{bae1.58} 
\HmDDU\brB{\cmfcUWo,\cmfcUWt} = \rDprf\cmfBv{ \xcUo\ycapU\xcUt
}{\;\yepiut(\xcWt)-\yepiuo(\xcWo)}.
\eee
An object $\ycF$ of the category\rx{bae1.58} generates a \FMtr\
functor
\ee
\xlabel{bae1.59}
\xymatrix@C=1.5cm{
\rDprfcUWo \ar[r]^-{\xPhcF} & \rDprfcUWt
}
\eee
which acts on an object $\ycEo\in\rDprfcUWo$ according to the
formula\rx{bae1.51}.
The composition of morphisms
%
$\ycFot\in\HmB{\cmfcUWo,\cmfcUWt}$ and
$\ycFth\in\HmB{\cmfcUWt,\cmfcUWh}$
%
is defined again
by the formula\rx{bae1.53},
where the maps $\yepiij$ are projections
\ee
\xlabel{bae1.61}
\xymatrix{ & \xcUo\ycapU\xcUt\ycapU\xcUh
\ar[dl]_{\yepiot}\ar[dr]^{\yepioh} \ar[d]^{\yepith}
\\
\xcUo\ycapU\xcUt & \xcUt\ycapU\xcUh &\xcUo\ycapU\xcUh
}
\eee
so that this composition agrees with the composition
of \FM\ functors.


The simplest \xzobj s of $\rDDprfU$ are \copfib s $\TrfUW$, where
\ee
\label{bae1.61b} 
\TrfU = 
\vcenter{
\xymatrix@C=1.5pc@R=1.5pc{
\{\text{1-point}\}\ar[r] & \xcA \ar[d]
\\
&\xcA
}
}
\eee
is a \opfib, that is, a fibration whose fiber is a single point,
and $\xcW$ is a holomorphic function on $\xcA$. We denote the
objects $\TrfUW$ simply as $\xcW$.
%
%
According to the general definition\rx{bae1.58}, the category of
morphisms between two \copfib s is the curved 2-periodic derived category
\ee
\label{bae1.63ap}
\Hom(\xcWo,\xcWt) = \rDsprfUWtmo,
\eee
and, in particular, the
endomorphisms of a \copfib\ $\xcW$ form the \tpd\ derived category of
$\xcA$:
\ee
\label{bae1.63a} 
\End(\xcW) = \rDprfU
\eee
containing the bounded derived category of coherent sheaves
 $\xrDX$ as a subcategory. The
composition of endomorphisms in $\End(\xcW)$ corresponds to the
tensor product\rx{bae1.45a2}, in other words, to the standard monoidal structure on $\rDprf(U)$.

The 2-category $\rDDprfU$ has a `\psmon' structure. For two objects
$\cmfcUWo$ and $\cmfcUWt$ we define
\ee
\xlabel{bae1.62}
\cmfcUWo\zcapU\cmfcUWt =
\Big(\xcUo\ycapU\xcUt,\yepiuo(\xcWo)+\yepiuo(\xcWt) \Big).
\eee
%
The
main property of the \psmon\ structure is that for the object $\TrfUZ$
consisting of the \opfib\ over $\xcA$ and $\xcW=0$, the following isomorphism holds:
\ee
\xlabel{bae1.63}
\HmB{\cmfcUWo\zcapU\cmfcUmWt,(\TrfU,0)} \cong
\HmB{\cmfcUWo,\cmfcUWt}.
\eee
%


The 2-category $\rDDprfU$ has a translation endo-2-functor
similar to\rx{bag2.1b}, which acts on a \zobj\ $\cmfcUW$ by
adding a space $\ICa$ (which is $\IC$ with the standard coordinate $\aa$) to
fibers of $\xcU$ and adding $\aa^2$ to the \crvng\ $\xcW$:
\ee
\xlabel{bae1.63a1}
\cmfcUW \ttrno = \cmfv{\xcU\times\ICa}{\xcW+\aa^2}.
\eee
%

For a holomorphic map $\mF\colon \xcAp\longrightarrow\xcA$
there exists a pull-back \xxtf
\ee
\xlabel{bae1.64}
\xymatrix@C=1cm{
\rDDprfUp & \rDDprfU \ar[l]-_{\dFua}
}
\eee
which acts on an object $\cmfcUW\in\rDDprfU$ by pulling back
the fibration $\xcU$ and the Dolbeault cohomology class $\xcW$:
\ee
\xlabel{bae1.65}
\dFua\cmfcUW =
\cmfBv{\mFua(\xcU)}{\mFua(\xcW)}
\eee
If $\xcAp$ is a \hlfb\ over $\xcA$ and the map $\mF$ is its projection,
then there is also a push-forward \xxtf
\ee
\xlabel{bae1.66}
\xymatrix@C=1.5cm{
\rDDprfUp \ar[r]^{\dFda}& \rDDprfU
}.
\eee
Indeed, for an object $\cmfcUWp\in\rDDprfUp$, a fibration
$\xcUp\xrarv{\xcp\p}\xcAp$ can be pushed forward
to a fibration $\mFda(\xcUp)$ over $\xcA$:
%
$\xcUp\xrarv{\mF\circ\xcp\p}\xcA$,
%
so we can keep
the form $\xcWp\in\xbOmecUp$
and  define
\ee
\xlabel{bae1.67}
\dFda\cmfcUWp = \Big( \mFda(\xcUp), \xcWp \Big).
\eee

\subsection{The 3-category of complex manifolds}
The 2-categories $\rDDprfU$ can be assembled  into a 3-category
$\rDDDprf$. It should be thought of as the 3-category of RW models with target-spaces of the form $\Ts\xcA$, for all complex manifolds $\xcA$.
Objects of $\rDDDprf$ are complex manifolds $\xcA$ (or,
equivalently, categories $\rDDprfU$). The category of morphisms
between two \xzobj s is
\ee
\label{bae1.68} 
\Hom(\xcAo,\xcAt) = \rDDprf(\xcAo\times\xcAt).
\eee
Objects of the category\rx{bae1.68} are called \emph{\corrs}. A
\corr\ $\cmfcUWot\in\Hom(\xcAo,\xcAt)$ determines a \xxtf
\ee
\label{bae1.69}
\xymatrix@C=2cm{
\rDDprf(\xcAo) \ar[r]^-{\xPhUWot} & \rDDprf(\xcAt)
}
\eee
acting on an object $\cmfcUWo\in\rDDprf(\xcAo)$ according to the
formula
\ee
\xlabel{bae1.70}
\xPhUWot\cmfcUWo = \yepidt\Big(\cmfcUWot\;\zcapv{(\xcAo\times\xcAt)}\;\yepiuo\cmfcUWo \Big),
\eee
where $\yepio$ and $\yepit$ are the projections\rx{bae1.49}.
The composition of correspondences
$\cmfcUWot\in\Hom(\xcAo,\xcAt)$ and
$\cmfcUWth\in\Hom(\xcAt,\xcAh)$
as morphisms in the 3-category $\rDDDprf$ is defined to
agree with the composition of their functors:
\ee
\xlabel{bae1.72}
\xcUth\circ\xcUot =
\yepidoh \Big(
\yepiuot\cmfcUWot\;\zcapv{(\xcAo\times\xcAt\times\xcAh)}\;\yepiuth\cmfcUWth
\Big),
\eee
where the maps $\yepiij$ are the projections\rx{bae1.54}.

The 3-category $\rDDDprf$ has a symmetric monoidal structure corresponding to that of
\ex{bag1.1}. The
product of \xzobj s corresponds to the product of the underlying complex manifolds
$\xcAo\times\xcAt$
and the unit
\xzobj\
is the complex manifold $\xcAopt$ consisting of a single point.
%
The duality endofunctor $\hve$ of \ex{bag1.3}
acts on $\rDDDprf$ as the identity.

\subsection{\xAug\ categories}
\label{ss.aug}

{}The curving $W$ enters the path-integral formulation of the RW model with boundaries only through it derivative $\partial W$. This suggests that one should define the 2-category of boundary conditions in such a way that $W$ is defined only up to addition of a locally constant function. Below we describe such a modification of the 2-category $\rDDprfaU$. It is also necessary for a geometric interpretation of $\rDDprfaU$ in terms of the cotangent bundle $\Ts U$, as we will see in the next section.

For an element
$\xcW\in\xbOmecU$, such that $\dlb\xcW=0$, define an \xaug\
category $\rDsprfaUW$ as a formal union over all locally constant
functions $\xcWlc$:
\ee
\xlabel{bae1.81a3}
\rDsprfaUW = \bigcup_{\del\xcWlc=0}
\rDsprfvv{\xcA}{\xcW+\xcWlc}.
\eee
We define the 2-category $\rDDprfaU$ in exactly the same way as
$\rDDprfU$, except that in the definition of morphisms\rx{bae1.58}
we replace the 2-periodic category $\rDprf$ with its \xaug\
version $\rDprfa$:
\ee
\label{bae1.58b} 
\HmDDaU\brB{\cmfcUWo,\cmfcUWt} = \rDprfa\cmfBv{ \xcUo\ycapU\xcUt
}{\;\yepiut(\xcWt)-\yepiuo(\xcWo)}.
\eee
Two objects $\cmfv{\xcU}{\xcWo}$ and $\cmfv{\xcU}{\xcWt}$ are
isomorphic within $\rDDprfaU$ if the difference $\xcWt-\xcWo$ is
locally constant on $\xcU$, that is, if $\del\xcWo=\del\xcWt$.

The \xaug\ 3-category $\rDDDprfa$ is defined in the same way as
$\rDDDprf$, except that we replace the categories $\rDDprf$
appearing in its definition with \xaug\ categories $\rDDprfa$.

The \xaug\ matrix factorization category is defined as the formal
union of categories
\ee
\xlabel{bae1.81a}
\xcMFabxW = \bigcup_{\shC\in\IC} \xcMFvv{\bax}{\xcW+\shC},
\eee
and the \xaug\ categories $\xcMFadbx$ and $\xcMFad$ are defined
similar to $\rDDprfaU$ and $\rDDDprfa$.

\def\xcWz{ \xcWv{0} }

\section{The 2-category $\ctLLTsU$: a geometric description and a
relation to $\rDDprfaU$}
\label{s.sct6}

%
\subsection{A geometric description of the 2-category
$\ctLLTsU$}
\label{ss.tcthlsmm}

Let $(X,\omega)$ be a holomorphic symplectic manifold. The RW model associates to $(X,\omega)$ a 2-category of boundary conditions $\ctLLXsom$.
Path integral arguments suggest that a certain part of $\ctLLXsom$ can be described in geometric terms.
In this subsection we consider the geometric description
when $\Xsom$ is the cotangent bundle of a complex manifold $\xcA$:
$\xX= \TsU$. The description uses only the \hlsm\
structure of $\xX$, and in our definitions we never refer to the cotangent
bundle structure. Conjecturally, this property should hold also for the whole 2-category $\ctLLTsU$, in the sense that it should be acted upon by the group of symplectic automorphisms of $\TsU$. This is far from obvious from the algebraic definition of $\ctLLTsU$ as the 2-category $\rDDprfaU$ given in the previous section.


\subsubsection{$\OnC$ bundles and matrix factorizations}

A holomorphic $\OnC$ vector bundle $\bnB$ over a complex manifold $\xcA$
determines  a `$\bnB$-twisted' version $\rDsprfUW\ytrnB$ of the
category $\rDsprfUW$.
%
The $\OnC$
structure determines, up to a non-zero constant factor,  a
holomorphic function $\xcWq$ on the total space of $\bnB$, which
is quadratic along the fibers. The
$\bnB$-twisting
replaces $\xcA$
with the total space of the bundle $\bnB$ and adds $\xcWq$ to the
\crvng:
\ee
\label{bae1.86a} 
\rDsprfUW\ytrnB = \rDsprfvv{\bnB}{\xcW+\xcWq}.
\eee
It is easy to see that the composition of \ytran s corresponds to
the sum of vector bundles:
$\ytrnv{\bnB_1}\ytrnv{\bnB_2}=\ytrnv{\bnB_1\oplus\bnB_2}$, and if
the bundle $\bnB$ is trivial, then
$
\rDsprfUW\ytrnB = \rDsprfUW\btrnv{\rnk\bnB},
$
where
$\btrno$
is the translation 2-functor\rx{bal1.1}.

A line bundle $\vcL\rightarrow\xcA$ has an $\rmO(1,\IC)$ structure
if and only if it is \sd, that is, if it is isomorphic to its
dual: $\vcL\cong\dulv{\vcL}$.
The top exterior power $\wdtpB$ of an $\OnC$ bundle $\bnB$ is \sd,
%
and there is a canonical equivalence of categories\footnote{We thank
M.~Kontsevich for pointing this out.}
\xlee{baf1.2}
\rDsprfUW\ytrnB = \rDsprfUW\ytrnv{\wdtpB}{\btrnv{\rnk\bnB-1}}.
\xeee
%


Similar to the untwisted \ytran\ discussed in
subsection\rw{ss.ttran}, the category
\rx{bae1.86a}
has an
alternative `intrinsic' description in terms of objects of $\rDsprfUW$. Consider the case when $\bnB$
is a  line bundle $\vcL$.
Let
$\ycL=\big( \xbOmbv{\vcL},\dlb\big)$ be the \xper\ \ZtDGM\
corresponding to $\vcL$.
The self-duality of $\vcL$
determines an isomorphism $\fvL$ between $\ycL\otimes\ycL$
and the structure sheaf $\strsU$:
\ee
\xlabel{bae1.87}
\fvL\in\Exttev\brb{\ycL\otimes\ycL,\strsU}.
\eee
An object of
$\rDsprfUW\ztrnL$
is a pair $\yobcE$, where $\ycE=\adgmE$
is a \xper\ \ZtDGM\rx{bae1.44} and the extension
$\fvcE\in\Exttod(\ycE,\ycE\otimes\ycL)$ satisfies the property
\ee
\xlabel{bae1.88}
\fvL\circ\fvcE\circ\fvcE = \xIdv{\ycE}.
\eee
Morphisms between two objects $\yobcE,\yobcEp\in\rDsprfmUWL$
are morphisms $g\in\Extub(\ycE,\ycEp)$ which intertwine the
extensions: $g\fvcE = \fvcEp g$.
We conjecture that the categories
$\rDsprfUW\ytrnL$
and $\rDsprfUW\ztrnL$
are equivalent.
%


\subsubsection{A holomorphic fibration with a lagrangian base as
a geometric object}

Let $\cnKU$ denote the canonical line bundle of a complex manifold
$\xcA$: $\cnKU\edfn \wdtpv{\TsU}$.
Let $\yY\subset\xX$ be a \hlgrsm, that is, $\yY$ is a
holomorphic submanifold of $\xX$, such that $\dim_{\IC} \yY =
\shlf\dim_{\IC}\xX$ and $\xom|_{\yY}=0$.
We are going to consider `geometric' objects of $\ctLLXsom$ which are pairs $\goYL$, where $\ycY$ is
a fibration
\ee
\label{bae1.90}
\xymatrix@C=1.5pc@R=1.5pc{
\yZ \ar[r] & \ycY \ar[d]^{\ypY} &
\\
& \yY \ar@{}[r]|<<<{\subset} & \xX
}
\eee
with a \hlgr\ base $\yY$, and $\vcLY$ is a holomorphic line bundle
on $\ycY$, whose square is the pull-back of the canonical bundle of $\yY$:
$\vcLY^{\otimes 2} = \ypY^\ast\cnKY.$

A particularly simple type of a holomorphic fibration\rx{bae1.90}
is a \opfib
\ee
\label{bae1.90b} 
\TrfY = 
\vcenter{
\xymatrix@C=1.5pc@R=1.5pc{
\yY \ar[d]
\\
\yY \ar@{}[r]|<<<{\subset}
&
\xX
}
}
\eee
Unless there is a danger of confusion, we denote such a fibration
simply as $\yY$. The pairs $\gYL$, where $\vLY\rightarrow\yY$ is a line bundle
such that $\yY\ott=\cnKY$,
are the simplest objects of the
type $\goYL$.


\subsubsection{Morphisms between geometric objects}

We say that two holomorphic submanifolds $\yYo,\yYt\subset\xX$
have a \emph{\gdint}, if any point $x\in\yYo\cap\yYt$ has an open
neighborhood $U_x$ which is isomorphic to a neighborhood of the
origin of the tangent space $\Tngx\xX$ in such a way that $\yYo$
and $\yYt$ correspond to their tangent spaces
$\Tngx\yYo,\Tngx\yYt\subset\Tngx\xX$.  This condition
guarantees that the intersection $\yYoct\edfn\yYo\cap\yYt$ is a disjoint union of holomorphic
submanifolds of $\xX$.

Define the $\xX$-product of two fibrations $\ycYo$ and $\ycYt$ as
a fibration over the intersection of their bases:
%
\ee
\xlabel{bae1.91}
\ycYo \ycapX \ycYt \edfn \ycYo|_{\yYoct} \ycapv{\yYoct}
\ycYt|_{\yYoct},\qquad
\xumap{\ycYo\ycapX \ycYt}{\xcpot}{\yYoct}.
\eee
There are obvious projections
\ylee{bae1.19b} 
\xymatrix@C=0.5pc{
& \ycYo \ycapX \ycYt \ar[dl]_{\pi_1} \ar[dr]^{\pi_2}
\\
\ycYo|_{\yYoct} &&\ycYt|_{\yYoct}
}
\yeee

Suppose that
the lagrangian bases $\yYo,\yYt\subset\xX$ of the fibrations of
two objects $\goYLo$, $\goYLt$
have a clean intersection, so the intersection $\yYot$ is a
complex manifold. The line bundle
\xlee{bae1.19b1}
\vcLot \edfn \pi_1^\ast(\vcLYo|_{\yYoct})\otimes
\pi_2^\ast(\vcLYt|_{\yYoct})
\otimes\xcpot^\ast(\cnKYot^{-1})
\xeee
is \sd. Indeed, for $i=1,2$ we have $\xcpot= \ypYi\circ\pi_i$, and since
$\vcLYi^2 = \ypYi^\ast\cnKYi$, the square of $\vcLot$ can be
presented as the pull-back of the product of canonical bundles:
\xlee{bae1.19b2}
\vcLot^2 =
\xcpot^\ast\brb{\cnKYo|_{\yYoct}\otimes\cnKYt|_{\yYoct}\otimes\cnKYot^{-2}}.
\xeee
Now consider
the quotient bundles
\xlee{bae1.90a} 
\vccLi = \Tng\yYi|_{\yYoct}/\Tng\yYoct,\qquad i = 1,2.
\xeee
The \hlsm\ form $\som$ produces a non-degenerate pairing between
$\vccLo$ and $\vccLt$, so their top exterior powers are dual to
each other and, as a result, the tensor product
$\wdtpv{\vccLo}\otimes\wdtpv{\vccLt}$ is a trivial line bundle. At the same
time,
$\cnKYi|_{\yYoct}=\wdtpv{\vccLdi}\otimes\cnKYot$, so
$\cnKYo|_{\yYoct}\otimes\cnKYt|_{\yYoct} = \cnKYot^2$ and the
tensor square\rx{bae1.19b2} is trivial, that is, $\vcLot$ is \sd.

Having established the self-duality of $\vcLot$, we propose
that the category of morphisms between the objects, whose
lagrangian bases have a \gdint, is the shifted 2-periodic derived category
of the $\xX$-product $\ycYo \ycapX \ycYt$:
\ee
\label{bae1.92} 
\!\!\!\!\!\!\!\!\!\Hom_{\ctLLX}\brB{\goYLo,\goYLt} =
\rDsprfcYoct\,\ytrnLot\,\btrnv{\shlf\dim \xX - \dim\yYot-1}.
\eee
Roughly speaking, the category of morphisms
$\Hom_{\ctLLX}\brB{\goYLo,\goYLt}$ is the 2-perio\-dic derived category of
coherent sheaves on the product $\ycYoct$. The origin of the shifts
in the \rhs of \ex{bae1.92} will become clear when we compare
\eex{bae1.92} and\rx{bae1.58}.

In the special case when $\ycYo$ and $\ycYt$ are \opfib
s\rx{bae1.90} with the same base $\yYo=\yYt=\yY$ and the accompanying line bundles
are the same, the formula\rx{bae1.92} becomes
\xlee{baj3.1}
\End_{\ctLLX}(\yY,\vcLY) = \rDprfY.
\xeee

\subsubsection{Composition of morphisms}

We describe the geometric composition of morphisms under the
simplifying assumption that
the lagrangian bases $\yYi$ of the fibrations $\ycYi$, $i=1,2,3$ are
\fCY\ and the accompanying line bundles $\vcLYi$ are trivial. This
implies that a clean intersection of two lagrangian submanifolds $\yYij=\yYi\cap\yYj$
is `semi-\fCY', that is, $\cnKYij\ott$ is trivial. Let us make a
stronger assumption that $\yYij$ is \fCY.
Then the general
formula\rx{bae1.92} simplifies:
\ee
\label{bae1.93}
\Hom(\ycYi,\ycYj) =
\rDprf\brb{\ycYi \ycapX \ycYj }
\eee
and we suggest that the composition of morphisms is a combination
of derived pull-backs, push-forwards and the tensor product: for two
morphisms $\cE_{12}\in\xMor(\ycY_1,\ycY_2)$ and
$\cE_{23}\in\xMor(\ycY_2,\ycY_3)$
their composition is
\ee
\label{bae1.97}
\cE_{23}\circ\cE_{12}=
(\yepi_{13})_\ast\Big(
\yepi_{12}^\ast(\cE_{12})\otimes \yepi_{23}^\ast(\cE_{23})
\Big),
\eee
where $\yepiv{ij}$ are the \xeps\
\ylee{bae1.98}
\xymatrix{ & \ycY_1\ycapX\ycY_2\ycapX\ycY_3
\ar[dl]_{\yepi_{12}}\ar[dr]^{\yepi_{23}} \ar[d]^{\yepi_{13}}
\\
\ycY_1\ycapX\ycY_2 &\ycY_1\ycapX\ycY_3& \ycY_2\ycapX\ycY_3
}
\yeee
%

In the special case when all $\ycYi$ are \opfib s\rx{bae1.90} with
the same base
$\yYo=\yYt=\yYh=\yY$, their categories of morphisms are given by
\ex{baj3.1} and the composition rule\rx{bae1.97} reduces to the
tensor product within $\rDprfY$: for $\cE,\cE\p\in\rDprfY$
\xlee{baj3.2}
\cE\circ\cE\p = \cE\otimes\cE\p.
\xeee

\subsection{Holomorphic lagrangian correspondences and
the  3-category of \hlsmm s}

In this subsection we describe part of the 3-category of RW models in geometric terms.
Throughout this subsection we will ignore the line bundles $\vcLY$
in the definition of the objects $\goYL$ of $\ctLLXsom$. To be
more precise, we may assume that all complex manifolds appearing
here are \fCY\ and all these bundles are trivial. Moreover, we
assume that all intersections are \xgd.

The statements of this subsection apply when $\Xsom$
are cotangent bundles: $\xX=\TsU$, but in the case of \opfib s the
statements apply to general \hlsmm s $\Xsom$.

\subsubsection{Lagrangian correspondence 2-functors}

Let us forget for a moment that the complex manifold $\xX$ has a
\hlsm\ structure and that the base $\yY$ of a fibration $\ycY$ must
be lagrangian. Then we may  define pull-back and
push-forward functors associated with a holomorphic map
$\mF\colon \xX\longrightarrow\xXp$. A pull-back of a fibration
$\ycYp\rightarrow\yYp\subset\xXp$ is a fibration $\mFua(\ycYp)\rightarrow \mF^{-1}(\yYp)$
constructed by pulling back $\ycYp$ by the restriction
$\mF|_{\mF^{-1}(\yYp)}$. In order to define a push-forward of a
fibration $\ycY\rightarrow\yY\subset\xX$ we assume that the
restricted map $\mF|_{\yY}\colon\yY\rightarrow\mF(\yY)$ represents a
holomorphic fibration. Then we define the fibration $\mFda(\ycY)$
as $\ycY\rightarrow\mF(\yY)$, whose projection is the
composition of projections
$\ycY\rightarrow\yY\rightarrow\mF(\yY)$.

A holomorphic fibration $\ycYaot\in\ctLLcrot$ determines a lagrangian correspondence 2-functor
\ee
\label{bae1.107}
\xymatrix@C=1.5cm{
\ctLLXsomo \ar[r]^-{\dPhYot} & \ctLLXsomt
}
\eee
defined by the formula
\ee
\label{bae1.108} 
\dPhYot= (\yepi_2)_\ast\Big(\ycYaot\ycapXott
\yepi_1^\ast
\Big),
\eee
where $\yepi_1$ and $\yepi_2$ are projections onto the factors
of the product $\xX_1\times\xX_2$:
\ee
\label{bae1.109} 
\xymatrix@C=0.5pc{
& \xX_1\times\xX_2 \ar[dl]_{\pi_1} \ar[dr]^{\pi_2}
\\
\xX_1 &&\xX_2
}
\eee
The action of the 2-functor\rx{bae1.108} on the bases of holomorphic
fibrations is described by a simple set-theoretic formula: if
$\ycYt=\xPhYot(\ycYo)$, then
\xlee{bah10.1}
\yYt = \yepit\brb{\yYot\cap\yepio^{-1}(\yYo)},
\xeee
where $\yYo$ and $\yYt$ are the bases of the corresponding
fibrations.

Although the operations $\yepi_1^\ast$ and
$(\yepi_2)_\ast$ do not correspond to well-defined 2-functors for 2-categories $\ctLL$,
their composition\rx{bae1.108} is \xwd, because if the base $\yYo$
of the fibration $\ycYo$ is lagrangian then so is the base $\yYt$
of its image determined by \ex{bah10.1}.


The \opfib\
\ee
\label{bae1.109b} 
\TrfDlX = 
\vcenter{
\xymatrix@C=1.5pc@R=1.5pc{
\DlX \ar[d]
\\
\DlX \ar@{}[r]|<<<{\subset}
&
\xX\times\xX
}
}
\eee
over the diagonal
$\DlX\subset\Xsomm\times\Xsom$ determines the identity endo-2-functor
$\dPhv{\TrfDlX}$ of the category $\ctLLXsom$.

\subsubsection{A geometric description of the 3-category $\ctLLL$}

As we mentioned in subsection\rw{ss.hcatintr}, objects of the
3-category $\ctLLL$ are \hlsmm s $\Xsom$. The duality functor $\hve$
switches the sign of the symplectic form: $\zdlv{\Xsom} = \Xsomm$,
and we define the 2-category of morphisms between two objects in
accordance with the general formula\rx{bag1.4}:
\ee
\xlabel{bae1.110}
\Hom_{\ctLLL}\brB{\Xsomo,\Xsomt } = \ctLLcrot.
\eee

The composition of morphisms represented by holomorphic fibrations
with lagrangian bases
%
is defined in such a way that it
would agree with the composition of their correspondence
functors\rx{bae1.108}: the composition of two morphisms
\ee
\xlabel{bae1.111}
\ycYaot\in\Hom\brB{\Xsomo,\Xsomt}
,\qquad\ycYath\in\Hom\brB{\Xsomt,\Xsomh}
\eee
is
\ee
\xlabel{bae1.112}
\ycYath\circ\ycYaot = (\yepi_{13})_\ast\Big(
\yepi_{12}^\ast(\ycYaot)\ycapv{\xXo\times\xXt\times\xXh} \yepi_{23}^\ast(\ycYath)
\Big),
\eee
where $\yepi_{ij}$ are the projections
\ee
\xlabel{bae1.113}
\xymatrix{ & \xXo\times\xXt\times\xXh
\ar[dl]_{\yepi_{12}}\ar[dr]^{\yepi_{23}} \ar[d]^{\yepi_{13}}
\\
\xXo\times\xXt &\xXo\times\xXh& \xXt\times\xXh
}
\eee

The identity endomorphism $\xIdXsom\in\End\Xsom$ is
the \opfib\rx{bae1.109b}:
\ee
\label{bae1.113a} 
\xIdXsom = \TrfDlX.
\eee

As we mentioned in subsection\rw{ss.intr},
the 3-category $\ctLLL$ has a symmetric monoidal structure:
%
\ee
\xlabel{bae1.114}
\Xsomo\times\Xsomt = \brb{\xXo\times\xXt,\yepiuo(\somo) +
\yepiut(\somt) },
\eee
where $\yepio$ and $\yepit$ are projections of the
diagram\rx{bae1.109}. The
unit object is the \hlsmm\ $\xXopt$ consisting of a single point.


\subsection{A relation between the geometric and the algebraic
descriptions}
\label{ss.reltan}

We have already outlined the equivalence of categories\rx{bah5.1}
%
%
in subsection\rw{ss.relout}.
Here we provide a more detailed description of the equivalence
\xxtf\ $\etfe$.
%
%
%

\subsubsection{Localization}
\label{ss.lcl}
Let us split the \crvng\ $\xcW$, which determines the category
$\rDsprfUW$, into zero degree and positive degree parts:
\ee
\label{bae1.89a} 
\xcW=\xcWz + \xcWpl,\qquad\xcWz\in\xbOmz(\xcA),\qquad
\xcWpl\in\bigoplus_{i\geq 1}\xbOmv{i}(\xcA).
\eee
We say that the set of critical points of $\xcWz$
\ee
\xlabel{bae1.89a1}
\xcAWz = \{ x\in\xcA\,|\, d\xcWz(x)=0\}
\eee
is \emph{\ygd} if it is
a smooth holomorphic submanifold of $\xcA$ and
the quadratic form induced by
the Hessian of
$\xcWz$ (that is, by the quadratic form of its second derivatives)
on the normal bundle $\Nrm\xcAWz$ is non-degenerate.
The non-degenerate Hessian gives rise to an $\OnC$ structure on
$\Nrm\xcAWz$,
and we conjecture the following equivalence of categories:
\ee
\label{bae.89a2}
\rDsprfUW =
\rDsprfmvv{ \xcAWz }{ \xcW|_{\xcAWz} }{\Nrm\xcAWz}.
\eee
In other words, we expect that the category $\rDsprfUW$
localizes to a small tubular neighborhood of $\xcAWz$ and that if
the dominant terms in the expansion of $\xcWz$ in the normal
directions to $\xcAWz$ are quadratic then lower degree terms do
not matter.

Since $\xcWz$ is locally constant on its critical set, the
relation\rx{bae.89a2} implies the following category equivalence:
\ee
\label{bae.89a2x}
\rDsprfaUW =
\rDsprfmvv{ \xcAWz }{ \xcWpl|_{\xcAWz} }{\Nrm\xcAWz}.
\eee
Note that we did not have to \xau\ the \rhs\ category, because
the connected parts of the critical set $\xcAWz$ will contribute
to it only when the constant value of $\xcW$ on them is zero.

In view of \ex{baf1.2}, the equivalence\rx{bae.89a2x} can be
simplified:
\xlee{bae.89a2y}
\rDsprfaUW = \rDsprfmvv{\xcAWz}{ \xcWpl|_{\xcAWz} }{
\wdtpv{\Nrm\xcAWz} }
\btrnv{\rnk\Nrm\xcAWz-1}.
\xeee
Finally, if $\xcWpl=0$, that is, if $\xcW$ is just a holomorphic function on $\xcA$,
then the equivalence takes the form
\ee
\label{bae.89a2x1} 
\rDsprfaUW =
\rDprf\lrbc{\xcAW}\ytrnv{\wdtpv{\Nrm\xcAW}}\btrnv{\rnk\Nrm\xcAW-1}.
\eee

\subsubsection{The relation between objects}
\label{ss.relobj}
Choose a
line bundle $\vcLz\rightarrow \xcA$ such that $\vcLz\ott=\cnKU$.
For simplicity we will consider only curved fibrations
$\cmfcUW\in\rDDprfaU$, for which $\xcWpl=0$ (see \ex{bae1.89a}),
that is, $\xcW$ is just a holomorphic function on $\xcU$. Recall that $\xcU$ is a
fibration\rx{bae1.55}: $\xcp\colon\xcU\rightarrow\xcA$. For a point $\xfu\in\xcU$ let
$\xfVu\subset\xcU$ be the fiber to which $\xfu$ belongs: $\xfVu=\xcp^{-1}(\xcp(\xfu))$.
Let
$\crcUW$ be the set of `fiber-critical' points of
$\xcW$: $\crcUW = \{\xfu\in\xcU\,\,|\,\,\,\del\xcW|_{\Tng_{\xfu}\xfVu} = 0
\}$. In other words, there is an exact sequence
$\Ts_{\xcp(\xfu)}\xcA\xrightarrow{a}\Ts_{\xfu}\xcU\xrightarrow{b}\Ts_{\xfu}\xfVu$
and $\xfu\in\crcUW$ if $b\big(\del\xcW(\xfu)\big) = 0$. The latter
condition means that
$\del\xcW(\xfu)$ is in the image of $a$. 
Hence
there is a map $\mpf\colon\crcUW\rightarrow\TsU$
such that $\mpf(\xfu) =
a^{-1}\big(\del\xcW(\xfu))\in\Ts_{\xcp(\xfu)}$.
The \emph{\spprt} of the object $\cmfcUW$ is defined to be the
image of $\mpf$, and we denote it suggestively
as $\yYUW$:
\xlee{bah8.2}
\yYUW=\xsupp\cmfcUW \edfn \mpf(\crcUW)\subset\TsU.
\xeee
%


Generally, $\yYUW$ is an isotropic submanifold with respect to the
symplectic structure $ \TsU$. Assume that for all $x\in\xcA$, the critical locus of the
function $\xcW|_{\xcp^{-1}(x)}$
is \xgd\ and that
$\crcUW$ is a complex manifold. Then $\yYUW\subset\TsU$ is a
lagrangian submanifold and the map $\mpf\colon\crcUW\rightarrow\yYUW$ is a
fibration. Let $\bnB_x$ denote the  normal bundle to the
critical set of $\xcW$ restricted to the fiber $\xcp^{-1}(x)$. This bundle
has an $\OnC$ structure given by the Hessian of $\xcW|_{\xcp^{-1}(x)}$, and all
these bundles together form a holomorphic $\OnC$ bundle $\bnB$ over
$\crcUW$.
Thus we define the action of the equivalence \xxtf\rx{bah5.1} on
the curved fibration $\cmfcUW$ as follows:
\ylee{bah5.2}
\etfe\cmfcUW = (\crcUW, (\xcp^\ast\vcLz)|_{\crcUW}\otimes\wdtpv{\bnB})\ttrnv{\rnk \bnB
-1},
\yeee
that is, the pair in the \rhs of this equation
%
is the object of $\ctLLTsU$ corresponding to the object
$\cmfcUW\in\rDDprfaU$.

%


The geometric object corresponding to the pair $\cmfcUW$ is particularly simple, if $\xcU$ is a \opfib\ over $\xcA$. Then
$\xcW$ is just a holomorphic function on $\xcA$ and the object of
$\ctLLTsU$ corresponding to $\xcW$ is the pair $(\yYW,\zzf^\ast\vcLz)$, where $\yYW$ is the
graph of $\del\xcW$:
%
\xlee{bah1.1}
\yYW = \{ p\in\Ts_x\xcA\,|\,x\in\xcA,\; p = \del\xcW|_x \}
\xeee
and
$\zzf\colon\yYW\rightarrow\xcA$ is the restriction of the
projection $\TsU\rightarrow\xcA$ to $\yYW$ (it establishes the
isomorphism between $\yYW$ and $\xcA$ as complex manifolds).

\subsubsection{The relation between categories of morphisms}

We will compare the categories of morphisms within 2-categories $\rDDprfaU$ and
$\ctLLTsU$ only for the simplest objects. Let $\xcWo$ and
$\xcWt$ be holomorphic functions on $\xcA$ such that their
difference $\xcWot=\xcWt-\xcWo$ has a \xgd\ set of critical
points. This is equivalent to saying that $\yYWo$ and $\yYWt$ have a \gdint.

According to the definition\rx{bae1.63ap} and the
equivalence\rx{bae.89a2x1}, the category of morphisms within
$\rDDprfaU$ is
\xlee{baf2.1}
\Hom_{\rDDprfaU}(\xcWo,\xcWt)   = 
\rDsprfaUWot =
\rDprf\lrbc{\xcAWot}\ytrnv{\Nrm\xcAWot}\btrnv{\rnk\Nrm\xcAWot-1}.
\xeee
At the same time, according to \ex{bae1.92},
\xlee{baf2.2}
\Hom_{\ctLLTsU}\brB{\yobYWo,\yobYWt} = \rDsprfYot\ytrnLot
\btrnv{\dim\xcA - \dim\yYot - 1},
\xeee
where $\yYot \edfn \yYWo\cap\yYWt$, the maps
$\zzfi\colon\yYWi\rightarrow\xcA$ are the restrictions of the
projection $\TsU\rightarrow\xcA$ and
\xlee{baf2.3}
\vcLot = (\zzf_1^\ast\vcLz)|_{\yYot}\otimes(\zzf_2^\ast\vcLz)|_{\yYot}\otimes
\cnKYot^{-1}.
\xeee
%

The maps $\zzfo$ and $\zzft$, as well another projection restriction
$\zzfot\colon \yYot\rightarrow \xcAWot$, establish isomorphisms
between the corresponding complex manifolds.
Therefore, there is an equivalence of categories
$$\rDsprfYot\ytrnLot =
\rDsprfv{\xcAWot}\ytrnv{\zzfv{12,\ast}\vcLot},$$ and according to
\ex{baf2.3}, the push-forward of the line bundle $\vcLot$ is
$$\zzfv{12,\ast}\vcLot = \vcLz^2\otimes\cnKYot^{-1}= \cnKU\otimes\cnKYot^{-1}
=\wdtpv{\Nrm\xcAWot}.$$
Since
$$ \rnk\Nrm\xcAWot = \dim\xcA - \dim\xcAWot,\qquad \xcAWot=  \yYot,$$
we established the equivalence of categories\rx{baf2.1}
and\rx{baf2.2} provided by the \xxtf\ $\etfe$.

\subsubsection{The relation between 2-functors}

Let $\som$ be the canonical symplectic form of the cotangent
bundle $\TsU$. The symplectomorphism $\xtsm\colon
(\TsU,\som)\rightarrow(\TsU,-\som)=\zdlv{(\TsU)}$ reverses the cotangent
vectors: for $p\in\Ts_{q}\xcA$ we define $\xtsm(p)\edfn-p$.
For
two cotangent spaces we define the symplectomorphism
$\xtsmo\colon(\TsUo,\somo)\times(\TsUt,\somt)\rightarrow
(\TsUo,-\somo)\times(\TsUt,\somt)$ which reverses the orientation
of the first cotangent vector: $\xtsmo\edfn\xtsm\times \xId$.

The symplectomorphism $\xtsm$ acts as a \xxtf\
$\ddxtsm\colon\ctLLTsU\rightarrow\ctLL\zdlv{(\TsU)}$ by transforming the bases
of holomorphic fibrations $\ycY\rightarrow\yY\subset\TsU$.
Similarly, $\xtsmo$ acts as a \xxtf\ $\ddxtsmo\colon\ctLL(\TsUo\times\TsUt)\rightarrow
\ctLL\brb{\zdlv{(\TsUo)}\times\TsUt}$. Now consider a composition
of 2-functors
\ylee{bah7.1}
\ddxtsmo\circ\etfeot\colon\rDDprfav{\TsUo\times\TsUt}\longrightarrow
\ctLL\brb{\zdlv{(\TsUo)}\times\TsUt}.
\yeee
We leave it for the reader to check that the 2-functors\rx{bae1.69} and\rx{bae1.107}
coming from the objects related by\rx{bah7.1} are intertwined by
\xxtf s $\etfe$, that is, the diagram
\ylee{bah7.2}
\xymatrix@C=3cm@R=1.2cm{
\rDDprf(\xcAo) \ar[d]^{\etfeo} \ar[r]^-{\xPhUWot} & \rDDprf(\xcAt)
\ar[d]^{\etfet}
\\
\ctLLv{(\TsUo)} \ar[r]^-{\dPhYot} & \ctLLv{(\TsUt)}
}
\yeee
is commutative, if $\ycYaot = \ddxtsmo\circ\etfeot \cmfcUW$. The
easiest part of this commutativity is the verification that
the support of a \xcfb  from $\rDDprf(\xcAo)$ is transformed as in
\ex{bah10.1}:
\xlee{bah7.3}
\yYv{\xPhUWot\cmfcUWo} =
\yepit\brb{
\yYv{\cmfcUWot}\cap
\yepio^{-1}\brb{\yYv{\cmfcUWo}}
},
\xeee
where $\yepio$ and $\yepit$ are the projections
\ylee{bah7.4}
\xymatrix@C=0.5pc{
& \TsUo\times\TsUt \ar[dl]_{\pi_1} \ar[dr]^{\pi_2}
\\
\TsUo &&\TsUt
}
\yeee
%



\section{The 2-category of a \zdctngb}
\label{s.sct5}
\label{ss.dlv}
\label{ss.3catdfrm}

\subsection{Outline}
\label{ss.outline}

The geometric description of the 2-category $\ctLLXsom$ provided
in subsection\rw{ss.tcthlsmm} is completely correct only when
$\Xsom$ is a cotangent bundle: $\xX=\TsU$. However, path-integral
considerations\cx{KRS1} suggest that the main feature of that
description holds true for a general \hlsmm\ $\Xsom$: a pair $\gYL$, where
$\yY\subset \xX$ is a lagrangian submanifold and
the
line bundle $\vLY$ is a square root
of the
canonical bundle of $\yY$, always represents an object in
$\ctLLXsom$.

Path integral arguments\cx{KRS1} also suggest that the 2-category
$\ctLLXsom$ is local. From the physics perspective, locality
means that there are no instanton corrections to the path
integral, that is, there are no 3-dimensional analogs of
A-model holomorphic disks which lie at the heart of the Floer homology and the
Fukaya category. From the mathematical perspective, locality means
that the category of morphisms between two objects of $\ctLLXsom$
is determined by the structure of $\xX$ in the formal neighborhood of
the intersection of their \spprt s, that is,
$\Hom_{\ctLLXsom}\brb{\gYLo,\gYLt}$ is determined by the formal neighborhood of $\yYo\cap\yYt\subset \xX$ (when $\xX$ is
a cotangent bundle and the intersection of supports is \xgd,
the category of morphisms is determined just by the intersection itself, as
suggested by \ex{bae1.92}). The composition of
morphisms is determined by the structure of $\xX$ near the triple
intersection of supports (\cf \ex{bae1.97}).

Let $\yY\subset\xX$ be a lagrangian submanifold of a general
\hlsmm\ $\Xsom$. Locality of $\ctLLXsom$ means that if we knew the
structure of the 2-category $\ctLL$ of the formal neighborhood of
$\yY$, we would know exactly all categories of morphisms involving
the objects
$\gYL$
as well as the compositions of
such morphisms.

In real symplectic geometry, a sufficiently small tubular
neighborhood of a lagrangian submanifold $\yY$
is symplectomorphic to a tubular neighborhood
of the zero section of the cotangent bundle $\TsY$. However, in
holomorphic case this is no longer true: there may be non-trivial
deformations of the holomorphic symplectic structure of (the formal
neighborhood of the zero section) of the cotangent bundle $\TsY$ and the
formal neighborhood of $\yY\subset\xX$ may be isomorphic to a
deformed formal neighborhood of the zero section of $\TsY$. Therefore, in
order to gain information about the morphisms involving the
object
$\gYL$
of $\ctLLXsom$
we have to explore the 2-category
corresponding to a deformed cotangent bundle $\mtdfTYSk$, where $\hdf$ is a
deformation parameter of the holomorphic symplectic structure of $\TsY$.

To understand the 2-category $\mtdfTYSk$ we follow the
algebraic approach:
for a complex manifold $\xcA$, which plays the role of $\yY$,
we construct a deformation $\rDDprfUSk$ of the 2-category
$\rDDprfU$. The construction of this deformation is based on two
path-integral-motivated assumptions: the deformation parameter $\hdf$ of the 2-category is
the same as the deformation parameter of the holomorphic symplectic structure
of $\TsU$ and the simplest objects of the deformed category
$\rDDprfUSk$ (that is, the objects corresponding to \opfib s and
described by holomorphic functions $\xcW$ on $\xcA$ in the
category $\rDDprfU$)
should be related to the lagrangian submanifolds of the deformed
cotangent bundle $\mtdfTUSk$ in the same way as in the undeformed
case discussed in subsection\rw{ss.reltan}.

\subsection{Differential Lie-Gerstenhaber algebras and the \CMe}

Let us review some well-known facts about algebras
governing the deformations of objects and categories appearing in
this paper.

\subsubsection{General definitions}

A \DPGa\ $\xdlL$ is a \Ztgrdd\  vector space endowed with
a differential $\pgD$ and a compatible graded Lie bracket $\pgbdd$  which may be even or odd. Let $\pgdg$ be the \Ztdgr\ of the Lie
bracket. If $\pgdg=\yev$ then $\xdlL$ is called a differential Lie algebra. If in addition $\xdlL$ has a supercommutative associative product compatible both with $\pgD$ and the Lie bracket, then $\xdlL$ is called a differential Poisson algebra. If $\pgdg=\yod$ and $\xdlL$ has a supercommutative associative product compatible both with $\pgD$ and $\pgbdd$, then $\xdlL$ is called a \DGa.

The graded Lie bracket of $\xdlL$ descends to its
$\pgD$-homology $\rmH_{\pgD}(\xdlL)$.

If an element $\pga\in\xdlL$ of \Ztdgr\ $\pgdg+1$ satisfies the
\CMe\
\xlee{bcea2.1}
\pgD\pga + \shlf\,\pgbvv{\pga}{\pga} = 0,
\xeee
then the operator
\xlee{bcea2.1a}
\pgDa = \pgD + \pgbvv{\pga}{\cdot}
\xeee
is also a
differential and it determines a deformed \DPGa\ $\xdlLa$. If two
\CMlt s are related by a `gauge transformation' $\pgb$, whose
infinitesimal form is
\xlee{bcea2.2}
\delta\pga = \pgDa\pgb,
\qquad \pgb\in\xdlL,\quad\dgZt\pgb=\pgdg,
\xeee
then the corresponding deformed algebras are isomorphic.

If a \CMlt\ is presented as a formal power series in a parameter
$\dfe$:
\ylee{bcea2.2a}
\pgae =
\sum_{i=1}^\infty \pgai\,\dfe^i,
\yeee
then the leading coefficient
$\pgao$ is $\pgD$-closed and, due to the gauge symmetry,
the corresponding deformation of $\xdlL$ is
determined by its class $\cpgao\in\rmH_{\pgD}(\xdlL)$. The $\dfe^2$ part
of the \CMe\rx{bcea2.1} says that
\ylee{bcea2.3}
\pgD\pgat + \shlf\,\pgbvv{\pgao}{\pgao} = 0,
\yeee
so $\cpgao$ satisfies the condition
\xlee{bcea2.4}
\pgbvv{\cpgao}{\cpgao}= 0. 
\xeee

The importance of \DPGa s stems from the fact that they
appear as deformation complexes of objects in categories and Hochschild complexes of categories and 2-categories.
Equivalence classes of \CMlt s parameterize deformations
of those objects, categories, and 2-categories.
We need three particular examples
:
the differential Lie algebra$\zdrcEnb$ which governs deformations of a \cqhlmvb\ $\adgmE$,
the \DGa\
$\zdgaU$ governing \tAinf\ deformations of the
2-periodic derived category
$\rDprfU$ of a complex manifold $\xcA$, and
the differential Poisson algebra $\zdpaXsom$ which, according to
path integral considerations\cx{KRS1}, governs deformations of the
2-category $\ctLLXsom$.



\subsubsection{The \DLa\ of a \cqhlmvb}

Let the pair $\adgmE$, where $\nbbE^2=\xcW\xIdv{\xE}$, be a \cqhlmvb\
defined in subsection\rw{tZtgdcs}.
The corresponding \darc\ $\zdrcEnb$
is the algebra of Dolbeault forms $\xbOmbzeE$ (with total grading),
the differential is the covariant Dolbeault differential $\nbbE$, and the
\xLie\ bracket is the supercommutator:
$\Lbrv{\zdfmo,\zdfmt} \edfn \zdfmo\zdfmt - (-1)^{\zdgt{\zdfmo}\zdgt{\zdfmt} }
\zdfmt\zdfmo$.

A \CMlt\ $\zdfm$ determines a deformed \qhlmvb\
$$\admgEzdfm\edfn(\xE,\nbbE+\zdfm)$$ (in the \rhs of this formula $\zdfm$
denotes an odd bundle map
$\zdfm\colon\xbOmbE\rightarrow\xbOmbE$).

\subsubsection{The \DGa\ of a complex manifold}
The well-known \DGa\ $\zdgaU$ of a complex
manifold $\xcA$ is the algebra $\xbOmbUUWT$, the
\Ztgrdng\ coming from the total degree of forms and
wedge-powers, with the differential $\dlb$ and the \NSb\
$\Lbrv{\cdot,\cdot}$. In fact, $\zdgaU$ has \Zgrdng, and its
degree 2 \CMlt s parameterize deformations of the derived category of
coherent sheaves $\rDbU$. More generally, even \CMlt s of $\zdgaU$
parameterize \tAinf\ deformations of $\rDprfU$.

A holomorphic function $\xcW\in\xbOmz(\xcA)$, $\dlb\xcW=0$
obviously satisfies the \CMe\ and hence determines a deformation of $\zdgaU$
which we denote as $\zdgalWU$. It has a new differential $\dlbW=\dlb+\Lbrv{\xcW,\cdot}$
(\cf \ex{bcea2.1a}). The corresponding deformation of $\rDprfU$ is
the curved category $\rDsprfUW$.

Define a \emph{\xrl} grading on $\zdgaU$ as
\ylee{beq2.1a1}
\dgrl\xbOmv{\ydgn}(\xcA,\wedge^{\ydgm}\TU) = \ydgn-\ydgm,
\yeee
then, obviously,
\ylee{beq2.1a2}
\dgrl \dlb = \dgrl [\cdot,\cdot ] = 1.
\yeee
We say that a \CMlt\ $\xdfm$ is \emph{\xrlnn} if $\dgrl\xdfm\geq
0$ and \emph{\xrlb} if $\dgrl\xdfm=0$.

\subsubsection{The \DPa\ of a \hlsmm}
\label{ss.dpahlsmm}
The
\DPa\ $\zdpaXsom$ of a \hlsmm\ $\Xsom$
is defined as the algebra $\xbOmbX$
of $(0,\bullet)$-forms on $\xX$ with the differential $\dlb$ and
with the Poisson bracket coming from $\som$. If $\xX=\TsU$, where $\xcA$ is a complex manifold,
then we may consider a simpler version of this algebra, which we
denote as $\zdpaTU$:
it is the algebra of $\SbTU$-valued $(0,\bullet)$ forms
$\xbOmbUUST$, where $\Sb$ is the symmetric algebra, and its
differential is $\dlb$. There is a natural injection
\xlee{bcea2.5}
\xbOmbUUST \hookrightarrow \xbOmbUST,
\xeee
which turns an element of $\xbOmbUUST$ into a $(0,\bullet)$ differential form
on $\TsU$ having a polynomial dependence on fiber coordinates and
restricting to zero on all fibers.
The bracket of $\zdpaTU$ is a well-defined restriction of the
Poisson bracket of $\xbOmbUST$, so
the injection\rx{bcea2.5} becomes an injection of \DPa s
\xlee{bcea2.6}
\zdpaTU\hookrightarrow \zdpaTsU.
\xeee

\subsubsection{\Ainfalg\ and its modules}

Let us recall the main definitions. An \Ainfalg\ is a \Ztgrdd\
vector space $\aA$ endowed with a series of $\xmn$-linear maps (\xmnmlt
s) $\xmmbf=\xmmz,\xmmo,\ldots$,
\ylee{beq2.1b}
\xmmn\colon \aA^{\otimes \xmn}\rightarrow \aA,\qquad \dgZt\xmmn = \xmn,
\qquad \xmn=0,1,2,\ldots,
\xeee
satisfying the relations
\xlee{beq2.2b}
\sum_{\substack{m_1,m_2,\xmn\geq 0 \\ m_1+m_2+\xmn=N}}
(-1)^{m_1 + \xmn m_2}\,
\xmmv{m_1+m_2+1}(\xId^{\otimes m_1}\otimes\xmmn\otimes\xId^{\otimes
m_2})=0
\xeee
for all $N\geq 0$.
Among other things, these relations indicate
that if $\xmmz=0$, then $\xmmo$ is a differential ($\xmmo^2=0$) and it satisfies the
usual Leibnitz rule with respect to the multiplication $\xmmt$.

If $\xmmz\neq 0$ then the \Ainfalg\ $\aA$ is called \emph{\xwk} or \emph{\xacrv}.

A module over an \Ainfalg\ $\aA$ is a vector space $\xmM$ endowed
with a series of $\xmn$-linear maps (\xmnact s) $\zmmbfM=\zmmMo,\zmmMt,\ldots$
\ylee{beq2.3b}
\zmmMn\colon \xmM\otimes\aA^{\otimes (\xmn-1)}\rightarrow \aA,\qquad \dgZt\zmmMn = \xmn,
\qquad \xmn=1,2,\ldots,
\xeee
satisfying the relations
\begin{multline}
\label{beq2.4b}
\sum_{\substack{m_2,\xmn\geq 0 \\ m_2+\xmn=N}}
(-1)^{\xmn m_2}\,\zmmMv{m_2+1}(\zmmMn\otimes\xId^{\otimes
m_2})
\\
+
\sum_{\substack{m_1\geq 1\\m_2,\xmn\geq 0 \\ m_1+m_2+\xmn=N}}
(-1)^{m_1 + \xmn m_2}\,
\zmmMv{m_1+m_2+1}(\xId^{\otimes m_1}\otimes\xmmn\otimes\xId^{\otimes
m_2})=0
\end{multline}
for all $N\geq 0$.

For two $\aA$-modules $\xmMo$ and $\xmMt$, the vector space
$\tHom(\xmMo,\xmMt)$ is formed by sequences $\zmbf=(\zmfo,\zmft,\ldots) $,
$\zmfn$ being $\xmn$-linear maps
\ylee{beq2.5b}
\zmfn\colon \xmMo\otimes \aA^{\otimes(\xmn-1)}\rightarrow\xmMt.
\yeee
Define a differential $d$ acting on $\tHom(\xmMo,\xmMt)$ by the
following
formula for each term in $d\zmbf$:
\begin{multline}
\label{beq2.6b}
(d\zmbf)_{N} = \sum_{\substack{m_2,\xmn\geq 0 \\ m_2+\xmn=N}}
(-1)^{\xmn m_2}\,
\zmmuMt_{m_2+1} (\zmfn\otimes\xId^{\otimes m_2} )
\\
- \sum_{\substack{m_1\geq 1\\m_2,\xmn\geq 0 \\ m_1+m_2+\xmn=N}}
(-1)^{m_1 + \xmn m_2}\,
\zmfv{m_1+m_2+1}(\xId^{\otimes m_1}\otimes\xmmn\otimes\xId^{\otimes
m_2})
\\
-
\sum_{\substack{m_2,\xmn\geq 0 \\ m_2+\xmn=N}}
(-1)^{\xmn m_2}\,
\zmfv{m_2+1} (\zmmuMo_{\xmn}\otimes\xId^{\otimes m_2} ).
\end{multline}

\subsubsection{A \fDAinfalg\ and the perfect homotopy category of its modules}
Now we adapt the general definitions to the Dolbeault setting.

First of all, note that the Dolbeault algebra
$(\xbOmbU,\dlb)$ associated with a complex manifold $\xcA$ has a
canonical sequence of \xmnmlt s
$\xmmbfuz$
which turn it into an
\Ainfalg: $\xmmuzo=\dlb$, $\xmmuzt$ is the wedge-product, and $\xmmuzn=0$ for $\xmn\neq 1,2$.

We define \emph{a \fDAinfalg} (\FDAi) on a complex manifold $\xcA$ as
the space $\xbOmbU$ endowed with \xmnmlt s $\xmmbf$ satisfying the
relation\rx{beq2.2b} and the following restriction on Dolbeault
degrees of their deviation from the canonical \xmnmlt s:
\xlee{beq2.7b}
\dgDlbv{(\xmmn-\xmmuzn)}\geq\xmn,\qquad \any\xmn\geq 0.
\xeee

If $\xdfm\in\xbOmbUUWT$ is a \xrlnn\ \CMlt, then the corresponding
deformation
$\ycdgaBUl$
of the Dolbeault algebra is a
\fDAinfalg, the relation between the components of $\xdfm$ and
\xmnmlt s being fairly complicated. If the deformation parameter
$\xdfm$ is \xrlb, that is,
\ylee{beq2.8b}
\xdfm = \sum_{\xmn=0}^\infty\xdfm_{\xmn},\qquad
\xdfm_{\xmn}\in\xbOmv{\xmn}(\xcA,\wedge^{\xmn}\TU),
\yeee
then the dominant part of the deviations $\xmmn-\xmmuzn$ is
determined by the formula
\xlee{beq2.9b}
(\xmmn-\xmmuzn)(\smu_1,\ldots,\smu_{\xmn})
= \xdfm_{\xmn}\spsmb(\del\smu_1,\ldots,\del\smu_{\xmn}) + \cdots,
\yeee
where $\smu_1,\ldots,\smu_{\xmn}\in\xbOmbU$ and the correction
terms have higher Dolbeault degree than the first term in the
\rhs of this equation.

A simple example of a \xrlb\ deformation of the Dolbeault algebra
$(\xbOmbU,\dlb)$ is $\xdfm = \xcW$, where $\xcW$ is a holomorphic
function on $\xcA$. Then the formula\rx{beq2.9b} has no correction terms and
the deformation results in the \xacrv\
Dolbeault algebra
$\ycdgaBU$ discussed already in subsection\rw{tZtgdcs}.

Consider again the Dolbeault algebra $(\xbOmbU,\dlb)$ as a
\fDAinfalg\ with \xmnmlt s $\xmmbfuz$. Its
\xper\ \ZtDGM\rx{bae1.44} with $\nbbE^2=0$ has a canonical
structure of an \tAinf-module if we endow it with
\xmnact s $\xmmbfuzE$ such that $\xmmuzEo=\nbbE$, $\xmmuzEt$ is the
standard multiplication and $\xmmuzEn=0$ for $n>2$. We define a
\xprmi\ over a \fDAinfalg\ as the vector space $\xbOmbE$ endowed
with a sequence of \xmnact s $\zmmbfE$ satisfying the
relations\rx{beq2.4b} and the restriction
\ylee{beq2.10b}
\dgDlbv{(\zmmEn-\xmmuzEn)}\geq\xmn,\qquad \any\xmn\geq 1.
\yeee

Finally, we define the 2-periodic perfect derived category of a \fDAinfalg\ $(\xbOmbU;\xmmbf)$:
its objects are \xper\ \tAinf-modules $(\xbOmbE;\zmmbfE)$
and morphisms between two modules are homologies of the
differential\rx{beq2.6b}.

If $\xdfm\in\xbOmbUUWT$ is a \xrlnn\ \CMlt, then $\rDsprfUd$
denotes the 2-periodic perfect derived category corresponding to the deformed \tAinf-algebra
$\ycdgaBUl$. In other words, $\rDsprfUd$ is the result of deforming $\rDsprfU$
with $\xdfm$. In particular, if $\xdfm=\xcW$, where $\xcW$ is a
holomorphic function on $\xcA$, then $\rDsprfUW$ is the
category\rx{bae1.45}.

\subsection{Deformations of \hlsms s}

\subsubsection{The general case}

Let $\Xsom$ be a \hlsmm. Deformations of its \hlsms\
which preserve the \dR\ cohomology class of $\som$ are
parameterized up to gauge equivalence by \CMlt s of the \DPa\
$\zdpaXsom$ defined in subsection\rw{ss.dpahlsmm}. Namely, if an
element
\ylee{beq2.11b}
\hdf\in\xbOmov{\xX}\subset\zdgav{\xX}
\yeee
satisfies the \CMe
\xlee{bbea2.3}
\label{bbea2.2a}
\dlb\hdf + \shlf\,\Pbrv{\hdf,\hdf} = 0,
\xeee
where $\Pbrv{\ ,\ }$ is the Poisson bracket corresponding
to the symplectic form $\som$, then the corresponding deformation
of the complex structure of $\xX$ is described by the
Beltrami differential
\xlee{bbea2.2a1}
\bdmu = \som^{-1}(\del\hdf),
\xeee
that is, the $(0,1)$ part of the deformed Dolbeault differential is
\xlee{beq2.12b}
\dlb\p = \dlb + \som^{-1}(\del\hdf)\spsmb\del,
\xeee
while the symplectic form $\som$ is replaced by
\ylee{beq2.13b}
\som\p=\som +
d\hdf,
\yeee
so that it remains of type $(2,0)$ relative to the new
complex structure. In the formula\rx{beq2.12b}
we defined
$\som^{-1}\colon\Gamma(\Ts\xX)\rightarrow\Gamma(\Tng\xX)$
as the inverse of $\iota_{\xblnk}\,\som$.


If the deformation of $\Xsom$ is perturbative, that is, if $\hdf$ is a formal power series
\xlee{bbea2.3a1}
\hdfe =
\sum_{i=1}^\infty \hdfi\,\dfe^i,
\xeee
then the relation\rx{bbea2.2a} says that the leading coefficient $\hdfo$
must be $\dlb$-closed, and its gauge equivalence class is
determined by its Dolbeault cohomology class
$\hcdfo\in\Hdlbo(\xX)$, while the relation
\ylee{bbea2.3a1}
\dlb\hdftw + \shlf\,\Pbrv{\hdfo,\hdfo} = 0
\yeee
implies the  `integrability condition' for $\hcdfo$:
%
%
\ylee{bbea2.4}
\Pbrv{\hcdfo,\hcdfo} = 0.
\yeee

\subsubsection{Deformations of the \hlsm\ structure of a cotangent bundle}

If $\xX=\TsU$, then we restrict ourselves to \CMlt s $\hdf$ belonging to the subalgebra
$\zdpaTU\hookrightarrow \zdpaTsU$ defined in
subsection\rw{ss.dpahlsmm}.
Moreover, we consider only the deformations which do
not deform the complex structure of the zero section $\xcAz\subset\TsU$,
so we impose the condition $\bdmu|_{\xcAz}=0$ on the Beltrami
differential\rx{bbea2.2a1}.
This condition is satisfied if $\hdf$ is at least quadratic as a
function of holomorphic coordinates on fibers of $\TsU$ :
\xlee{bbea2.5} 
\label{bbea2.4a}
\hdf = \hdf_2 + \hdf_3 +\cdots,
\qquad
\hdf_i\in \xbOmov{\xcA,\Stai\TU}.
%
\xeee
The first two terms in this sum play a particularly important role in what follows and we give
them special names:
\xlee{bbea2.5b}
\hdf_2 = \hdfbt,\qquad\hdf_3 = \hdfgm.
\xeee
According to \ex{bbea2.2a}, they satisfy the equations
\xlee{beq2.25b}
\dlb\hdfbt=0,\qquad\dlb\hdfgm+ \shlf\Pbrv{\hdfbt,\hdfbt} = 0.
\xeee
Thus $\hdfbt$ is $\dlb$-closed, and its gauge equivalence class is
determined by the Dolbeault cohomology class that it represents:
\ylee{beq2.25b1}
\cdfbt\in\Hdlbo(\xcA,\Stat\TU),\qquad \Pbrv{\cdfbt,\cdfbt}=0.
\yeee

The class $\cdfbt$ has a simple geometric interpretation. For a
holomorphic submanifold $\yY\subset\xX$ of a complex manifold
$\xX$, the exact sequence of vector bundles on $\yY$
\xlee{beq2.25b2}
\TY\longrightarrow\Tng\xX|_{\yY} \longrightarrow\NY
\xeee
determines an extension class $\tdfbtY\in\Ext^1(\NY,\TY)$.
\ftnt{If $\xX$ is \Khl, then the class $\tdfbtY$ may be
represented by the anti-holomorphic part of the second fundamental
form of $\yY$ contracted with the \Khl\ metric and with its
inverse in order to turn two anti-holomorphic indices on the second fundamental form
into the holomorphic ones.}
If
$\yY$ is a lagrangian submanifold of a \hlsmm\ $\xX$ then the symplectic form $\som$
establishes an isomorphism $\NY\simeq\TsY$, so
$\tdfbtY\in\Ext^1(\TsY,\TY)$. The zero-section $\xcA\subset\TsU$
is a lagrangian submanifold and its exact sequence\rx{beq2.25b2}
splits, so in this case $\tdfbtU=0$. However, if we consider the
zero-section of the deformed bundle $\TsUhdf$ then the
sequence\rx{beq2.25b2} does not have to split. The injection
$\Hdlbo(\xcA,\Stat\TU)\hookrightarrow\Ext^1(\TsU,\TU)$ turns
$\cdfbt$ into an extension class within $\Ext^1(\TsU,\TU)$ and, in
fact,
\xlee{beq2.25b3}
\cdfbt = \tdfbtU.
\xeee
In other words, the leading coefficient $\cdfbt$ in the
expansion\rx{bbea2.4a} of $\hdf$ reflects the fact that the
sequence\rx{beq2.25b2} for the zero-section of the deformed
cotangent bundle $\TsUhdf$ does not split.



The injection\rx{bcea2.5} turns
an element $\hdf\in\xbOmbUUST$
into a $\bar{\Tng}^\vee\xcA$-valued function
(or, rather, a formal power series)
on the total space of $\TsU$. We denote this function by the same letter
$\hdf$.
The evaluation of $\hdf$ on a section of $\TsU$ gives a map
%
\xlee{bae2.0c1}
\hdf\colon\Gamma(\TsU)\rightarrow \xbOmoU.
\xeee

The restriction of the $(1,0)$ part of the differential $\del\hdf$
of an element $\hdf\in\xbOmbUST$
to the fibers of
$\TsU$ determines a vertical holomorphic differential map
\xlee{bae2.0c2}
\dfib\hdf\colon\Gamma(\TsU)\rightarrow \xbOmv{1}(\xcA,\TU).
\xeee

Recall that if $\xcW$ is a function on $\xcA$, then
$\yYW\subset\TsU$ denotes the graph of $\del\xcW$ defined by \ex{bah1.1}.
If $\xcW$ is holomorphic, then $\yYW$ is a \hlgrsm.
Let $\mtdfTUh$ denote the total space of the cotangent bundle
$\TsU$ whose \hlsm\ structure is deformed by
the \CMlt\rx{bbea2.4a}.
The deformed cotangent bundle $\mtdfTUh$ is canonically
diffeomorphic to $\TsU$, so for an arbitrary function $\xcW$, the graph
$\yYW$ is still a submanifold in $\mtdfTUh$, but this time $\yYW$ is
lagrangian if $\xcW$ satisfies the equation
%
\xlee{bae2.0a3}
\dlb\xcW = \hdf(\del\xcW),
\xeee
where $\hdf(\cdot)$ is the map\rx{bae2.0c1}.
If we consider a perturbative deformation\rx{bbea2.3a1}, then the
generating function becomes a formal power series
\xlee{bae2.0b1}
\xcWe = \xcWzz + \sum_{i=1}^{\infty}\xcWi\,\dfe^i.
\xeee
The leading term $\xcWzz$ is a holomorphic function describing
a \hlgrsm\  $\yY\subset\xX$ and it has a special property
\ylee{bae2.0b2}
\hcdfo(\del\xcWzz) = 0,
\yeee
which guarantees that $\yY$ can be deformed to the first order in $\dfe$, while the first order perturbation $\xcWo$
satisfies the equation
\ylee{bae2.0b3}
\dlb\xcWo - \hdfo(\del\xcWzz) = 0.
\yeee

The complex structure of the \hlgrsm\ $\yYW$ determined by the
function $\xcW$ satisfying the condition\rx{bae2.0a3} can be
described by saying that the bundle projection of $\TsU$
establishes an isomorphism between $\yYW$ and the base $\xcA$,
whose complex structure is deformed by the Beltrami differential
\xlee{bae2.0b4}
\mu_{\xcW} = -\dfib\hdf(\del\xcW),
\xeee
where $\dfib\hdf$ is the map\rx{bae2.0c2}.

\subsection{Deformation of the 2-category of \copfib s:
deformation of the category of morphisms}
\label{ss.dfmtcmr}

\subsubsection{Objects of the deformed category}
\label{ss.odc}

Following the outline of subsection\rw{ss.outline},
we conjecture that the \CMlt\  $\hdf$ parameterizing the
deformations of the \hlsms\ of $\TsU$, parameterizes also the
deformations of the 2-category $\rDDprfU$. We are going to discuss the
structure of the deformed category $\rDDprfUSk$, but we will limit
ourselves to simplest \xzobj s in it, which are the analogs of \copfib\
objects $\TrfUW$ denoted simply as $\xcW$.


Recall that
in the undeformed case the function $\xcW$ labeling an object of
$\rDDprfU$ is holomorphic,
so that the graph of its holomorphic
differential\rx{bah1.1} is a lagrangian submanifold of $\TsU$.
%
We conjecture that in the case of $\rDDprfUSk$, a
similar  object
is (parameterized by) a function $\xcW$ on $\xcA$, which satisfies
the equation\rx{bae2.0a3}, because then the graph of its
holomorphic differential $\yYW$ defined by the same
equation\rx{bah1.1} is again a lagrangian submanifold of the
deformed cotangent bundle $\mtdfTUSk$, and this is in line with
our conjecture that lagrangian submanifolds represent the objects
of $\ctLLXsom$ not only when $\Xsom$ is an undeformed cotangent bundle,
but also when it is a general \hlsmm.

%

\subsubsection{The universal \CMlt}

Consider the tensor algebra over $\IC$ of \emph{Dolbeault tensor fields}
\ylee{bah3.1}
\xTflU \edfn \bigoplus_{\ixk,\ixl=1}^\infty \TflklU,\qquad
\TflklU \edfn \xbOmb(\Tng^{\ixk}\xcA\otimes\Tng^{\vee,\ixl}\xcA).
\yeee
Fix a $\del$-connection
$\nabla\colon\TflbbU\rightarrow\Tflvv{\bullet}{\bullet+1}(\xcA)$.
For a tensor field $\tnsct\in\xTflU$, let $\bnbl\tnsct$ denote a sequence of
multiple covariant derivatives:
$\bnbl\tnsct\edfn\tnsct,\nabla\tnsct,\nabla^2\tnsct,\ldots$. For
tensor fields $\tnsct_1,\ldots,\tnsct_k$, let
$\Tflnv{\tnsct_1,\ldots,\tnsct_k}$ denote a subalgebra of $\xTflU$
generated by all tensor fields
$\bnbl\tnsct_1,\ldots,\bnbl\tnsct_k$ and by all possible contractions
within their tensor products.

Let $\hdf_i$ denote an element of $\xbOmov{\xcA,\Stai\TU}$. The Poisson-Schouten
bracket $\Pbrv{\hdf_i,\hdf_j}$ of two such elements
is an element of $\Tflnv{\hdf_i,\hdf_j}$
and it is universal in the sense that the coefficients of its expression in terms
of the appropriate contractions of $\hdf_i\otimes\nabla\hdf_j$ and
$\hdf_j\otimes\nabla\hdf_i$ are universal constants. The same
holds true for elements
$\xdfm_{ij}\in\xbOmv{\xmi}(\xcA,\wedge^{\xmj}\TU)$ and their
\NSb\ $[\xblnk,\xblnk]$.

For a \CMlt\ $\hdf\in\xbOmv{1}(\xcA,\Sb\TU)$ satisfying
\ex{bbea2.2a} and for two functions $\xcWo$, $\xcWt$ on
$\xcA$ satisfying \ex{bae2.0a3}, denote
\xlee{bah4.1}
\Tflnbot\edfn\Tflnv{\del\xcWo,\del\xcWt,\crR,\hdf_2,\hdf_3,\ldots},
\xeee
where $\hdf_i$ are the components\rx{bbea2.4a}, while
$\crR\edfn[\dlb,\nabla]\in\xbOmov{\rmS^2\TsU\otimes\TU}$
is the curvature of the tangent bundle $\TU$ corresponding
to the connection $\nabla$. The Dolbeault differential $\dlb$ acts
universally on the elements of $\Tflnbot$. Indeed, its action on
the components of $\hdf$ is prescribed by the \CMe\rx{bbea2.2a},
its action on $\del\xcWi$ is prescribed by \ex{bae2.0a3}, $\dlb\crR=0$, and
a permutation of $\dlb$ and $\nabla$ generates the curvature tensor $\crR$.

\begin{conjecture}
\label{cnj.univ}
For a \CMlt\ $\hdf\in\xbOmv{1}(\xcA,\Sb\TU)$ satisfying
\ex{bbea2.2a} and for two functions $\xcWo$, $\xcWt$ on
$\xcA$ satisfying \ex{bae2.0a3} there exists a universal \xrlb\ \CMlt\
\begin{gather}
\label{beq2.15b}
\xdfmot = \sum_{\xmi=0}^\infty \xdfmoti,\qquad
\xdfmoti\in\xbOmv{\xmi}(\xcA,\wedge^{\xmi}\TU),
\\
\label{beq2.17b}
\dlb\xdfmot + \shlf\,[\xdfmot,\xdfmot] = 0
\end{gather}
%
%
such that
\xlee{beq2.16b}
\xdfmotz=\xcWot\edfn
\xcWt-\xcWo
\xeee
and $\xdfmoti\in\Tflnbot$ for $i\geq 1$,
where $\hdf_i$ are the components\rx{bbea2.4a}.
The universality of $\xdfmot$ means
that the coefficients in the expression of its components
$\xdfmoti$ in terms of the tensor fields and their derivatives
are constants that do not depend on $\xcA$. The universal element
$\xdfmot$ is unique up to gauge equivalence\rx{bcea2.1a}, and different choices of
$\nabla$ also lead to gauge equivalent elements $\xdfmot$.
\end{conjecture}

Simply put, if two functions $\xcWo$, $\xcWt$
satisfy \ex{bbea2.2a} then their difference is not necessarily
holomorphic and hence it cannot serve as a deformation parameter
of the category $\rDprfU$. However, we conjecture that there is a
special unique correction to $\xcWot$ which turns it into a
\CMlt\ suitable for deforming  $\rDprfU$.
Hence we conjecture that the category of morphisms between the
objects of $\rDDprfUSk$ represented by $\xcWo$ and $\xcWt$ is the
deformed category $\rDprfU$:
%
%
\xlee{beq2.14b}
\Hom_{\rDDprfUSk}(\xcWo,\xcWt) =
 \rDsprfUdot,
\xeee
where $\xdfmot$ is the unique universal deformation parameter of
Conjecture\rw{cnj.univ}. Note that for a fixed manifold $\xcA$,
the sum in \ex{beq2.15b} is effectively finite, the highest value
of $i$ being the complex dimension of $\xcA$.

\subsubsection{Perturbative construction of the universal \CMlt}
\label{ss.pcucmlt}

We construct the universal expression for $\xdfmot$ perturbatively
in the Dolbeault degree $\dgDlb$ defined in an obvious way:
\xlee{e.dgdlb}
\qquad\dgDlb \xcWi=0,\qquad \dgDlb\crR = 1,\qquad
\dgDlb\hdfi=1,\qquad\dgDlb\nb=0.
\xeee
We substitute the
expansion\rx{beq2.15b} into the \CMe\rx{beq2.17b} and find the
equation determining $\xdfmotn$ in terms of $\xdfmoti$ with
$\xmi<\xmn$:
\xlee{beq2.18b}
[\xcWot,\xdfmotn] = - \dlb\xdfmotnmo - \hlf\sum_{i=1}^{\xmn-1}
[\xdfmoti,\xdfmotv{\xmn-i}].
\xeee

We will find expressions for the first three corrections in the
expansion\rx{beq2.15b}. We introduce a special notation for them:
\ylee{beq2.18ba1}
\xdfmoto = \bdmu,\qquad
\xdfmott = \bdnu,\qquad
\xdfmoth = \bdxi.
\yeee
(the indices $12$ at $\bdmu$, $\bdnu$ and $\bdxi$ are dropped temporarily in this subsection).
We are particularly interested in the third correction, because,
as we will see shortly, this is the dominant (that is, the lowest Dolbeault degree)
term in $\xdfmot$
when $\xcWo=\xcWt=0$, that is, from the
$\ctLLTUSk$ perspective
it describes the leading deformation
of the category of endomorphisms of the zero-section of $\TsU$.

According to \ex{beq2.18b}, the first three corrections are determined by
the equations
\begin{align}
\label{beq2.18ba2}
[\xcWot,\bdmu] & = - \dlb\xcWot,
\\
\label{beq2.18ba3}
[\xcWot,\bdnu] & = -\dlb \bdmu - \shlf\,[\bdmu,\bdmu],
\\
\label{beq2.18ba4}
[\xcWot,\bdxi] & = -\dlb \bdnu - [\bdmu,\bdnu].
\end{align}
%
%

First of all, we find an exact universal solution for $\bdmu$.
Equation\rx{beq2.18ba2} can be rewritten simply as
\xlee{beq2.19b}
\dlb\xcWot + \bdmu \spsmb \del \xcWot = 0.
\xeee
%
In order to solve it,
we
have to introduce \ddf\ notations. Let $\xaV$ and $\xbV$ be  vector
spaces. An element $\xaa\in\Sb\xaV\,\otimes\,\xbV$ determines a
polynomial function (or a formal power series) $\xaa\colon
\Vd\rightarrow\xbV$. The first \ddf\ of $\xaa$ is
a symmetric map $\fsd\xaa\colon \Vd \times\Vd\rightarrow\xaV\otimes\xbV$
defined by the property
\ylee{beq2.20b}
\fsd\xaa(\xavo,\xavt) \spsmb (\xavt - \xavo) = \xaa(\xavt) -
\xaa(\xavo),\qquad\forall\,\xavo,\xavt\in\xaV.
\yeee
The second \ddf\ of $\xaa$ is a totally symmetric
map $\ssd\xaa\colon\xaV\times\xaV\times\xaV\rightarrow
\Stav{2}\xaV\,\otimes\,\xbV$ defined by the property
\ylee{beq2.21b}
\ssd\xaa(\xavo,\xavt,\xavh) \spsmb (\xavh-\xavt)=
\fsd\xaa(\xavo,\xavh) - \fsd\xaa(\xavo,\xavt)
,\qquad\forall\,\xavo,\xavt,\xavh\in\xaV.
\yeee
In application to an element
$\hdf\in\xbOmov{\xcA,\;\Sb\TU}$ \ddfs\ produce symmetric maps
\ylee{beq2.22b}
\fsd\hdf\colon \Gamma(\TsU)^{\times 2}\rightarrow
\xbOmov{\xcA,\TU},\qquad
\ssd\hdf\colon \Gamma(\TsU)^{\times 3}\rightarrow
\xbOmov{\xcA,\Stav{2}\TU}.
\yeee
According to the condition\rx{bae2.0a3},
\ylee{beq2.24b}
\dlb\xcWot = \hdf(\del\xcWt) - \hdf(\del\xcWo) =
\fsd\hdf(\del\xcWo,\del\xcWt)\spsmb \del\xcWot.
\yeee
Hence the universal solution to \ex{beq2.19b} is the \ddf\
of $\hdf$ evaluated on the differentials $\del\xcWo$ and
$\del\xcWt$:
\xlee{beq2.23b}
\bdmu = - \fsd\hdf(\del\xcWo,\del\xcWt).
\xeee
%

We will find the universal expressions for $\bdnuot$  and
$\bdxiot$ only approximately, up to certain powers of $\xcWo$ and
$\xcWt$:
\xlee{beq2.23c}
\bdmu = \bdmupr + \OWh,\qquad
\bdnu = \bdnupr + \OWt,\qquad
\bdxi = \bdxipr + \OWo,
\xeee
where $\OWn$ denotes an expression which is at least of
combined degree $n$ in $\xcWo$ and $\xcWt$. According to
\ex{beq2.23b},
\ylee{beq2.24c}
\bdmupr = - \fsd\hdfbt(\del\xcWo,\del\xcWt)
- \fsd\hdfgm(\del\xcWo,\del\xcWt),
\yeee
where $\hdfbt$
and $\hdfgm$ are
defined by \ex{bbea2.5b}.
By
substituting \eex{beq2.23c} into \eex{beq2.18ba3}
and\rx{beq2.18ba4} we find
\begin{align}
\label{beq2.18bb3}
[\xcWot,\bdnupr] & = -\dlb \bdmupr - \shlf\,[\bdmupr,\bdmupr], 
\\
\label{beq2.18bb4}
[\xcWot,\bdxipr] & = -\dlb \bdnupr,
\end{align}
%
which determine the universal expressions for
$\bdnupr$ and
$\bdxipr$.

In order to simplify the calculations required to derive the
universal formula for $\bdnupr$ and $\bdxipr$, we present the
formula\rx{beq2.23b} in a different form. Consider the relation
between the Lie bracket on a manifold, and
the Poisson bracket on the total space of its cotangent bundle.
For a function $\xcW$ on a complex manifold $\xcA$, let $\tcW$
denote its pull-back to the total space of the cotangent bundle
$\TsU$. For an element $\mu\in\xbOmbv{\xcA,\TU}$ let $\tmu$ denote
the corresponding $(0,\bullet)$-form on the total space of $\TsU$ which is
linear along the fibers. In our previous notations,
$\mu=\dfib\tmu$. The Lie bracket on $\xcA$ and the Poisson bracket
on $\TsU$ are related as follows:
\ylee{beq2.28b}
\widetilde{\Lbrv{\mu,\xcW} } = -\Pbrv{\tmu,\tcW},\qquad
\widetilde{\Lbrv{\mu,\mu\p}} = -\Pbrv{\tmu,\tmu\p}.
\yeee

For a function $\xcW$ on $\xcA$, let $\hcW=\Pbrv{\tcW,\cdot}$ be a
linear operator acting on functions on the total space of $\TsU$.
The operators $\hcW$ commute with each other:
for two functions $\xcWo$ and $\xcWt$
\ylee{beq2.29b}
[\hcWo,\hcWt] = \widehat{ \Pbrv{\xcWo,\xcWt}   } = 0.
\yeee

For $\hdf\in\xbOmbUUST$ let $\hdf\xcmi$ denote its component in
$\xbOmbv{\xcA,\Stai\TU}$.
It is easy to see that in our new notations the \rhs of the
equations\rx{bae2.0a3} and\rx{beq2.23b} can be presented as
\begin{gather}
\label{beq2.30b}
\dlb\xcW = \hdf(\del\xcW) = \left.\lrbc{e^{\hcW}\hdf}\right|_{0},
\\
\bdmuot = - \fsd\hdf(\del\xcWo,\del\xcWt)
=-\dfib\left.\lrbc{\frac{e^{\hcWt}-e^{\hcWo}}{\hcWt-\hcWo}\,\hdf}\right|_1
\end{gather}
(in the \rhs of these formulas $\hdf$ is considered to be a
$\bar{\Tng}^\vee\xcA$-valued function
on the total space of $\TsU$).
According to the first formula,
\xlee{beq2.30c}
\dlb\xcW = \hdfbt(\dlb\xcW) + \OWh = \shlf\,\hcW^2\hdfbt + \OWh.
\xeee
According to the second formula,
\ylee{beq2.31c}
\bdmupr = -\dfib\Big(\shlf\, (\hcWo + \hcWt)\hdfbt
+ \tfrac{1}{6}\,(\hcWo^2 + \hcWo\hcWt + \hcWt^2)\,\hdfgm
\Big).
\yeee

Both sides of \ex{beq2.18bb3} are elements of
$\xbOmv{2}(\xcA,\TU)$,
so applying $\;\tilde{ }\;$ to them (that is, turning them into
$(0,2)$-forms on the total space of $\TsU$, which are linear along
fibers), we find
%
\begin{equation}
\label{beq2.35b}
\begin{split}
\widetilde{\Lbrv{\xcWot,\bdnupr} }
& =
- \dlb\tdmupr + \shlf\,\Pbrv{\tdmupr,\tdmupr} + \OWh
\\
& = \shlf\,\Pbrv{(\dlb\xcWo+\dlb\xcWt),\hdfbt}
+\tfrac{1}{6}\,(\hcWo^2 + \hcWo\hcWt + \hcWt^2)\,\dlb\hdfgm
\\
&\qquad + \tfrac{1}{8}\,
\Pbrv{ (\hcWo + \hcWt)\hdfbt,(\hcWo + \hcWt)\hdfbt } + \OWh
\\
& = - \tfrac{1}{8}\,\Pbrv{\hcWot\hdfbt,\hcWot\hdfbt}
+ \tfrac{1}{24}\,\hcWot^2\Pbrv{\hdfbt,\hdfbt}.
\end{split}
\end{equation}
We used the formulas\rx{beq2.25b}
and\rx{beq2.30c}
as well as the
Jacobi identity for the Poisson bracket in order to derive
the last line in this equation.

In order to solve the equation\rx{beq2.35b} for $\bdnupr$, we
express its \rhs in terms of a torsionless covariant $(1,0)$-differential on
the tangent bundle $\TU$:
$$\nb\colon\Gamma(\TU)\rightarrow\Omega^{1,0}(\xcA,\TU).$$
We use index notations: let $\vrxI$,
$\inI=1,\ldots,\dim_{\IC}\xcA$ be local holomorphic coordinates on
$\xcA$. The corresponding frames in $\TU$ and $\TsU$ are formed by
$\dlI$ and $d\vrxI$.
In our formulas we assume summation over repeated indices
appearing on opposite levels (this corresponds to applying
contraction to tensor products). Anti-holomorphic indices are
hidden. Thus in our notations
%
\begin{gather}
\nonumber
\hdfbt = \hdfbt^{\inI\inJ}\,\dlI\dlJ,
\qquad \hdfbt^{\inI\inJ} = \hdfbt^{\inJ\inI},
\qquad
\bdnu = \bdnu^{\inI\inJ}\,\dlI\wedge\dlJ,
\qquad
\bdnu^{\inI\inJ} = - \bdnu^{\inJ\inI}
\\
\label{beq2.37b}
\Pbrv{\xcW,\hdfbt} =
2\hdfbt^{\inI\inJ}\,(\dlI\xcW)\,\dlJ,\qquad
\Pbrv{\hdfbt,\hdfbt} =
-4\hdfbt^{\inI\inL}\,(\nbcL\hdfbt^{\inJ\inK})\,\dlI\dlJ\dlK,
\qquad
\Lbrv{\xcW,\bdnu} = 2\bdnu^{\inI\inJ}(\dlJ\xcW)\,\dlI.
\end{gather}
%
A straightforward computation shows that if we lift the tilde in
\ex{beq2.35b} by applying $\dfib$ to both sides, then its
\rhs can be rewritten as
\begin{multline}
\nonumber
\Lbrv{\xcWot,\bdnupr}
=
\\
\Bigg(
\hdfbt^{\inJ\inK}\hdfbt^{\inI\inL}\,(\dlJ\xcWot)\,(\nbcK\dlL\xcWot)
-\tfrac{1}{3}\Big( \hdfbt^{\inI\inL}\,(\nbcL\hdfbt^{\inJ\inK})
-
\hdfbt^{\inJ\inL}\,(\nbcL\hdfbt^{\inI\inK})\Big)(\dlJ\xcWot)(\dlK\xcWot)
\Bigg)\dlI.
\end{multline}
Comparing this expression with the last equation of\rx{beq2.37b}
we find a formula for the leading term in $\bdnu$:
%
\ylee{beq2.38b}
\bdnupr = \frac{1}{2}\;
\Bigg(
\hdfbt^{\inJ\inK}\hdfbt^{\inI\inL}\,(\nbcK\dlL\xcWot)
-\frac{2}{3}\; \hdfbt^{\inI\inL}\,(\nbcL\hdfbt^{\inJ\inK})
(\dlK\xcWot)
\Bigg)\dlI\wedge \dlJ .  
\yeee

Finally, we solve \ex{beq2.18bb4} for $\bdxipr$.
In order to compute its \rhs, we introduce the
Riemann curvature tensor $\crR\in\xbOmov{\xcA,\Stav{2}\TsU\otimes\TU}$
of the connection $\nb$ by the
formula $\crR = [\dlb,\nb]$. In index notations
\ylee{beq2.40b}
\crR = \crR^{\inI}_{\inJ\inK}\, d\vrxJ d\vrxK\,\dlI,\qquad
\crR^{\inI}_{\inJ\inK}=\crR^{\inI}_{\inK\inJ}.
\yeee
Then a straightforward computation shows that
\ylee{beq2.41b}
\dlb \bdnupr
 = \frac{1}{3}
\Big(
\hdfbt^{\inI\inL}\hdfbt^{\inJ\inM}\crR^{\inK}_{\inL\inM}
+
\hdfbt^{\inJ\inL}\hdfbt^{\inK\inM}\crR^{\inI}_{\inL\inM}
+
\hdfbt^{\inK\inL}\hdfbt^{\inI\inM}\crR^{\inJ}_{\inL\inM}
\Big)(\dlK\xcWot)\,\dlI\wedge\dlJ.
\yeee
Since in index notations
\ylee{beq2.42b}
\bdxi = \bdxi^{\inI\inJ\inK}\,\dlI\wedge\dlJ\wedge\dlK,\qquad
\Lbrv{\xcW,\bdxi} = 3\,
\bdxi^{\inI\inJ\inK}(\dlK\xcW)\,\dlI\wedge\dlJ,
\yeee
we find that
\xlee{beq2.42b}
\bdxipr = \frac{1}{3}\,\hdfbt^{\inI\inL}\hdfbt^{\inJ\inM}\crR^{\inK}_{\inL\inM}
\;\dlI\wedge\dlJ\wedge\dlK.
\xeee

\subsubsection{Semi-classical grading}
\label{ss.scgr}

The algebra $\Tflnbot$ of \ex{bah4.1} has an important
\emph{\smcl} grading defined as follows:
\xlee{bah4.2}
\spdg \xcW_i = -2,\qquad\spdg\crR=0 ,\qquad\spdg
\hdf_i=i-3,\qquad\spdg\nabla=1.
\xeee
The defining relation $\crR=[\dlb,\nabla]$ implies $\spdg\dlb=-1$.
If we set
\xlee{bah4.4}
\spdg\xdfmot=-2,
\xeee
then all three defining
relations of our deformation construction:
\ylee{bah4.3}
\dlb\hdf + \shlf\,\Pbrv{\hdf,\hdf} = 0,
\qquad
\dlb\xcW = \hdf(\del\xcW),
\qquad
\dlb\xdfmot + \shlf\,[\xdfmot,\xdfmot] = 0
\yeee
respect the \smclgrd. Hence the universal \CMlt, determined recursively
by \ex{beq2.18b}, must have the \smcldgr\ -2: its zeroth term
$\xdfmz=\xcWot$ has degree -2 in virtue of $\spdg\xcWot=-2$ and the
degree of higher terms is expressed through \ex{beq2.18b} in
terms of the degrees of the lower terms.

Conjecture\rw{cnj.univ} together with the relation\rx{bah4.4} has
three easy corollaries.

Consider a new grading on the algebra
$\Tflnbot$:
\xlee{bah2.4}
\dgdW\bnbl\hdf_i=\dgdW\bnbl\crR=0,\qquad
\dgdW \nabla^k \del \xcW =
\begin{cases}
1, & \text{if $k=0$,}
\\
0, & \text{if $k\geq 1$.}
\end{cases}
\xeee
\begin{corollary}
\label{crl1}
If $\hdfbt=0$, then $\dgdW\xdfmi\geq 2$ for $i\geq 1$.
\end{corollary}
\begin{corollary}
\label{crl2}
If $\hdfbt=0$ and $\xcWo=\xcWt=0$, then $\xdfm=0$.
\end{corollary}
Let $\ccrR\in\Ext^1(\rmS^2\TU,\TU)$ denote the Atiyah class of the
tangent bundle $\TU$ (it is the class represented by the curvature tensor
$\crR$).
\begin{corollary}
\label{crl3}
If $\ccrR=0$ and $\xcWo=\xcWt=0$, then $\xdfm = 0$.
\end{corollary}
%

Observe that among all generators
$\bnbl\hdf_i$, $\bnbl\crR$ and $\bnbl\del\xcWi$ of $\Tflnbot$ only
$\hdfbt$ and
$\del\xcW_i$ have negative \smcldgr s:
$\spdg\hdfbt=\spdg\del\xcW_i=-1$. At the same time, $\spdg\xdfm=-2$, so if $\hdfbt=0$,
then each
term in the expression of $\xdfmi$, $i\geq 1$ must have at least
two powers of $\del\xcW$. This proves Corollary\rw{crl1}. If $\xcWo=\xcWt=0$, then, obviously,
$\xdfmi=0$ for $i\geq 1$ and, at the same time,
$\xdfmz=\xcWot=0$. This proves Corollary\rw{crl2}.

In order to prove Corollary\rw{crl3}, we introduce two more
gradings on the algebra
$\Tflnbot$. The first grading reflects the difference between the
total numbers of upper and lower indices in a tensor field:
\ylee{bah6.2}
\dgbl\nabla=-1,\qquad\dgbl\hdf_i=i-1,\qquad\dgbl\del\xcWi=-1,
\qquad\dgbl\crR=-2.
\yeee
Since the universal deformation parameter $\xdfm$ is \xrlb, its
degree must be zero: $\dgbl\xdfm=0$.

The second degree is the sum: 
$\dgdf=\spdg + \dgbl$, so
\ylee{bah6.3}
\dgdf\nabla=0,\qquad\dgdf\hdf_i=2i-4,\qquad\dgdf\del\xcWi=-2,
\qquad\dgdf\crR=-2.
\yeee
Since $\spdg\xdfm=-2$ and $\dgbl\xdfm=0$, then $\dgdf\xdfm=-2$.
However, $\dgdf\hdf_i\geq 0$ for $i\geq 2$, so
if we set $\xcWo=\xcWt=0$ and $\crR=0$, then all remaining tensor
fields in $\Tflnbot$ must have non-negative degree, so $\xdfm$
must be zero. This proves Corollary\rw{crl3}.


\subsection{Deformation of the 2-category of \copfib s:
deformation of the composition of morphisms}

Describing the deformation of the composition functor requires the
following steps: first, for a pair of \xzobj s $\xcWo$, $\xcWt$ we
have to describe the \fDAinfalg\ $\brb{\xbOmbU,\dlb,\xdfmot}$
corresponding to the \CMlt\ $\xdfmot$ found in
subsection\rw{ss.dfmtcmr},
by working out the expressions for its \xmnmlt s.
Second, we have to describe the objects and morphisms of the
\tpdprc\ of $\brb{\xbOmbU,\dlb,\xdfmot}$. After that we can
define the composition between objects of
$\Hom(\xcWo,\xcWt)$ and $\Hom(\xcWt,\xcWh)$.

\subsubsection{A first order deformation of the \hlsm\ structure}

The procedure is simplified if we stay within the realm of
categories $ \rDsprfUWtmo$ which are very similar to derived
categories of coherent sheaves. We achieve this by considering only infinitesimal
deformations of $\rDDprfU$ by $\hdf$. In other words, we introduce the algebra
$\ICv{\dfe}/(\dfe^2)$ and consider a deformation of the \hlsms\ of
$\TsU$ by an element $\dfe\hdf$, where $\hdf$ is of the
form\rx{bbea2.4a}. The quadratic term in the \CMe\rx{bbea2.2a}
drops out, hence $\hdf$ must be holomorphic:
\ylee{\bbf1}
\dlb\hdf = 0.
\yeee
A function $\xcW$ describing an \xzobj\ of the deformed category has
the form
\ylee{bbf2}
\xcW=\acWz+\dfe\acWo.
\yeee
It must satisfy \ex{bae2.0a3}, which reduces to two
relations
\xlee{bbf3}
\dlb\acWz = 0,\qquad\dlb\acWo = \hdf(\del\acWz).
\yeee

Let us introduce two shortcut notations related to
functions $\xcW_i$ satisfying the conditions\rx{bbf3}:
\xlee{bbf3a}
\shdfij \edfn \fsd\hdf(\del\xcWvz{i},\del\xcWvz{j}),\qquad
\shdfijk \edfn \ssd\hdf(\del\xcWvz{i},\del\xcWvz{j},\del\xcWvz{k}).
\xeee

For two functions $\xcWo$, $\xcWt$ satisfying the
conditions\rx{bbf3} and thus defining \xzobj s of
the deformed category $\rDDprfUhae$, the corresponding \CMlt\ has
the form
\xlee{bbf4}
\xdfmot = \xcWotz + \dfe(\xcWoto + \bdmuot),
\yeee
where, according to \ex{beq2.23b},
\ylee{bbf5}
\bdmuot =
-\shdfot
= - \fsd\hdfbt(\del\xcWoz+ \del\xcWtz)
+ \OWt.
\yeee
This \xpBd\ satisfies the equations
\ylee{bbf6}
\dlb\bdmuot = 0,\qquad\dlb\xcWoto + \bdmuot \spsmb \del\xcWotz=0,
\yeee
but generally does not satisfy the integrability condition
$\Lbrv{\bdmuot,\bdmuot}=0$.

\subsubsection{Deformation of the \tZtgdcs}
\def\xcWz{ \xcW }
According to \ex{bbf4}, the category of morphisms between two \xzobj s $\xcWo$, $\xcWt$
of the deformed 2-category $\rDDprfUSek$ is
\xlee{bbf7}
\Hom_{\rDDprfUSek}(\xcWo,\xcWt) = \rDsprfBvv{\xcA}{\xcWotz + \dfe(\xcWoto +
\bdmuot)}
= \rDsprfUeWot,
\yeee
where $\xcAekot$ denotes the complex manifold $\xcA$ whose complex
structure is deformed by the Beltrami differential $\dfe\bdmuot$.
Hence the category\rx{bbf7} is defined along the lines of
subsection\rw{tZtgdcs}. However, before we go into specifics, let
us recall the definition of the \cAc.

Recall that a \xper\ object\rx{bae1.44} of a 2-periodic category $\rDsprfUW$
defined in
subsection\rw{tZtgdcs} is determined by a pair $\cE=(\xE,\nbbE)$, where
$\xE$ is a
\cqhlmvb\
%
over a complex manifold $\xcA$,
while $\nbbE$ is its curved differential satisfying the
properties\rx{bae1.42}--(\ref{bae1.42a1}).
Let us endow $\xE$ also with a (possibly curved) \ydh\ covariant
differential
\ee
\xlabel{bae2.2b7}
\xbOmxvb{i}(\xE) \xrarv{\nbE} \xbOmxvb{i+1}(\xE),\qquad
\zdgt{\nbE} = \yev,
\eee
which satisfies the analog of \ex{bae1.42a}:
\ee
\label{bae2.2b8}
\nbE\,(\smu\wedge\sigma) = (\del\smu)\wedge \sigma +
(-1)^{\dgZt\smu}\smu\wedge (\nbE\,\sigma).
\eee
The choice of $\nbE$ is not unique, and the difference of two differentials
$\nbE$,  $\nbEp$ satisfying\rx{bae2.2b8} is a
differential form
\ee
\label{bae2.2b9} 
\nbEp = \nbE + \sea,\qquad\sea
\in
\xbOmxob(\End\xE).
\eee
In other words, all possible differentials $\nbE$ form an affine
space based on a vector space $\xbOmoeE$. The commutator of
\ydh\ and \ydah\ differentials is the $(1,\hat 1)$-curvature
\ylee{bbf8}
[\nbbE,\nbE ]=\acFE\in
\xbOmxob(\End\xE),
\yeee
which satisfies the curved Bianchi identity
\xlee{bbf9}
\nbbE\acFE = -(\del\xcW)\,\xIdv{\xE}.
\yeee
The curvature $\acFE$ is determined by the object $\cE$ (that is, by
$\nbbE$) up to a $\nbbE$-exact element: if we replace $\nbE$ by
$\nbEp$ of \ex{bae2.2b9}, then $\acFE$ is replaced by
$\acFpE= \acFE + \nbbE\sea$. Hence $\cE$ determines the
\emph{\cAc}
\ylee{bbf10}
\hacFE\in
\xbOmxob(\End\xE)
\Big/ \nbbE \lrbc{
\xbOmxob(\End\xE)
}.
\yeee
If $\xcW=0$, then the Bianchi identity\rx{bbf9} implies that
\ylee{bbf11}
\hacFE\in\HnbbE(\End\xE) = \Ext(\xcE,\xcE).
\yeee

A perfect object of the deformed category $\rDsprfUeWot$ is a pair
\ylee{bbf12}
\cE=\adgmE,\qquad
\nbbE = \nbbzE + \dfe\nbboE,
\yeee
(\cf
\ex{bae1.45b}), such that the pair $\xcEbz=(\xE,\nbbzE)$ is an
object of the undeformed category $\rDsprfUWot$, while
\xlee{bae2.2a10}
\nbboE = \ydfmotv{\nbE} + \seb,\qquad
\seb\in\xbOmbv{\End\xE}
\xeee
and $\seb$ satisfies the condition
\ee
\label{bae2.2a9} 
\nbbzE\,\seb = \xcWoto\,\xIdE - \ydfmotv{\acFE}.
\eee
Here $\acFE$ is the $(1,1)$-curvature of $\xcEbz$, so $\nbbE$
satisfies the condition\rx{bae1.42a1}: $\nbbE^2 =
\xcWot\xIdv{\xE}$.
A change\rx{bae2.2b9} in the choice of $\nbE$ is
compensated by the corresponding replacement of $\seb$ by
$\seb\p=\seb - \ydfmotv{\sea}$.

The space of morphisms 
$\Hom_{\rDsprfUeWot}(\xcEo,\xcEt)$
between two \xper\ objects is defined by means of an obvious deformation of the general
formula\rx{bae1.7}. A morphism between two objects
\wlee{bae2.2a11}
\msgm
\in\Hmo{\cEo,\cEt}
\weee
is represented by a $\nbbE$-closed sum
%
\begin{gather}
\nonumber
\msgm  = \msgmbz + \dfe\msgmbo,\qquad
\msgmbz,\msgmbo\in\xbOmbv{\xEt\otimes\xEo^\ast},
\\
\label{bae2.2a14}
\nbbzcxE\,\msgmbz = 0,\qquad
\nbbzcxE\,\msgmbo = - \nbbocxE\,\msgmbz
\end{gather}
up to $\nbbE$-exact elements. Note that the dominant component $\msgmbz$
defines a morphism between the undeformed objects $\xcEobz$ and
$\xcEtbz$.

%

\subsubsection{Deformation of the composition of morphisms}
\label{ss.dfmcmpf}
The composition of morphisms between three \xzobj s
$\xcWo$, $\xcWt$ and $\xcWh$ of the deformed
2-category $\rDDprfUSek$
 is a
bi-functor
\ee
\label{bae2.8} 
\rDsprfUeWot \times \rDsprfUeWth \longrightarrow \rDsprfUeWoh.
\eee
The composition of two morphisms $\xcEot\in\rDsprfUeWot$ and
$\xcEth\in\rDsprfUeWth$
is the appropriately deformed tensor product:
%
\xlee{bae2.10} 
\xcEth\circ\xcEot=
\adgmeEth
\circ
\adgmeEot =
(\xEth\otimes\xEot,\nbbv{\xEot\otimes\xEth} + \dfe\nboth),
\xeee
where the deformation term is
\xlee{bae2.11}
\begin{split}
\nboth & =
\shdfoth
\spsmb
\brB{(\del\xcWotz)\nbEth - (\del\xcWthz)\nbEot+ \acFEot\acFEth}
=
\\
& =
\hdfbt\spsmb (\acFEot\acFEth) + \OWo,
\end{split}
\xeee
and $\shdfoth$ is a shortcut notation defined by \ex{bbf3a}.

The first two terms in the \rhs of this equation are related to the fact
that each of three categories in\rx{bae2.8}
has its own deforming \xpBd. Hence we had
to add  correction terms to $\nbbv{\xEoh}$
so that it
would satisfy the
condition\rx{bae2.2b8}
or, more precisely, the
condition\rx{bae2.2a10}, that is, that the difference
%
\ee
\xlabel{bae2.12}
\nbbv{\xEoh}
-\dfe\bdmuoha{\nabla_{E_{12}\otimes E_{23}}}
\eee
must be just an odd element of $\xbOmbv{\End\xEoh}$. The third correction term
in \ex{bae2.11} is required to comply with the
condition\rx{bae2.2a9}.

The composition of morphisms $\xcEth\circ\xcEot$ defined by \ex{bae2.10}
is independent (up to an isomorphism) of the choice of \ydh\
differentials $\nbEot$ and $\nbEth$. Indeed, if we replace
$\nbEot$ with $\nbEot\p=\nbEot+\sea$ as in \ex{bae2.2b9},
then $\nboth$ is replaced by
%
%
\begin{equation}
\nonumber
\begin{split}
\nboth\p
&=
\nboth -
\shdfoth
\spsmb
\brB{  (\del\xcWthz)\sea - (\nbbzEot\;\sea)\acFEth}
\\
&=
\nboth + \nbbzEot\, \brB{ \shdfoth \spsmb (\sea\,\acFEth) },
\end{split}
\end{equation}
%
%
%
so it changes by a $\nbbzEoh$-exact term.

The bi-functorial nature of the map\rx{bae2.8} means that an
object
%
$\xcEth\in\rDsprfUeWth$
%
determines a functor
\ee
\xlabel{bae2.14}
\xymatrix@C=2cm{
\rDsprfUeWot
\ar[r]^-{\xPhEth} &
\rDsprfUeWoh
}.
\eee
Its action on objects is defined by the composition\rx{bae2.10}.
Consider now its action on morphisms. For
$\msgmot\in\Hom(\xcEot,\xcEotp)$, where
$\xcEot,\xcEotp\in\rDsprfUeWot$,
we define
\ee
\xlabel{bae2.16}
\xPhEth(\msgmot) = \msgmot\otimes\xIdv{23}
+ \dfe\, (\shdfoth) \spsmb \brb{\nbEot\msgmotz\otimes\acFEth}
.
\eee
The correction term
is required to satisfy the
relation\rx{bae2.2a14}.
Similarly, an object $\xcEot\in\rDsprfUmWot$
determines a functor
\ee
\xlabel{baez2.17}
\xymatrix@C=2cm{
\rDsprfUeWth
\ar[r]^-{\xPhEot} &
\rDsprfUeWoh
},
\eee
which maps a morphism $\msgmth\in\Hom(\xcEth,\xcEth\p)$,
where $\xcEth,\xcEth\p\in\rDsprfUeWth$
into a morphism
\ee
\xlabel{bae2.19}
\xPhEot(\msgmth) = \xIdv{12}\otimes\msgmth
- \dfe\, (\shdfoth) \spsmb \brb{\acFEot\otimes\nbEth\msgmthz}.
\eee

The images of morphisms $\msgmote$ and $\msgmthe$ commute in the
following sense:
\begin{multline}
\xlabel{bae2.20}
\xPhEpot(\msgmth)\circ\xPhEth(\msgmot) -
\xPhEpth(\msgmot)\circ\xPhEot(\msgmth)
=0
\\
\text{in}\quad
\Hom\brB{
\xcEth
\circ
\xcEot,\xcEth\p
\circ
\xcEot\p
},
\end{multline}
because
\begin{multline}
\xlabel{bae2.21}
\xPhEpot(\msgmth)\circ\xPhEth(\msgmot) -
\xPhEpth(\msgmot)\circ\xPhEot(\msgmth)
\\
=
 \dfe\, \nbbv{\xEoh|0}
\brB{(\shdfoth) \spsmb \brb{\nbEot\msgmotz\otimes\nbEth\msgmthz}}.
\end{multline}

In the special case $\xcWo=\xcWt=\xcWh=0$, the
formula\rx{bae2.11} says the following. Let $0$ denote the object of $\rDDprfUSek$ corresponding to the trivial fibration over $U$ with $W=0$. The endomorphism category
$\End_{\rDDprfUSek}(0)$ is a monoidal category which is equivalent to $\rDprfU$ as a category, but with a monoidal structure given by the deformed tensor product
\xlee{baj4.1}
\adgmeE
\circ
\adgmeEp =
\brb{\xE\otimes\xEp,\nbbv{\xE\otimes\xEp} + \dfe\hdfbt\spsmb (\acFE\acFv{\xEp})},
\xeee

\subsubsection{Deformation of the monoidal structure beyond the
first order}
\label{ss.dmsbfo}

Now we return to the 2-category $\rDDprfUSk$ for a general \CMlt\
$\hdf$ and consider the category
$\xEndsp$ of endomorphisms of the \opfib\ with $\xcW=0$ denoted here
as $\zsU$. This
category has a monoidal structure corresponding to the composition
of endomorphisms. According to the general formula\rx{beq2.14b},
the endomorphism category itself is a deformation of the category
$\rDDprfU$: $\xEndsp=\rDsprfUd$, where
$\xdfm = \xdfmv{3}+\xdfmv{4}+\cdots$ and
$\xdfmi\in\xbOmv{\xmi}\brb{\xcA,\wedge^{\xmi}\TU}$, while its
monoidal structure is a deformation of the monoidal structure of
$\rDDprfU$, the latter being the tensor product\rx{bae1.45a2}.

Let us assume that the Atiyah class $\ccrR$ of the tangent bundle
$\TU$ is zero. Then, according to Corollary\rw{crl3}, the deformation
parameter is zero: $\xdfm=0$,
so we have an equivalence of categories
\xlee{bai5.1}
\xEndsp\simeq\rDprfU.
\xeee
However, as the study of the first order perturbation in subsection\rw{ss.dfmcmpf}
demonstrated, the monoidal structure of $\xEndsp$ is still a
non-trivial deformation of the tensor product monoidal structure
of $\rDprfU$. The relatively simple nature of the
category\rx{bai5.1} allows us to discuss the properties of this
deformation without invoking \Ainfalg s and their modules.

A deformation of the monoidal structure of the category $\rDprfU$
is described by two sets of data. First, for every pair of \qhlmvb s
$\adgmeEo$, $\adgmeEt$ there is a \CMlt\
$\zdfmot\in\xbOmbz\brb{\End(\xEo\otimes\xEt)}$,
\xlee{bai5.1a}
\dlb\zdfmot + \shlf\,[\zdfmot,\zdfmot]=0,
\xeee
which determines the deformed monoidal bifunctor of the composition within the endomorphism
category $\xEndsp$:
\xlee{bai5.2}
\adgmeEo\ycrc\adgmeEt = (\xEo\otimes\xEt,\bar\nabla_{E_1\otimes E_2}  +
\zdfmot).
\xeee
Second, for every triple of \qhlmvb s $\adgmeEo$, $\adgmeEt$,
$\adgmeEh$ there is an associator
$\zasoth\in\xbOmbz\brb{\End(\xEo\otimes\xEt\otimes\xEh)}$ which
establishes the associativity isomorphism
\ylee{bai5.3}
\zasoth\colon \brB{\adgmeEo\ycrc\adgmeEt}\ycrc\adgmeEh
\xmapta{\cong}
\adgmeEo\ycrc\brB{\adgmeEt\ycrc\adgmeEh}.
\yeee
If both sides of the associativity isomorphism have the presentation
\begin{equation}
\label{bai5.4}
\begin{split}
\brB{\adgmeEo\ycrc\adgmeEt}\ycrc\adgmeEh
&=
(\xEo\otimes\xEt\otimes\xEh,\nbboth+\zdfmoth),
\\
\adgmeEo\ycrc\brB{\adgmeEt\ycrc\adgmeEh}
&=
(\xEo\otimes\xEt\otimes\xEh,\nbboth+\zdfmothp),
\end{split}
\end{equation}
where $\nbboth \edfn \bar\nabla_{E_1\otimes E_2\otimes E_3}$,
then $\zasoth$ is an invertible element satisfying the equation
\xlee{bai5.5}
\nbboth\zasoth + \zdfmothp\,\zasoth - \zasoth \, \zdfmoth = 0.
\xeee

We conjecture that there exist unique universal formulas for the element
$\zdfmot$ and for the associator $\zasoth$ related to the deformation of the tensor product
monoidal structure of $\rDprfU$ into the monoidal structure
of the endomorphism category $\xEndsp$. These universal formulas
express
$\zdfmot$ and $\zasoth$ in terms of the deformation parameter
$\hdf$, \ooc s $\acFEi$ and their holomorphic covariant
derivatives:
\xlee{bai5.6}
\zdfmot\in
\Tflnv{\acFEo,\acFEt,\hdf_2,\hdf_3,\ldots}
 , \qquad
\zasoth\in
\Tflnv{\acFEo,\acFEt,\acFEh,\hdf_2,\hdf_3,\ldots}.
\xeee

We propose to derive the universal formulas perturbatively. Define
the Dolbeault degree by \ex{e.dgdlb} and by the additional formula
$\dgDlb\acFEi=1$ (note that generally $\acFEi\in\xbOmb(\End\xEi)$
and its Dolbeault degree coincides with $j$ of $\xbOmv{j}$ only when
$\xEi$ is a holomorphic vector bundle with the operator $\nbbv{\xEi}$ not containing forms of degree other than $1$).
We present the deformation parameter $\zdfmot$
and the associator $\zasoth$ as the sums
\ylee{bai5.7}
\zdfmot = \sum_{i=1}^{\infty}\zdfmoti,\qquad
\zasoth = \xIdEm + \sum_{i=2}^\infty\zasothi,
\yeee
where
\ylee{bai5.8}
\dgDlb\zdfmoti = 2i + 1,\qquad \dgDlb\zasothi = 2i
\yeee
(the reason for assuming that the Dolbeault degree of $\zdfmot$ is
even and the Dolbeault degree of $\zasoth$ is odd will become
clear shortly).
The \CMe\rx{bai5.1a} splits:
\xlee{bai5.8}
\dlb\zdfmotn + \sum_{i=1}^{\xmn-1} \zdfmoti\,\zdfmotv{\xmn-i}=0,
\xeee
and the associativity equation\rx{bai5.5} splits:
\xlee{bai5.9}
\dlb\zasothn + \zdfmothnp - \zdfmothn +
\sum_{i=2}^{\xmn-1} (\zdfmothv{\xmn-i}\p\,\zasothi -
\zasothi\,\zdfmothv{\xmn-i})=0.
\xeee
The action of the Dolbeault differential $\dlb$ on the
elements of the algebras\rx{bai5.6} follows from its action on the
elementary tensor fields prescribed by the \CMe\rx{bbea2.2a} and
by the Bianchi identity $\nbbEi\acFEi=0$
\ftnt{We assume for simplicity that the curvatures of the $\partial$-connections $\nbv{\xEi}$
are zero.}, and from the defining
equation of the curvature tensor $\acFEi = [\nbbEi,\nbv{\xEi}]$.

We introduce the notation $\zdfmdv{\xblnk,\xblnk}$ to emphasize
the dependence of the universal deformation parameter $\zdfmot$ on
curvatures: $\zdfmot = \zdfmdv{\acFEo,\acFEt}$. The
parameters $\zdfmoth$ and $\zdfmothp$ of \ex{bai5.4} can be
expressed in terms of $\zdfmdv{\xblnk,\xblnk}$:
\begin{equation}
\xlabel{bai5.10}
\begin{split}
\zdfmoth &= \zdfmdv{\acFEo,\acFEt} + \zdfmdv{\acFEo+\acFEt+\nbv{\xEo\otimes\xEt}
\zdfmdv{\acFEo,\acFEt},\acFEh},
\\
\zdfmothp &=\zdfmdv{\acFEt,\acFEh} + \zdfmdv{\acFEo,\acFEt+\acFEh
+\nbv{\xEo\otimes\xEt}\zdfmdv{\acFEt,\acFEh} }.
\end{split}
\end{equation}
These formulas allow us to present the difference $\zdfmothnp -
\zdfmothn$ appearing in \ex{bai5.9} in the form
\xlee{bai5.11}
\zdfmothnp - \zdfmothn = \adltothzn + \tzdfmoth,
\xeee
where
\ylee{bai5.12}
\zdfmothn = \zdfmdnv{\acFEt,\acFEh} +
\zdfmdnv{\acFEo,\acFEt+\acFEh} -
\zdfmdnv{\acFEo,\acFEt} -\zdfmdnv{\acFEo+\acFEt,\acFEh}
\yeee
and the expression $\tzdfmoth$ contains the deformation parameter
components $\zdfmi$ only with $i<n$.

After the substitution\rx{bai5.11}, the associativity
equation\rx{bai5.9} becomes
\xlee{bai5.13}
\dlb\zasothn + \adltothzn + \tzdfmoth +
\sum_{i=2}^{\xmn-1} (\zdfmothv{\xmn-i}\p\,\zasothi -
\zasothi\,\zdfmothv{\xmn-i})=0.
\xeee

We conjecture that the \CMe\rx{bai5.8} together with the associativity
equation\rx{bai5.13} can be solved perturbatively
over the Dolbeault degree, thus producing the unique universal solutions
$\zdfmot$ and $\zasoth$ if, following \ex{baj4.1}, we set
\xlee{bai5.14}
\zdfmoto=\hdfbt\spsmb (\acFEo\acFEt)=
\hdfbt^{\inI\inJ} \acFv{\xEo,\inI}\acFv{\xEt,\inJ},
\xeee
where
we used the index notations explained at the end of
subsection\rw{ss.pcucmlt}, as well as the notation $\acFE =
\acFv{\xE,\inI}\,dx^{\inI}$ for the \ooc\ tensor components.
The parity of Dolbeault degrees of $\zdfmoti$ and $\zasothi$ is
dictated by these equations.


It is easy to verify that the expression\rx{bai5.14} satisfies
\eex{bai5.8} and\rx{bai5.13} for $\xmn=1$.
We leave it for the reader to verify that the following
expressions
%
\begin{align}
\xlabel{bai5.15}
\zdfmott &=
\tfrac{1}{3}\,
\hdfbt^{\inJ\inL}(\nb_{\inL}\hdfbt^{\inI\inK})
\acFv{\xEo,\inI}(\acFv{\xEo,\inJ}-\acFv{\xEt,\inJ})\acFv{\xEh,\inK}
\\
&\qquad\qquad
+\shlf\,\hdfbt^{\inI\inJ}\hdfbt^{\inK\inL}\brB{
(\nb_{\inI}\acFv{\xEo,\inK})\, \acFv{\xEt,\inJ} \acFv{\xEt,\inL}
+
(\nb_{\inI}\acFv{\xEt,\inK})\, \acFv{\xEo,\inJ} \acFv{\xEo,\inL}
},
\\
\label{bai5.15a}
\zasotht &= \tfrac{2}{3}\,\hdfgm^{\inI\inJ\inK}\acFv{\xEo,\inI}\acFv{\xEt,\inJ}\acFv{\xEh,\inK}.
\end{align}
satisfy these equations
for $\xmn=2$.

We extend the \smcl\ grading of subsection\rw{ss.scgr} to the
algebras\rx{bai5.6} by setting $\spdg \acFEi=0$.
Solving the system of equations\rx{bai5.8} and\rx{bai5.9}
recursively over $\zn$ with the initial condition\rx{bai5.14}
determines the \smcl\ degrees of the deformation parameter
$\zdfmot$ and of the associator $\zasoth$:
\ylee{bai5.16}
\qquad \spdg\zdfmot = -1,\qquad \spdg\zasoth = 0.
\yeee

Let us consider what happens if we set $\hdfbt=0$ in the universal formulas for $\zdfmot$ and
$\zasoth$. Since $\hdfbt$ is the only generator of the
algebras\rx{bai5.6} with negative \smcl\ grading, then
$\spdg\zdfmot=-1$ implies
\ylee{bai5.17}
\zdfmot|_{\hdfbt=0}=0,
\yeee
that is, the composition part\rx{bai5.2} of the monoidal structure
remains undeformed.
However, the equation\rx{bai5.15a} indicates that if $\hdfgm\neq 0$, then
$\zasoth|_{\hdfbt=0}\neq \xId$. This means that if for a complex
manifold $\xcA$ there exists a non-trivial \CMlt\
$\hdf\in\xbOmbUUST$ such that $\hdfbt\equiv\hdf_2=0$ while
$\hdfgm\equiv\hdf_3\neq 0$, then the tensor product monoidal structure
of the category $\rDprfU$ has a non-trivial associator
$\zasoth\neq \xId$ in addition to the standard one. This situation is realized, for example, when $U$ is a holomorphic symplectic manifold $X$ \cx{Kapranov}.
The element $\hdf$ in this case describes the formal neighborhood of the diagonal in $X\times X$. It follows that the $\ZZ_2$-graded derived category of any holomorphic symplectic manifold admits a non-trivial monoidal structure with a deformed associator. This provides an underlying reason  for the results of J. Roberts and S. Willerton \cx{RobWil}.

If $\hdfbt=0$, then all remaining generators of the
algebras\rx{bai5.6} have non-negative \smcl\ degrees. Among them, only
$\acFEi$ and $\hdfgm$ have zero degrees, and all others, including
holomorphic derivatives, have positive degrees. Since
$\spdg\zasoth=0$, this means that $\zasoth$ belongs to
the algebra generated by $\hdfgm$ and $\acFEi$:
\ylee{bai5.18}
\zasoth|_{\hdfbt=0}\in\Tfl[\hdfgm,\acFEo,\acFEt,\acFEh].
\yeee
In fact, we conjecture that if $\hdfbt=0$, then the associator is a pure exponential:
\ylee{bai5.19}
\zasoth|_{\hdfbt=0} =
\exp\brb{\tfrac{2}{3}\,\hdfgm^{\inI\inJ\inK}\acFv{\xEo,\inI}\acFv{\xEt,\inJ}\acFv{\xEh,\inK}}.
\yeee

\subsection{A geometric description of the 2-category $\ctLLXsom$}
\label{ss.gdgs}

Following the outline of subsection\rw{ss.outline}, we apply the
results of the previous subsection to formulate conjectures about
a geometric description of the category $\ctLLXsom$, where $\Xsom$
is a general \hlsmm. Our goal is to
 explain how the statements
of subsection\rw{ss.tcthlsmm} referring to the case of $\xX=\TsU$, should be modified for
a general $\Xsom$.

The pairs $\gYL$, where $\yY\subset\xX$ is a lagrangian submanifold
and $\vLY\rightarrow \yY$ is a line bundle such that
$\vLY\ott=\cnKY$, are still objects of the 2-category
$\ctLLXsom$. We conjecture that the analogs of holomorphic
fibration
objects $\goYL$ also appear, but this time $\ycY\rightarrow\yY$ is
not a holomorphic fibration, as in the case of $\xX=\TsU$, but
rather a special `non-holomorphic' deformation of a holomorphic fibration. The
reason for this deformation is similar to the non-holomorphicity
of the functions $\xcW$ which solve the equation\rx{bae2.0a3}, but
we will not explore this subject further.

Suppose that two lagrangian submanifolds $\yYo,\yYt\subset\xX$
have a \gdint. We conjecture that the category of morphisms between them
is the deformed and shifted 2-periodic category of their
intersection:
\ee
\label{bae1.92z} 
\!\!\!\!\!\!\!\!\!\Hom_{\ctLLXsom}\brB{\gYLo,\gYLt} =
\rDsprfYadoct\ytrnLot\btrnv{\shlf\dim \xX - \dim\yYot-1},
\eee
where $\yYoct\edfn\yYocct$, the line bundle $\vcLot\rightarrow\yYoct$ is defined by
\ex{bae1.19b1} adapted to the case of \opfib s:
\ylee{bae1.92z1}
\vcLot \edfn
\vLYo|_{\yYoct}\otimes
\vLYt|_{\yYoct}
\otimes\cnKYot^{-1}
\yeee
and
 $\adfmot\in\xbOmb\bro{\xtbw^{\bullet} \Tng\yYoct}$
is a special \CMlt\ which determines the \tAinf-deformation
of the category $\rDsprfv{\yYoct}$.

Based on the results of the previous subsection, we make
the following conjectures about $\adfmot$:
\begin{enumerate}

\item

The \CMlt\ $\adfmot$ is \xrlb\ and $\deg \adfmot\geq 2$:
%
\ylee{bah2.1a}
\adfmot = \sum_{\xmi=2}^\infty \adfmoti,\qquad
\adfmoti\in\xbOmv{\xmi}(\xcA,\wedge^{\xmi}\yYoct).
\yeee

\item

If at least one of the classes $\tdfbtYo$, $\tdfbtYt$ determined by the exact
sequences\rx{beq2.25b2} is zero
and  the other lagrangian submanifold has a presentation of
subsection\rw{ss.odc} as the
graph of a differential $\del\xcW$, where $\xcW$ is a function on the first
lagrangian surface, then $\xdfmot=0$.

\item

If $\yYo=\yYt=\yY$, then
$\adfmotv{2}=0$
and $\adfmotv{3}$
is given by the formula\rx{beq2.42b}, where $\hdfbt$ represents
$\tdfbtY$ and $\crR$ is the curvature of the tangent bundle $\Tng\yY$.

\item If $\yYo=\yYt=\yY$ and the Atiyah class $\ccrR$ of the
tangent bundle $\TY$ is zero, then $\adfmot=0$.

\end{enumerate}

In order to derive these conjectures from \ex{beq2.18b},
we consider a tubular neighborhood of
$\yYo$ (or $\yYt$) as a tubular neighborhood of the zero section of
a deformed cotangent bundle $\mtdfTYoSk$ with an appropriate
deformation parameter $\hdf$. Then the object $\yYo$ corresponds to
the zero section and hence it is represented by the holomorphic function
$\xcWo = 0$. We assume that within $\mtdfTYoSk$ the second object $\yYt$
is of the form $\yYW$ for an appropriate function $\xcW$ on
$\yYo$. Generally, this is not true, but we expect that conjecture 1 holds
true independently of whether such a presentation
exists, while conjectures 3 and 4 correspond to the case $\xcW=0$.
Finally, we assume that the line bundle $\vLYt$ is the
pull-back of the deformation of $\vLYo$ under the projection of $\yYt$ onto the
zero-section of $\TsY_1$
(recall that the complex structure of the projection of $\yYt$
onto the base $\yYo$ of the cotangent bundle has a complex
structure corresponding to the Beltrami differential\rx{bae2.0b4},
hence the bundle $\vLYo$ has to be deformed in order to be
holomorphic with respect to it). Under these assumptions
\xlee{bah2.1}
\Hom_{\ctLLXsom}\brB{\gYLo,\gYLt} =
\Hom_{\rDDprfYoSk}(0,\xcW) = \rDsprfYoxdoct
\xeee
(\cf \ex{beq2.14b}), where $\xdfmot$ is the deformation parameter
determined by \ex{beq2.18b}, in which we set $\xcWo=0$,
$\xcWt=\xcW$ and, consequently, $\xcWot=\xcW$. Hence $\xdfmot$ has
the expansion\rx{beq2.15b}:
$\xdfmot = \sum_{\xmi=0}^\infty \xdfmoti$,
$\xdfmoti\in\xbOmv{\xmi}(\xcA,\wedge^{\xmi}\TY_1)$. Here $\xdfmotz
= \xcW$, while all other terms $\xdfmoti$ depend on $\xcW$ by
being polynomials in $\del\xcW$ and its covariant holomorphic
differentials $\nabla^k\del\xcW$, $k\geq 1$.

Since we assumed that
$\yYo$ and $\yYt$ have a \gdint, it follows that $\xcW$ has a
\xgd\ critical locus $\CrW$ which is isomorphic to the intersection
$\yYoct$. We conjecture that the category $\rDsprfYoxdoct$
localizes to $\CrW$:
\ylee{bah2.3}
\rDsprfYoxdoct = \rDsprfYadoct\ytrnLot\btrnv{\shlf\dim \xX -
\dim\yYot-1},
\yeee
and the deformation parameter $\adfmot$ is determined somehow by the
restriction $\xdfmot\rCrW$. We do not understand this relation
precisely, but we can still make conjectures about $\adfmot$ based
on the properties of $\xdfmot\rCrW$.

Consider the degree $\dgdW$ defined by \ex{bah2.4}.
%
The critical locus $\CrW$ is determined by
the condition $\del\xcW=0$, hence if $\dgdW\xdfmoti\geq 1$, then
$\xdfmoti\rCrW=0$.
%
Then explicit formula\rx{beq2.23b} for $\xdfmoto=\mu$
implies that $\xdfmo\rCrW=0$. Since $\xcW$ is locally constant at
$\CrW$, we may also assume that $\xdfmotz=0$. Thus our first
conjecture is that $\adfmotz=\adfmoto=0$. We also conjecture that
$\adfmot$ is \xrlb, because the same is true for $\xdfmot$.

The formula\rx{beq2.25b3} states that $\cdfbt = \tdfbtYo$, so if
$\tdfbtYo=0$, then $\cdfbt=0$ and we can use the gauge transformation of the
\CMlt\ $\hdf$
in order to set $\hdfbt=0$. Now
Corollary\rw{crl1} says that $\dgdW\xdfmi\geq 2$ for $i\geq 2$, so
$\xdfmoti\rCrW=0$ for all $i$. Hence we conjecture that if
$\tdfbtYo=0$ and $\yYt$ has a presentation as the graph of $\del\xcW$, then $\adfmot=0$.

If $\hdfbt=0$ and $\yYt$ is presented as the graph of $\del\xcW$
for a function $\xcW$ satisfying the equation\rx{bae2.0a3}, then
the normal bundle $\vccLo\edfn\Tng\yYo|_{\yYoct}/\Tng(\yYoct)$
appearing in \ex{bae1.90a} admits an $\OnC$ structure. Indeed,
$\yYoct=\CrW$, so $\del\xcW|_{\yYoct}=0$ and there is a
well-defined Hessian
$\hsnW\in\Gamma(\rmS^2\dulv{\vccLo})$. This Hessian is
non-degenerate, because we assumed that $\yYo$ and $\yYt$ have a
\gdint. Equation\rx{bae2.0a3} implies that generally it satisfies
the equation
\ylee{bai2.1}
\dlb\hsnWb = \hdfbt\spsmb \brb{\hsnWb\,\hsnWb},
\yeee
but since we assumed that $\hdfbt=0$ we find that the Hessian is
holomorphic: $\dlb\hsnWb=0$. The holomorphic non-degenerate
Hessian provides the $\OnC$ structure for the bundle $\vccLo$.
In fact, we suspect that the converse is also true: if $\hdfbt=0$
and the bundle $\vccLo$ has an $\OnC$ structure then the
lagrangian submanifold $\yYt$ has a presentation as the graph of
$\del\xcW$ at least in a tubular neighborhood of $\yYocct$.

If $\yYo=\yYt=\yY$, then $\yYoct=\yY$ and $\xcW=0$, so
\ex{bah2.1} says that $\adfmot=\xdfmot$. Hence
$\deg \adfmot\geq 3$ and $\adfmotv{3}$
is given by the formula\rx{beq2.42b}. Also, if $\tdfbtYo=0$, then
$\adfmot=0$ follows directly from Corollary\rw{crl1} without any
further conjectures regarding the localization properties of the
deformed category $\rDsprfYoxdoct$.

Finally, if $\yYo=\yYt=\yY$ and the Atiyah class $\ccrR$ of $\TY$
is zero, then Corollary\rw{crl3} says that $\xdfmot=0$. Since in
this case $\adfmot=\xdfmot$, then $\adfmot=0$.

We cannot say much about the deformation of the
composition\rx{bae1.93} except that when all
$\ycYi$ are \opfib s with the same base $\yY=\yYo=\yYt=\yYh$, and the Atiyah class
$\ccrR$ of $\TY$ is zero, then the deformation of
the composition rule\rx{baj3.2} is described by the formulas
of subsection\rw{ss.dmsbfo} in which we replace $\xcA$ with $\yY$.

\section{\Mlc\ definition of the 2-category $\ctLLXsom$}
\label{s.mcl}

\subsection{Symplectic rectangles}

Let $\xcAbx$ denote an $\zn$-dimensional Stein complex manifold $\xcA$ equipped with
holomorphic coordinate functions $\bax=\ax_1,\ldots\ax_n$. The
functions $\bax$ determine an embedding
$\xcAbx\hookrightarrow\ICnbax$, where $\ICnbax$ is the affine
space $\IC^{\zn}$ equipped with the standard coordinates $\bax$.
In other words, $\xcAbx$ is just an open subspace of $\ICnbax$
with inherited coordinates.

A \emph{\smrc} is a product $\UxVy$
with the \hlsms\ determined by the 2-form
$\som = \sum_{i=1}^{\zn} d\ay_i\wedge d\ax_i$.
The identity map establishes an isomorphism between $\UxVy$ and
$\UmxVmy$. The permutation map
$\prms\colon\xcA\times\xcV\rightarrow \xcV\times\xcA$ establishes
the isomorphisms $\UxVy\rightarrow\VmyUx$ and
$\UxVy\rightarrow\VyUmx$.

A \smrc\
$\UxVy$ has a pair of transversal lagrangian fibrations:
a \qfib\ $\au\times\xcVby$ for $\au\in\xcA$ and a \pfib\
$\xcAbx\times\av$ for $\av\in\xcV$.
A \emph{\qemb}
\xlee{bah8.0}
\sfq\colon\UxVypp\xemq\UxVy
\xeee
is a
symplectic embedding such that there exists an embedding
$\rfq\colon\xcAp\hookrightarrow\xcA$ for which the diagram
\ylee{bah8.1}
\xymatrix{
\sUVpp \ar@{^{(}->}[r]^{\sfq}_>{q} \ar[d]& \sUV \ar[d]
\\
\xcAp \ar@{^{(}->}[r]^{\rfq} & \xcA
}
\yeee
is commutative. In other words, a \qemb\ must preserve the \qfib.
A composition of \qemb s 
\xlee{bah8.1b}
\xymatrix@C=1.2cm{
\UxVyh
\ar@{^{(}->}[r]^-{\sfqth}_>{q}
&
\UxVyt
\ar@{^{(}->}[r]^-{\sfqot}_>{q}
&
\UxVyo
}
\xeee
is a \qemb.

The cotangent bundle $\TsUbx$ has a canonical structure of a \smrc,
because
the holomorphic differentials $\del\bax$ form a frame of the
cotangent bundle $\TsUbx$ thus providing an isomorphism
\xlee{bah8.1a}
\TsUbx\xmapta{\cong}\xcAbx\times\ICnbay.
\xeee
Moreover, an embedding
$\xcVby\hookrightarrow\ICnbay$ generates an embedding
$\UxVy\hookrightarrow\TsUbx$, which preserves the symplectic
structure as well as both lagrangian fibrations.

\subsection{2-categories of \smrc s and their functors}
\label{ss.tctf}

\hyphenation{ca-te-go-ry}
We define the 2-cate\-go\-ry $\rDDprfaUxVy$ as a full subcategory of
$\rDDprfaU$. A \xcfb\ $\cmfcUW\in\rDDprfaU$ is an object of
$\rDDprfaUxVy$ if its support\rx{bah8.2} lies within
$\UxVy$ as embedded into $\Ts\xcAbx$:
\ylee{bah8.3}
\yYUW\subset\UxVy\subset\Ts\xcAbx.
\yeee
The isomorphism\rx{bah8.1a} implies the equivalence of categories
\xlee{bah8.3a}
\ctLLbv{\xcAbx\times\ICnbay} \simeq \rDDprfU.
\xeee

A curved fibration $\cmfcUWot\in\rDDprfav{\xcAo\times\xcAt}$
determines the \xxtf\ $\xPhUWot$ of \ex{bae1.69}. The formula\rx{bah7.3} describing
the transformation of the support of a curved fibration under the
action of $\xPhUWot$ implies that if the support of $\cmfcUWot$
fits within the product of \smrc s:
\xlee{bah8.4}
\yYv{\cmfcUWot}\subset
(\UxVyo)\times(\xcAxt\times\xcVyt),
\xeee
then the \xxtf\ $\xPhUWot$ restricts to the \xxtf
\xlee{bah8.5}
\xPhUWot\colon\ctLLbv{\UxVyo}\longrightarrow\ctLLbv{\UxVyt}.
\yeee

A particular example of the \xxtf\rx{bah8.5} is the analog of
\Ldrt s\rx{bag3.1a}. This time the \Ldrt s are the \xxtf s
\xlee{bah8.6}
\dLtp\colon\rDDprfaUxVy \longrightarrow \rDDprfaVyUmx,
\qquad
\dLtm\colon\rDDprfaUxVy\longrightarrow\rDDprfaVmyUx
\xeee
determined through \ex{bae1.69} by the \opfib\ and the curving
$\bax\cdot\bay\edfn  \sum_{i=1}^\zn\ax_i\ay_i$:
\wlee{bah8.7}
\dLt_{\pm} \edfn \dPhv{\pm\bax\cdot\bay}.
\weee
It is easy to see that the curvings $\pm\bax\cdot\bay$ satisfy the
condition\rx{bah8.4} and, moreover, the \Ltf s essentially do not change the
supports of objects: for a \xcfb\ $\cmfcUW\in\rDDprfaUxVy$
\ylee{bah8.8}
\yYv{\dLt_{\pm}\cmfcUW} = \prms\brb{ \yYUW }.
\yeee

We conjecture that the composition of \Ltf s yields the identity
\xxtf:
\wlee{bah8.9}
\dLtp\circ\dLtm \simeq \dLtm\circ\dLtp \simeq \xIdv{\rDDprfaUxVy},
\weee
so the \Ltf s themselves establish equivalences of 2-categories in
\ex{bah8.6}.

An important class of \xxtf s related to 2-categories
$\rDDprfaUxVy$ are restrictions. Suppose that $\xcAp$ is a
submanifold of $\xcA$ of the same dimension. Then there is the restriction
\xxtf\ $\rsf\colon\rDDprfaU\rightarrow\rDDprfaUp$, which acts on
\xcfb s and their morphisms by restricting them from $\xcA$ to
$\xcAp$. The \xxtf\ $\rsf$ can be restricted to the subcategory:
\xlee{bah8.10}
\rsf\colon \rDDprfaUxVy\longrightarrow \rDDprfaUp.
\xeee
If the subset $\xcAp\subset\xcA$ inherits the coordinates $\bax$
then the image of the \xxtf\rx{bah8.10} lies within
$\ctLLbv{\UpxVy}$:
\xlee{bah8.11}
\rsf\colon \rDDprfaUxVy\longrightarrow \rDDprfaUpxVy.
\xeee

For a \qemb\rx{bah8.0} we define the restriction \xxtf
\ylee{bah8.12}
\rsfe\colon\rDDprfaUxVy\longrightarrow\rDDprfaUxVypp
\yeee
as the composition of five \xxtf s:
%
\ylee{bah8.13}
\xymatrix@C=1.5cm@R=1.5cm{
%
\rDDprfaUxVy
\ar[r]^-{\rsfo}
\ar@/^3pc/[rrr]^-{\rsfe}
&
\rDDprfaUp
\ar[r]^-{\xId}
&
\ctLLbv{\xcApbxp\times\IC_{\bayp}}
\ar[d]^-{\dLtp}
&
\rDDprfaUxVypp
\\
& &
\ctLLbv{\IC_{\bayp}\times\xcApmbxp}
\ar[r]^-{\rsft}
&
\ctLLbv{\xcVpbyp\times\xcApmbxp}
\ar[u]^-{\dLtm}
}
\yeee
In this diagram the restriction \xxtf\ $\rsfo$ is of the
type\rx{bah8.10}, the restriction \xxtf\ $\rsft$ is of the
type\rx{bah8.11} and the equivalence $\xId$ is of the
type\rx{bah8.3a}.

We conjecture that the restriction \xxtf\ of the composition of
\qemb s is isomorphic to the composition of individual restrictions, that is,
for a chain of \qemb s\rx{bah8.1b}
\xlee{bah8.13}
\rsfv{\sfqot\circ\sfqth}\simeq\rsfv{\sfqot}\circ\rsfv{\sfqth}.
\xeee


\subsection{A \prshf\ of 2-categories}

A \emph{\rcch} in a \hlsmm\ $\Xsom$ is a symplectic map
\xlee{bah8.14}
\mch\colon\UxVy\rightarrow\xX
\xeee
To every \rcch\ we associate the
2-category $\rDDprfaUxVy$. The relation\rx{bah8.13} suggests that these
chart 2-categories form a \prshf\ $\prsfXsom$. An object $\zzO$ of the category
$\ctLLXsom$ is defined to be a global section of this \prshf: to
every \rcch\rx{bah8.14} we associate an object
$\zzOf\in\rDDprfaUxVy$ with two conditions:
for any commutative triangle
\ylee{bah8.15}
\xymatrix{
\UxVy \ar[rr]^-{\prms}
\ar[rd]_-{\mch}
&&
\VyUmx
\ar[ld]^-{\mch\p}
\\
&\xX
}
\yeee
there is a relation $\zzOfp \cong \dLtp \zzOf$, and for any
commutative triangle
\ylee{bah8.16}
\xymatrix{
\UxVypp
\ar[rd]_-{\mch\p}
\ar@{^{(}->}[rr]^-{\sfq}_>{q}
&&
\UxVy
\ar[ld]^-{\mch}
\\
&\xX
}
\yeee
there should be a relation $\zzOfp \cong \rsfe\zzOf$.

Two global sections $\zzOo,\zzOt\in\ctLLXsom$ determine a \prshf\
of categories $\mathcal{H}om(\zzOo,\zzOt)$: to every
\rcch\rx{bah8.14} we associate the category
$\Hom(\zzOfo,\zzOft)$ and we define
$\Hom_{\ctLLXsom}(\zzOo,\zzOt)$ as the category of global sections
of this \prshf.

\section{Categorified algebraic geometry and the RW model}\label{sec:sheaves}

\subsection{RW model of a graded cotangent bundle}

In the case when $X=\Ts Y$ one can promote the RW model from a
$\ZZ_2$-graded \TQFT\ to a $\ZZ$-graded one, as explained in \cite{KRS1}.
To this end, one assigns cohomological degree $2$ to linear coordinates on the
fiber of the projection $\Ts Y\raa Y$. From the physical viewpoint,
the degree is the weight with respect to a $U(1)$ ghost number
symmetry. We will call the resulting graded manifold $\Ts Y[2]$. The
sheaf of holomorphic functions on $\Ts Y[2]$ is a quasicoherent
sheaf of graded algebras on $Y$:
$$
\cO_X=\oplus_p \Sym^p {\rm T} Y.
$$
The RW model with the target $\Ts Y[2]$ has $U(1)$ ghost number
symmetry, and it is natural to consider boundary conditions and
topological defects which preserve this symmetry. This gives a
$\ZZ$-graded version of the model.

The 2-category of boundary conditions supported on $Y$ has a
distinguished object: the zero section of $\Ts Y[2]$. It is easy to
see that this boundary condition is invariant with respect to the
$U(1)$ ghost number symmetry. The corresponding endomorphism
category is $\sDb(Y)$, the bounded derived category of coherent
sheaves on $Y$. From the physical viewpoint, it arises as the
homotopy category of a DG-category $\Dperf(Y)$. Objects of
$\Dperf(Y)$ are perfect DG-modules over the $\ZZ$-graded Dolbeault
DG-algebra $(\obul(Y),\bpartial)$, with morphisms being the usual
morphisms of DG-modules. $\sDb(Y)$ is a symmetric monoidal
DG-category; as discussed in \cite{KRS1}, the monoidal structure is
the standard one (this is easy to see on the classical level, but it
takes some work to show that there are no quantum corrections). The
algebra of boundary local operators for the distinguished boundary
condition (i.e. the endomorphism algebra of the unit object in the
endomorphism category) is isomorphic to $H^*(\cO_Y)$.

In the $\ZZ$-graded case, infinitesimal deformations of a boundary
condition correspond to degree-2 elements in the algebra of local
boundary operators. Thus infinitesimal deformations of the
distinguished boundary condition are parameterized by $H^2(\cO_Y)$.
If $Y$ is compact and K\"ahler, such deformations are unobstructed.
Indeed, we can choose a harmonic representative $B$ of a class in
$H^2(\cO_Y)$, and then the deformation of the boundary action is
simply
$$
\int_{\partial M} \phi^*B,
$$
where $\phi$ is a map from the space-time $M$ to the target $X$.
Since the form $B$ is closed, such deformation is obviously
BRST-invariant and does not affect BRST-transformations of any
fields. We will call such a deformation a B-field deformation, by
analogy with the 2d sigma-models.

Let $(Y,B)$ denote the distinguished boundary condition deformed by
$B$. The category of morphisms from $(Y,B_1)$ to $(Y,B_2)$ is the
bounded derived category of twisted coherent sheaves on $Y$, where
the twist is given by the class of $B_2-B_1$. We will denote this
category $\sD(Y,B_2-B_1)$. The composition of morphisms is the
obvious one (tensor product of twisted coherent sheaves).
Physically, $\sDb(Y,B)$ arises as the homotopy category of a certain
DG-category which we denote $\Dperf(Y,B)$. It is the category of
perfect CDG-modules over the CDGA $(\obul(Y),\bpartial,B)$.

More complicated boundary conditions can be obtained by considering
complex fibrations $\Z\raa Y$ equipped with a B-field $B\in
H^2(\cO_\Z)$. The category of morphisms from $(\Z_1,B_1)$ to
$(\Z_2,B_2)$ is the bounded derived category of twisted coherent
sheaves on $\Z_1\times_Y \Z_2$ with the twist given by
$\pi_2^*B_2-\pi_1^* B_1$, where $\pi_s$ is the projection from
$\Z_1\times_Y \Z_2$ to $\Z_s$, $s=1,2$.

To understand the resulting 2-category better, note that an object
of $\sDb(\Z_1\times_Y\Z_2,B_2-B_1)$ defines a functor from
$\sDb(\Z_1,B_1)$ to $\sDb(\Z_2,B_2)$. Composition of morphisms in
the 2-category of boundary conditions is simply the composition of
functors. Moreover, this functor intertwines the natural action of
$\sDb(Y)$, regarded as a monoidal category, on $\sDb(\Z_1,B_1)$ and
$\sDb(\Z_2,B_2)$. That is, if we regard the categories
$\sDb(\Z_s,B_s),$ $s=1,2$ as modules over the monoidal category
$\sDb(Y)$, then this functor defines a morphism in the 2-category of
modules.

Note that for a $\CC$-linear (or DG) monoidal category $\cC$ there
are two very different notions of a module: a module over $\cC$
regarded simply as a $\CC$-linear (or DG) category, and a module
over $\cC$ regarded as a monoidal $\CC$-linear (or monoidal DG)
category. The former is a functor from $\cC$ to the category of
complex vector spaces $\Vect_\CC$ (or the category of differential
graded complex vector spaces); the latter is a $\CC$-linear (or DG)
category which is acted upon by $\cC$. To avoid confusion, we will
call the latter notion a 2-module over $\cC$. This terminology is
not standard,\footnote{The more standard name for a 2-module is a
module category.} but natural, if we think about a monoidal category
as a 2-algebra, i.e. a categorification of an algebra. 2-modules
over a monoidal category $\cC$ form a 2-category.

One could hope that any morphism in the 2-category of 2-modules is
represented by an object of $\sDb(\Z_1\times_Y\Z_2,B_2-B_1)$. Then
the 2-category of boundary conditions in the RW model would be a
full sub-2-category of the 2-category of 2-modules. This statement
is incorrect as formulated, however, apparently it does become
correct if we replace the derived category of (twisted) coherent
sheaves with its enhancement. Recall that an enhancement of a
triangulated category $\cC$ is a DG-category $\frC$ whose homotopy
category $H^0(\frC)$ is triangulated and an equivalence of
$H^0(\frC)$ and $\cC$. From the physical viewpoint, a natural
enhancement of $\sDb(Y)$ is the DG-category $\Dperf(Y)$. Similarly,
a natural enhancement of $\sDb(Y,B)$ is the DG-category
$\Dperf(Y,B)$ of perfect CDG-modules over the CDGA
$(\obul(Y),\bpartial,B)$. The category $\Dperf(Y)$ is a monoidal
DG-category which acts by DG-functors on the DG-category
$\Dperf(\Z,B)$. Any object of $\Dperf (\Z_1\times_Y\Z_2,B_2-B_1)$
determines a DG-functor from $\Dperf (\Z_1,B_1)$ to
$\Dperf(\Z_2,B_2)$ which intertwines the action of $\Dperf(Y)$. The
improved version of the conjecture is that any such DG-functor is
represented by an object of $\Dperf (Y)$.

In \cite{Toen} this conjecture was proved for $Y$ being a point. In
\cite{BFN} the proof was extended to the case when $Y$ is a more
general scheme.

The conclusion is that the 2-category of boundary conditions in the
$\ZZ$-graded version of the RW model with target $\Ts Y[2]$ is the
homotopy category of a full sub-2-category in the 2-category of
2-modules over the monoidal DG-category $\Dperf (Y)$.

\subsection{Derived categorical sheaves}

Complex fibrations over $Y$ play a role in the RW model similar to
that played by holomorphic vector bundles in the B-model with target
$Y$. But it is well-known that more general coherent sheaves also
arise as B-branes, and it is natural to ask if boundary conditions
in the RW model can be similarly generalized.

It is convenient to take a more algebraic viewpoint and replace
complex fibrations over $Y$ with families of algebras or DG-algebras
over $Y$. Likely this entails no essential loss of generality. For
example, it is known that for any sufficiently nice (quasi-compact
and quasi-separated) scheme $Z$ the derived category of complexes of
sheaves on $Z$ with quasicoherent cohomology is equivalent to the
derived category of modules over some DG-algebra with bounded
cohomology \cite{BvB}. Thus we will replace the fibration $\Z$ with
a sheaf of DG-algebras over $Y$. One may conjecture that any sheaf
of DG-algebras over $Y$ can be interpreted as a boundary condition
in the RW model.

To test this conjecture, we need to have a reasonable definition of
the category of morphisms between sheaves of DG-algebras. A natural
definition has been sketched by B.~Toen and G.~Vezzosi \cite{ToVe}. They
work with more general objects called derived categorical sheaves
over $Y$. A derived categorical sheaf is a sheaf of DG-categories
over $Y$. This means that to any affine open subscheme $\Spec\,
A=U\subset Y$ one attaches a DG-category $\frC(U)$ over $A$, to any
inclusion of affine open subschemes $U'\subset U$ one attaches a
morphism of DG-categories $r_{U'U}:\frC(U)\raa\frC(U')$, and to any
inclusion of affine open subschemes $U''\subset U'\subset U$ one
attaches an invertible 2-morphism from $r_{U''U'}\circ r_{U'U}$ to
$r_{U''U}$. These data must satisfy a number of conditions which are
spelled out, for example, in \cite{catsheaf}. A sheaf of DG-algebras
can be thought of as a special case of this, with the DG-category
$\frC(U)$ having a single object for any $U$.

The category of morphisms from the derived categorical sheaf $\sT_1$
to the derived categorical sheaf $\sT_2$ is defined as follows.
First of all, one can define the derived tensor product of two
derived categorical sheaves which is again a sheaf of DG-categories. In \cite{ToVe} it is denoted $T_1\otimes^{\mathbb L} T_2$.
The category of morphisms from $\sT_1$ to $\sT_2$ is defined to be
the derived category of modules over the sheaf of DG-categories $\sT_1^{op}\otimes^{\mathbb L} T_2$, where $T_1^{op}$ is the opposite of $T_1$. In
this way one gets a 2-category of derived categorical sheaves over
$Y$. There are versions of this definition which depend on which
modules precisely one considers.

The simplest object in this 2-category is the structure sheaf
$\cO_Y$ regarded as a sheaf of DG-algebras with zero differential.
It corresponds to the distinguished boundary condition in the
2-category of boundary conditions for the RW model with target $\Ts
Y[2]$. Its endomorphism category is $\sDb(Y)$; this agrees with the
endomorphism category of the distinguished boundary condition in the
RW model. Given any other derived categorical sheaf $\sT$ over $Y$,
the category of morphisms from the distinguished object to $\sT$ is
a 2-module over the monoidal category $\sDb(Y)$. Thus the 2-category
of derived categorical sheaves over $Y$ is embedded into the
2-category of 2-modules over $\sDb(Y)$.

A simple but interesting example of a derived categorical sheaf is a
skyscraper sheaf, i.e. a sheaf of DG-categories such that $\frC(U)$
is quasi-equivalent to the trivial category if $U$ does not contain
a point $p\in Y$, and is quasi-equivalent to a fixed DG-category
$\frC_0$ otherwise. We may call this a skyscraper sheaf with stalk
$\frC_0$. We now explain how to construct the corresponding boundary
condition in the RW model by allowing the fibration $\Z$ over $Y$ to
carry a nontrivial curving $W\in H^0(\cO_\Z)$. Such boundary
conditions should be regarded as 3d analogues of 0-branes. For
simplicity we will assume that the DG-category $\frC_0$ is simply a
DG-algebra $\wA$, and is moreover of a geometric origin, i.e. its
derived category of modules $\sD(\wA)$ is equivalent to the derived
category of coherent sheaves on some complex manifold $V$.

First we note that in order for a curving $W$ to preserve the ghost
number symmetry, we have to allow the fiber $\Z$ to be a graded
manifold with a nontrivial $\CC^*$ action. Then the space
$H^0(\cO_Z)$ is also graded, and $W$ must sit in its degree-2
component. We will call such $W$ a superpotential. A graded
fibration $\Z\raa Y$ equipped with a superpotential $W$ of degree 2
defines a boundary condition for the RW model.

The category of morphisms from the distinguished boundary condition
to the boundary condition $(\Z,W)$ is $H^\xbulx(\Dperf
(\Z,W))=\sDb(\Z,W)$. Note that this category is equivalent to a
trivial one if $W$ has no critical points \cite{KLi,Orlov:MF}. This
is a local statement: given an open set $U\subset Y$ we may consider
the restriction $(\Z_U,W_U)$ of $(\Z,W)$ to $U$ and the category
$\sDb(\Z_U,W_U)$; this category is trivial if $\Z_U$ does not
contain critical points of $W$. Therefore a natural candidate for an
analogue of a skyscraper sheaf is a pair $(\Z,W)$ such that all
critical points of $W$ are contained in the fiber over a point $p\in
Y$.

To be concrete, let us consider the case when $Y$ is the
$n$-dimensional affine space $\AA_n$ with coordinates
$y^1,\ldots,y^n$. We will describe a boundary condition in the RW
model with target $\Ts Y[2]$ which corresponds to a skyscraper sheaf
over $Y$ with the stalk at $y=0$ being a DG-algebra $\wA$ of a
geometric origin. Let $\Z=\AA_n[2]\times Y\times V$, where
$\AA_n[2]$ denotes the affine space with linear coordinates
$a_1,\ldots,a_n$ of cohomological degree 2. The graded manifold $\Z$ is a trivial fibration over $Y$.
The superpotential will be
$$
W=\sum_i y^i a_i,
$$
We may think of $y^i$ and $a_i$ as coordinates on $\Ts Y[2]$.

The category of morphisms from the distinguished boundary condition
to $(\Z,W)$ is
\\
$\sDb(\Z,W)$. By Kn\"orrer periodicity, it is
equivalent to $\sDb(V)$. Furthermore, $W$ has a single critical
point $y=a=0$, so the category $\sDb(\Z_U,W_U)$ is equivalent to a
trivial one if $U$ does not contain the point $y=0$.

We propose that this boundary condition corresponds to the
skyscraper sheaf with the stalk $\wA$ at $y=0$. By definition, the
category of morphisms from the derived categorical sheaf $\cO_Y$ to
this skyscraper is the category $\sDb(\wA)\simeq\sDb(V)$, which
agrees with the RW model.

To check this proposal further, let us compare the endomorphism
category of $(\Z,W)$ regarded as a boundary condition in the RW
model and the endomorphism category of the skyscraper sheaf. The
former is the category $\sDb(\AA_n[2]\times \AA_n[2]\times V\times
V\times Y,\tilde W)$. The superpotential $\tilde W$ is given by
$$
\tilde W= y^i(a_i-\tilde a_i),
$$
where $\tilde a_i$ denote the coordinates on the second copy of
$\AA_n[2]$. By Kn\"orrer periodicity, this category is equivalent to
$\sDb(\AA_n[2]\times V\times V)$.

To compute the endomorphism category of a skyscraper sheaf, we first
need to compute its derived tensor product with itself. Since the
base is an affine space $\AA_n$, we can think about sheaves of
DG-algebras in algebraic terms, i.e. as DG-algebras over the ring
$\CC[y^1,\ldots,y^n]$. From this point of view, the skyscraper sheaf
with a stalk $\wA$ is simply the DG-algebra $\wA$ made into a
$\CC[y^1,\ldots,y^n]$-module by letting all $y^i$ act trivially.
Equivalently, it is a tensor product over $\CC$ of the DG-algebra
$\wA$ over $\CC$ and $\CC$ regarded as a DG-algebra over
$\CC[y^1,\ldots,y^n]$ with a trivial action of $y^i$ for all $i$ and
a trivial differential.

Since such a module is not flat over $\CC[y^1,\ldots,y^n]$, to
compute its derived tensor product with itself we need a flat
resolution for it.\footnote{We are grateful to Dima Orlov for
explaining to us the content of this paragraph.} Consider a
DG-algebra
$$
\KK_n=\left(\CC[y^1,\ldots,y^n\vert
\theta^1,\ldots,\theta^n],Q\right),
$$
where $\theta^1,\ldots,\theta^n$ are anticommuting odd variables of degree $-1$,
and the differential $Q$ is the Koszul differential
$$
Q=y^i \frac{\partial}{\partial \theta^i}
$$
It is quasi-isomorphic to $\CC$ regarded as a DG-algebra over
$\CC[y^1,\ldots,y^n]$. Hence we can obtain the desired flat
resolution by tensoring over $\CC$ the DG-algebra $\wA$ with
$\KK_n$. The derived tensor product is now computed by tensoring
with $\wA$ over $\CC[y^1,\ldots,y^n]$. The result is a DG-algebra
\begin{equation}\label{wA}
\wA_\theta=\wA\otimes_\CC\wA\otimes_\CC
\CC[\theta^1,\ldots,\theta^n],
\end{equation}
with a trivial action of the variables $y^i$. By definition, the
endomorphism category of the skyscraper is a suitable version of the
derived category of modules over this DG-algebra.

The algebra $\CC[\theta^1,\ldots,\theta^n]$ is Koszul-dual to the
algebra
\begin{equation}\label{CCalg}
\CC[a_1,\ldots,a_n],
\end{equation}
where the variables $a_i$ are even and have degree $2$, and the
differential is zero. Consequently, suitably defined derived
categories of the DG-algebra (\ref{wA}) and the DG-algebra
\begin{equation}\label{wwA}
\wA_a=\wA\otimes_\CC\wA\otimes_\CC \CC[a_1,\ldots,a_n].
\end{equation}
are equivalent. This agrees with what we got from the RW model and
Kn\"orrer periodicity.

Note that the resolution of the skyscraper categorical sheaf used
above is in some sense Koszul-dual to the trivial fibration $\Z=\AA_n[2]\times Y\times V$;
the role of the Koszul differential is played by the superpotential
$W$.

Let us consider a slightly more complicated example: a sheaf of
algebras over $\AA_1=\Spec\, \CC[y]$ which in algebraic terms is the
algebra $\CC[y]/y^k$ over the ring $\CC[y]$. For $k=1$, this is a
special case of the previous example (with $n=1$). We will argue
that there exists a boundary condition in the RW model equivalent to
such a sheaf of algebras.

Note first that the above sheaf of DG-algebras can be deformed into
a collection of $k$ skyscrapers by replacing $y^k$ with $P_k(y)$,
where $P_k$ is a degree-$k$ polynomial without multiple roots. This
corresponds to the following boundary condition in the RW model:
$\ZZ=\AA_1\times \CC[2]$, $W=a P_k(y)$. In the limit when $P_k(y)$
degenerates to $y^k$, we get $W=a y^k$. Therefore we propose that
the boundary condition with $\ZZ=\AA_1\times \CC[2]$, $W=a y^k,$
corresponds to the DG-algebra $\CC[y]/y^k$ over $\CC[y]$.

The category of morphisms from the distinguished boundary condition
to this one is the category of $\CC^*$-equivariant matrix
factorizations of $W=a y^k$. If the proposal is correct, then this
category must be equivalent to the derived category of DG-modules
over $\CC[y]/y^k$. The equivalence presumably arises from the
following matrix factorization:
\begin{align}
D=\begin{pmatrix} 0 & a \\ y^k & 0\end{pmatrix}
\end{align}
Its endomorphism algebra is a DG-algebra quasi-isomorphic to
$\CC[y]/y^k$ with the zero differential. Thus we get a bimodule
which defines a functor from the category of equivariant matrix
factorizations to the derived category of DG-modules over
$\CC[y]/y^k$. With suitable definitions, this functor should be an
equivalence of categories \cite{Orlov:priv}.

We note in passing that this construction allows one to think about
the derived category of modules over $\CC[y]/y^k$ as a category of
B-branes in some physical theory (namely, the Landau-Ginzburg model
on $\CC\times\CC[2]$ with the superpotential $W=a y^k$). In other
words, the Landau-Ginzburg model whose target is a graded manifold allows one to give meaning to such
a singular-looking theory as a sigma-model with target
$\Spec(\CC[y]/y^k)$. Similarly many other graded Landau-Ginzburg models can
be thought of as representing sigma-models whose targets are
singular schemes or even DG-schemes.

\end{document}